\documentclass{emsprocart}

\usepackage[all]{xy}

\usepackage{makeidx}

\usepackage{tikz}
\usetikzlibrary{decorations.pathreplacing}
\usetikzlibrary{matrix,arrows}

\renewcommand{\footnote}[1]{}

\contact[keller@math.jussieu.fr]{Bernhard Keller, Universit\'e Denis Diderot -- Paris~7,
Institut universitaire de France, U.F.R.~de math\'ematiques, Institut de math\'ematiques de Jussieu, U.M.R.~7586
du CNRS, Case 7012, B\^atiment Chevaleret, 75205 Paris, France}

\newtheorem{theorem}{Theorem}[section]
\newtheorem{corollary}[theorem]{Corollary}
\newtheorem{lemma}[theorem]{Lemma}
\newtheorem{proposition}[theorem]{Proposition}
\newtheorem{conjecture}[theorem]{Conjecture}
\newtheorem{main-conjecture}[theorem]{Main Conjecture}

\theoremstyle{definition}

\newtheorem{remark}[theorem]{Remark}

\newcommand{\ie}{{\em i.e.}\ }

\newcommand{\ko}{\: , \;}

\newcommand{\opname}[1]{\operatorname{\mathsf{#1}}}

\newcommand{\nil}{\opname{nil}\nolimits}

\newcommand{\aff}{\hat{\mathbb{A}}}

\newcommand{\per}{\opname{per}\nolimits}
\newcommand{\add}{\mbox{{\rm add}}}

\newcommand{\Gr}{\mbox{{\rm Gr}}}

\newcommand{\tw}{\opname{tw}\nolimits}
\newcommand{\pr}{\mbox{{\rm pr}}}

\newcommand{\Z}{\mathbb{Z}}
\newcommand{\N}{\mathbb{N}}
\newcommand{\Q}{\mathbb{Q}}
\newcommand{\C}{\mathbb{C}}

\renewcommand{\P}{\mathbb{P}}
\newcommand{\E}{\mathbb{E}}

\newcommand{\iso}{\stackrel{_\sim}{\rightarrow}}

\newcommand{\Hom}{\mbox{{\rm Hom}}}

\newcommand{\RHom}{\mbox{{\rm RHom}}}

\newcommand{\Ext}{\mbox{{\rm Ext}}}

\newcommand{\ten}{\otimes}

\newcommand{\ca}{{\mathcal A}}

\newcommand{\cc}{{\mathcal C}}
\newcommand{\cd}{{\mathcal D}}

\newcommand{\cf}{{\mathcal F}}

\newcommand{\cl}{{\mathcal L}}

\newcommand{\cp}{{\mathcal P}}

\newcommand{\ct}{{\mathcal T}}
\newcommand{\cu}{{\mathcal U}}
\newcommand{\cv}{{\mathcal V}}

\newcommand{\eps}{\varepsilon}
\renewcommand{\phi}{\varphi}
\newcommand{\del}{\partial}

\renewcommand{\hat}[1]{\widehat{#1}}

\newcommand{\sgn}{\mbox{sgn}}

\renewcommand{\tilde}[1]{\widetilde{#1}}

\newcommand{\Ad}{\mbox{{\rm Ad}}}

\renewcommand{\sf}{\opname{Sf}}

\title{Cluster algebras and derived categories}

\author[Bernhard Keller]{Bernhard Keller
}

\makeindex

\begin{document}

\begin{abstract}
This is an introductory survey on cluster algebras and their (additive)
categorification using derived categories of Ginzburg algebras. After
a gentle introduction to cluster combinatorics, we
review important examples of coordinate rings admitting a cluster
algebra structure. We then present the general definition of
a cluster algebra and describe the interplay between cluster
variables, coefficients, $c$-vectors and $g$-vectors. We
show how $c$-vectors appear in the study of quantum
cluster algebras and their links to the quantum dilogarithm.
We then present the framework of additive categorification
of cluster algebras based on the notion of quiver with potential and
on the derived category of the associated Ginzburg algebra.
We show how the combinatorics introduced previously
lift to the categorical level and how this leads to proofs,
for cluster algebras associated with quivers, 
of some of Fomin--Zelevinsky's fundamental conjectures.
\end{abstract}



\maketitle

\tableofcontents

\section{Introduction}
\label{s:introduction}

Cluster algebras, invented  \cite{FominZelevinsky02} by Sergey Fomin
and Andrei Zelevinsky around the year 2000, are commutative algebras
whose generators and relations are constructed in a recursive manner.
Among these algebras, there are the algebras of homogeneous
coordinates on the Grassmannians, on the flag varieties and on many
other varieties which play an important role in geometry and representation
theory. Fomin and Zelevinsky's main aim was to set up a combinatorial
framework for the study of the so-called canonical
bases which these algebras possess \cite{Kashiwara90} \cite{Lusztig90}
and which are closely related to the notion of total positivity
\cite{Lusztig96} \cite{Fomin10} in the associated varieties. It has rapidly
turned out that the combinatorics of cluster algebras also
appear in many other subjects, for example in
\begin{itemize}
\item Poisson geometry
\cite{GekhtmanShapiroVainshtein03}
\cite{GekhtmanShapiroVainshtein05}
\cite{GekhtmanShapiroVainshtein08}
\cite{GekhtmanShapiroVainstein10}
\cite{BerensteinZelevinsky05} \ldots ;

\item discrete dynamical systems
\cite{DiFrancescoKedem09}
\cite{FominZelevinsky03b}
\cite{InoueIyamaKunibaNakanishiSuzuki10}
\cite{Kedem08}
\cite{Keller10b}
\cite{Keller10a}
\cite{KunibaNakanishiSuzuki11}
\ldots ;

\item higher Teichm\"uller spaces
% \cite{FockGoncharov03}
% \cite{FockGoncharov07b}
\cite{FockGoncharov06a}
\cite{FockGoncharov07}
\cite{FockGoncharov09}
\cite{FockGoncharov10}
\cite{FockGoncharov09a}
\ldots ;

\item combinatorics and in particular the study of polyhedra
like the Stasheff associahedra
\cite{Chapoton04}
\cite{ChapotonFominZelevinsky02}
\cite{FominReading05}
\cite{FominShapiroThurston08}
\cite{IngallsThomas09}
\cite{Krattenthaler06}
\cite{Musiker11}
\cite{MusikerSchifflerWilliams11a}
\cite{MusikerSchifflerWilliams11}
\cite{Schiffler10}
\ldots ;

\item commutative and non commutative algebraic geometry and
in particular the study of stability conditions in the sense of Bridgeland
\cite{Bridgeland07},
Calabi-Yau algebras
\cite{Ginzburg06}
\cite{IyamaReiten08} ,
Donaldson-Thomas invariants in geometry
% \cite{Kontsevich07a} \cite{Kontsevich07}
\cite{JoyceSong09}
\cite{KontsevichSoibelman08}
\cite{KontsevichSoibelman10}
\cite{Reineke11}
\cite{Szendroi08}
\ldots\
and in string theory
\cite{CecottiEtAl11a}
\cite{CecottiEtAl11}
\cite{CecottiCordovaVafa11}
\cite{CecottiNeitzkeVafa10}
\cite{CecottiVafa11}
\cite{GaiottoMooreNeitzke10}
\cite{GaiottoMooreNeitzke09}
\cite{GaiottoMooreNeitzke10a}
\ldots ;
\item in the representation theory of quivers and finite-dimensional
algebras, cf. for example the survey articles
\cite{Amiot11}
\cite{BuanMarsh06}
\cite{GeissLeclercSchroeer08a}
\cite{Keller10b}
\cite{Leclerc10}
\cite{Reiten10a} \cite{Reiten10}
\cite{Ringel07}
\ldots
\end{itemize}
as well as in mirror symmetry \cite{GrossHackingKeel11},
KP solitons \cite{KodamaWilliams11}, hyperbolic $3$-manifolds
\cite{NagaoTerashimaYamazaki11}, \ldots\ .
We refer to the introductory articles
\cite{Fomin10}
\cite{FominZelevinsky03a}
\cite{Zelevinsky04}
\cite{Zelevinsky02}
\cite{Zelevinsky05}
\cite{Zelevinsky07a}
and to the cluster algebras portal \cite{Fomin07}
for more information on cluster algebras and their
links with other subjects in mathematics (and physics).

In these notes, we give a concise introduction to cluster
algebras and survey their (additive) categorification via derived categories
of Ginzburg dg (=differential graded) algebras. We prepare the
ground for the formal definition of cluster algebras by
giving an approximate description and the first examples
in section~\ref{s:description-examples}.
In section~\ref{s:cluster-algebras-associated-with-quivers},
we introduce the central construction of quiver mutation
and define the cluster algebra associated with a quiver
and, more generally, with a valued quiver (section~\ref{ss:valued-quivers}).
We extend the definition to that of cluster algebras
of geometric type and present several examples
in section~\ref{s:cluster-algebras-of-geometric-type}. Here
we also review results on ring-theoretic properties of
cluster algebras (finite generation and factoriality).
In section~\ref{s:general-cluster-algebras}, we give the general definition
of cluster algebras with coefficients in an arbitrary
semifield. In this general framework, the symmetry between cluster
variables and coefficients becomes apparent, for example in
the separation formulas in Theorem~\ref{thm:separation-formula} but
also, at the `tropical level', in the duality Theorem~\ref{thm:tropical-Langlands-duality}.
In section~\ref{s:quantum-cluster-algebras-and-quantum-dilogarithms},
we present the construction of quantum cluster algebras and its
link with the quantum dilogarithm function. We show how
cluster algebras allow one to construct identities between products
of quantum dilogarithm series. This establishes the link
to Donaldson--Thomas theory, as we will see for example in
section~\ref{ss:proof-quantum-dilog}.

In section~\ref{s:categorification}, we turn to
the (additive) categorification of cluster algebras. In section~5
of \cite{Keller10b}, the reader will find a gentle introduction to
this subject along the lines of the historical development.
We will not repeat this here but restrict ourselves to
a description of the most recent framework, which applies to
arbitrary symmetric cluster algebras (of geometric type).
The basic idea is to lift the cluster variables in the
cluster algebra associated with a quiver $Q$
to suitable representations of $Q$. These representations
have to satisfy certain relations, which are encoded in
a potential on the quiver. We review quivers with potentials
and their mutations following Derksen-Weyman-Zelevinsky
\cite{DerksenWeymanZelevinsky08} in section~\ref{ss:quivers-with-potential}.
A conceptual framework for the study of the representations
of a quiver with potential is provided by the derived category
of the associated Ginzburg dg algebra (section~\ref{ss:der-cat-Ginzburg-alg}).
Here mutations of quivers with potential yield
equivalences between derived categories of Ginzburg
algebras (section~\ref{ss:der-eq}). In fact, each mutation
canonically lifts to {\em two} equivalences. Thus, in trying to
compose the categorical lifts of $N$ mutations, we are forced to choose between
$2^N$ possibilities. The canonical choice was discovered
by Nagao \cite{Nagao10} and is presented in section~\ref{ss:patterns}.
The framework thus created allows for the categorification
of all the data associated with a commutative cluster algebra
(Theorem~\ref{thm:decategorification}). A recent extension
to quantum cluster algebras (under suitable technical 
assumptions\protect{\footnote{cf. the first two lines on page~33 of \cite{Efimov11},
where the potential $W$ is assumed to be well-mutable.}})
is due to Efimov \cite{Efimov11}. Surprisingly, the combinatorial
data determine the categorical data to a very large extent
(sections~\ref{ss:reign-of-the-tropics} and \ref{ss:rigid-objects-and-cluster-monomials}).
We end by linking our formulation of the `decategorification Theorem'~\ref{thm:decategorification}
to the statements available in the literature
(sections~\ref{ss:proof-of-decategorification}) and by proving
Theorem~\ref{thm:quantum-dilog-identity} on quantum dilogarithm identities
(section~\ref{ss:proof-quantum-dilog}).

This introductory survey leaves out a number of
important recent developments, notably monoidal
categorification, as developed by Hernandez-Leclerc
\cite{HernandezLeclerc10} \cite{Leclerc10} and Nakajima \cite{Nakajima11},
the theory of cluster algebras associated with marked surfaces
\cite{FominShapiroThurston08} \cite{MusikerSchifflerWilliams11a}
\cite{MusikerSchifflerWilliams11} \cite{CerulliLabardini11} 
\ldots\ and recent progress on the links between (quantum) cluster
algebras and Lie theory  \cite{GeissLeclercSchroeer12}
\cite{GeissLeclercSchroeer11b} \cite{GeissLeclercSchroeer11} \ldots\ .

\section*{Acknowledgment}
This survey expands on a talk given by 
the author at the GCOE Conference `Derived Categories 2011 Tokyo',
where Theorem~\ref{thm:quantum-dilog-identity} was presented.
He deeply thanks the organizers 
for their invitation and their kind hospitality. He is very grateful
to B.~Leclerc for correcting an error in a previous version
of example~\ref{ss:Gr-3-6} and to H.~Nakajima for a helpful
conversation on the results of \cite{Efimov11}.

\section{Description and first examples of cluster algebras}
\label{s:description-examples}
\index{cluster algebra!examples}
\subsection{Description}\index{cluster algebra!description} \label{ss:description}
A {\em cluster algebra}\index{cluster algebra} is a commutative $\Q$-algebra endowed
with a set of distinguished generators (the {\em cluster variables})
grouped into overlapping subsets (the {\em clusters}) of constant
cardinality (the {\em rank})\index{cluster algebra!rank of} which are constructed recursively
via {\em mutation}\index{mutation!in cluster theory} from an {\em initial cluster}.
The set of cluster variables can be finite or infinite.

\begin{theorem}[\cite{FominZelevinsky03}] The cluster algebras having
only a finite number of cluster variables are parametrized by
the finite root systems\index{root system}.
\end{theorem}

Thus, the classification is analogous to the one of semi-simple
complex Lie algebras. We will make the theorem more precise in
section~\ref{s:cluster-algebras-associated-with-quivers} below (for simply laced root systems).

\subsection{First example} \label{ss:first-example}
In order to illustrate the description and the theorem, we present
\cite{Zelevinsky07a} the cluster algebra $\ca_{A_2}$ associated
with the root system $A_2$. By definition, it is generated as
a $\Q$-algebra by the cluster variables $x_m$, $m\in\Z$, submitted
to the {\em exchange relations}\index{exchange relations}
\[
x_{m-1} x_{m+1} = 1+x_m \ko m\in \Z.
\]
Its clusters are by definition the pairs of consecutive cluster variables
$\{ x_m, x_{m+1}\}$, $m\in \Z$. The initial cluster is $\{x_1, x_2\}$ and
two clusters are linked by a mutation if and only if they share exactly
one variable.

The exchange relations allow one to write each cluster variable as
a rational expression in the initial variables $x_1$, $x_2$ and
thus to identify the algebra $\ca_{A_2}$ with a subalgebra of
the field $\Q(x_1, x_2)$. In order to make this subalgebra explicit,
let us compute the cluster variables $x_m$ for $m\geq 3$. We have:
\begin{align}
x_3 &= \frac{1+x_2}{x_1} \\
x_4 &= \frac{1+x_3}{x_2} = \frac{x_1 + 1 + x_2}{x_1 x_2}\\
x_5 &= \frac{1+x_4}{x_3}
= \frac{x_1 x_2 + x_1 +1 +x_2}{x_1 x_2}\; \div \; \frac{1+x_2}{x_1}
= \frac{1+x_1}{x_2} \; . \label{eq:laurent}
\end{align}
Notice that, contrary to what one might expect, the denominator
in~(\ref{eq:laurent}) remains a monomial! In fact, each cluster variable
in an arbitrary cluster algebra is a Laurent polynomial, cf.
Theorem~\ref{thm:class-finite-type}
below. Let us continue the computation:
\begin{align}
x_6 &= \frac{1+x_5}{x_4}
     =  \frac{x_2+1+x_1}{x_2}\; \div \;\frac{x_1+1+x_2}{x_1 x_2} = x_1 \\
x_7 &= (1+x_1)\; \div \; \frac{1+x_1}{x_2}  = x_2.
\end{align}
It is now clear that the sequence of the $x_m$, $m\in\Z$, is $5$-periodic
and that the number of cluster variables is indeed finite and equal to $5$.
In addition to the two initial variables $x_1$ and $x_2$, we have three
non initial variables $x_3$, $x_4$ and $x_5$. By examining their denominators
we see that they are in natural bijection with the positive roots
$\alpha_1$, $\alpha_1+\alpha_2$, $\alpha_2$ of the root system $A_2$. This
generalizes to an arbitrary Dynkin diagram, cf. Theorem~\ref{thm:class-finite-type}.

\subsection{Cluster algebras of rank $2$}\index{cluster algebra!of rank $2$}
To each pair of positive integers $(b,c)$,
there is associated a cluster algebra $\ca_{(b,c)}$. It is defined in analogy
with $\ca_{A_2}$ by replacing the exchange relations with
\[
x_{m-1} x_{m+1} = \left\{ \begin{array}{ll} x_m^b +1 & \mbox{if $m$ is odd, } \\
                                            x_m^c +1 & \mbox{if $m$ is even.} \end{array} \right.
\]
The algebra $\ca_{(b,c)}$ has only a finite number of cluster variables
if and only if we have $bc\leq 3$. In other words, if and only if the
matrix
\[
\left[ \begin{array}{cc} 2 & -b \\ -c & 2 \end{array} \right]
\]
is the Cartan matrix of a finite root system $\Phi$ of rank $2$.
The reader is invited to check that in this case, the non initial cluster
variables are still in natural bijection with the positive
roots of $\Phi$.

\section{Cluster algebras associated with quivers}
\label{s:cluster-algebras-associated-with-quivers}

\subsection{Quiver mutation} A {\em quiver}\index{quiver} is an oriented graph, i.e. a
quadruple $Q=(Q_0, Q_1, s, t)$ formed by a set of vertices $Q_0$, a set
of arrows $Q_1$ and two maps $s$ and $t$ from $Q_1$ to $Q_0$ which send
an arrow $\alpha$ respectively to its source $s(\alpha)$ and its
target $t(\alpha)$. In practice, a quiver is given by a picture
as in the following example
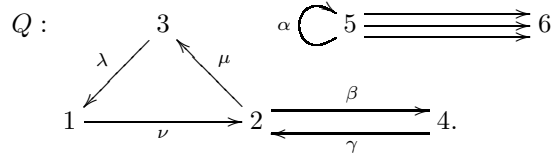
\[ Q:
\xymatrix{ & 3 \ar[ld]_\lambda & & 5 \ar@(dl,ul)[]^\alpha \ar@<1ex>[rr] \ar[rr] \ar@<-1ex>[rr] & & 6 \\
  1 \ar[rr]_\nu & & 2 \ar@<1ex>[rr]^\beta \ar[ul]_\mu & & 4.
  \ar@<1ex>[ll]^\gamma }
\]
An arrow $\alpha$ whose source and target coincide is a {\em loop}\index{loop}; a
{\em $2$-cycle}\index{$2$-cycle} is a pair of distinct arrows $\beta$ and $\gamma$ such that
$s(\beta)=t(\gamma)$ and $t(\beta)=s(\gamma)$. Similarly, one defines
{\em $n$-cycles} for any positive integer $n$. A vertex $i$ of a quiver
is a {\em source}\index{source of a quiver} (respectively a {\em sink}\index{sink of a quiver}) if there is no arrow
with target $i$ (respectively with source $i$).

By convention, in the sequel, by a quiver we always mean a finite
quiver without loops nor $2$-cycles whose set of vertices is the
set of integers from $1$ to $n$ for some $n\geq 1$. Up to an
isomorphism fixing the vertices such a quiver $Q$ is given by
the {\em skew-symmetric matrix $B=B_Q$} whose coefficient $b_{ij}$ is
the difference between the number of arrows from $i$ to $j$ and
the number of arrows from $j$ to $i$ for all $1\leq i,j\leq n$.
Conversely, each skew-symmetric matrix $B$ with integer coefficients
comes from a quiver.

Let $Q$ be a quiver and $k$ a vertex of $Q$. The {\em mutation $\mu_k(Q)$}\index{mutation!of a quiver}
is the quiver obtained from $Q$ as follows:
\begin{itemize}
\item[1)] for each subquiver $\xymatrix{i \ar[r]^\beta & k \ar[r]^\alpha & j}$,
we add a new arrow $[\alpha\beta]: i \to j$;
\item[2)] we reverse all arrows with source or target $k$;
\item[3)] we remove the arrows in a maximal set of pairwise disjoint $2$-cycles.
\end{itemize}
For example, if $k$ is a source or a sink of $Q$, then the
mutation at $k$ simply reverses all the arrows incident with $k$. In general,
if $B$ is the skew-symmetric matrix associated with $Q$ and $B'$ the one
associated with $\mu_k(Q)$, we have
\begin{equation} \label{eq:matrix-mutation}
b'_{ij} = \left\{ \begin{array}{ll}
-b_{ij} & \mbox{if $i=k$ or $j=k$~;} \\
b_{ij}+\sgn(b_{ik})\max(0, b_{ik} b_{kj}) & \mbox{else.}
\end{array} \right.
\end{equation}
This is the {\em matrix mutation rule}\index{mutation!of a matrix} for skew-symmetric
(more generally: skew-symmetri\-zable) matrices introduced by
Fomin-Zelevinsky in \cite{FominZelevinsky02},
cf. also \cite{FominZelevinsky07}.

One checks easily that $\mu_k$ is an involution. For example,
the quivers
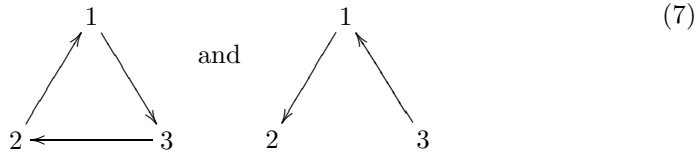
\begin{equation} \label{quiver1}
\begin{xy} 0;<0.3pt,0pt>:<0pt,-0.3pt>::
(94,0) *+{1} ="0",
(0,156) *+{2} ="1",
(188,156) *+{3} ="2",
"1", {\ar"0"},
"0", {\ar"2"},
"2", {\ar"1"},
\end{xy}
\begin{minipage}{1cm}
\vspace*{1cm}
\begin{center} and \end{center}
\end{minipage}
\begin{xy} 0;<0.3pt,0pt>:<0pt,-0.3pt>::
(92,0) *+{1} ="0",
(0,155) *+{2} ="1",
(188,155) *+{3} ="2",
"0", {\ar"1"},
"2", {\ar"0"},
\end{xy}
\end{equation}
are linked by a mutation at the vertex $1$. Notice that these
quivers are drastically different: The first one is a cycle,
the second one the Hasse diagram of a linearly ordered set.

Two quivers are {\em mutation equivalent}\index{mutation equivalent} if they are linked
by a finite sequence of mutations. For example, it is an easy
exercise to check that any two orientations of a tree are
mutation equivalent. Using the quiver mutation applet \cite{KellerQuiverMutationApplet}
or the package \cite{MusikerStump11} one can check
that the following three quivers are mutation equivalent
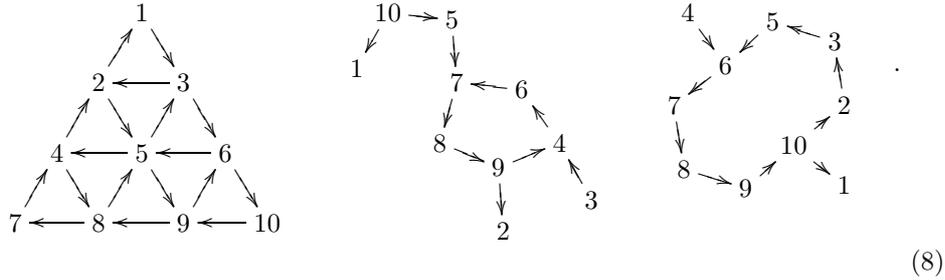
\begin{equation} \label{quiver3}
\begin{xy} 0;<0.6pt,0pt>:<0pt,-0.6pt>::
(79,0) *+{1} ="0",
(52,44) *+{2} ="1",
(105,44) *+{3} ="2",
(26,88) *+{4} ="3",
(79,88) *+{5} ="4",
(131,88) *+{6} ="5",
(0,132) *+{7} ="6",
(52,132) *+{8} ="7",
(105,132) *+{9} ="8",
(157,132) *+{10} ="9",
"1", {\ar"0"},
"0", {\ar"2"},
"2", {\ar"1"},
"3", {\ar"1"},
"1", {\ar"4"},
"4", {\ar"2"},
"2", {\ar"5"},
"4", {\ar"3"},
"6", {\ar"3"},
"3", {\ar"7"},
"5", {\ar"4"},
"7", {\ar"4"},
"4", {\ar"8"},
"8", {\ar"5"},
"5", {\ar"9"},
"7", {\ar"6"},
"8", {\ar"7"},
"9", {\ar"8"},
\end{xy}
\quad\quad
\begin{xy} 0;<0.3pt,0pt>:<0pt,-0.3pt>::
(0,70) *+{1} ="0",
(183,274) *+{2} ="1",
(293,235) *+{3} ="2",
(253,164) *+{4} ="3",
(119,8) *+{5} ="4",
(206,96) *+{6} ="5",
(125,88) *+{7} ="6",
(104,164) *+{8} ="7",
(177,194) *+{9} ="8",
(39,0) *+{10} ="9",
"9", {\ar"0"},
"8", {\ar"1"},
"2", {\ar"3"},
"3", {\ar"5"},
"8", {\ar"3"},
"4", {\ar"6"},
"9", {\ar"4"},
"5", {\ar"6"},
"6", {\ar"7"},
"7", {\ar"8"},
\end{xy}
\quad\quad
\begin{xy} 0;<0.3pt,0pt>:<0pt,-0.3pt>::
(212,217) *+{1} ="0",
(212,116) *+{2} ="1",
(200,36) *+{3} ="2",
(17,0) *+{4} ="3",
(123,11) *+{5} ="4",
(64,66) *+{6} ="5",
(0,116) *+{7} ="6",
(12,196) *+{8} ="7",
(89,221) *+{9} ="8",
(149,166) *+{10} ="9",
"9", {\ar"0"},
"1", {\ar"2"},
"9", {\ar"1"},
"2", {\ar"4"},
"3", {\ar"5"},
"4", {\ar"5"},
"5", {\ar"6"},
"6", {\ar"7"},
"7", {\ar"8"},
"8", {\ar"9"},
\end{xy}
\begin{minipage}{1cm}
\vspace*{1.5cm}
\begin{center} . \end{center}
\end{minipage}
\end{equation}
The common {\em mutation class}\index{mutation class} of these quivers contains 5739 quivers
(up to isomorphism). The mutation class of `most' quivers is infinite.
The classification of the quivers having a finite mutation class
was achieved by by Felikson-Shapiro-Tumarkin
\cite{FeliksonShapiroTumarkin08} \cite{FeliksonShapiroTumarkin10}:
in addition to the quivers associated
with triangulations of surfaces (with boundary and marked points,
cf. \cite{FominShapiroThurston08})
the list contains $11$ exceptional quivers, the largest of which is
in the mutation class of the quivers~(\ref{quiver3}).

\subsection{Seed mutation, cluster algebras} \label{ss:seed-mutation}
Let $n\geq 1$ be an integer and $\cf$ the field $\Q(x_1, \ldots, x_n)$
generated by $n$ indeterminates $x_1, \ldots, x_n$. A {\em seed}\index{seed}
(more precisely: {\em $X$-seed}\index{$X$-seed}) is a pair $(R,u)$, where $R$
is a quiver and $u$ a sequence $u_1, \ldots, u_n$ of elements
of $\cf$ which freely generate the field $\cf$. If $(R,u)$ is
a seed and $k$ a vertex of $R$, the {\em mutation $\mu_k(R,u)$}\index{mutation!of a seed} is
the seed $(R',u')$, where $R'=\mu_k(R)$ and $u'$ is obtained
from $u$ by replacing the element $u_k$ by the element $u_k'$
defined by the {\em exchange relation}\index{exchange relation}
\begin{equation} \label{eq:echange}
u_k' u_k= \prod_{s(\alpha)=k} u_{t(\alpha)} + \prod_{t(\alpha)=k} u_{s(\alpha)} \ko
\end{equation}
where the sums range over all {\em arrows} of $R$ with source $k$
respectively target $k$. Notice that, if $B$ is the skew-symmetric matrix associated
with $R$, we can rewrite the exchange relation as
\begin{equation} \label{eq:exchange-B}
u_k' u_k= \prod_{i} u_{i}^{[b_{ik}]_+} + \prod_{i} u_{i}^{[-b_{ik}]_+} \ko
\end{equation}
where, for a real number $x$, we write $[x]_+$ for $\max(x,0)$.
One checks that $\mu_k^2(R,u)=(R,u)$. For example, the mutations of the
seed
\[
(\xymatrix{1 \ar[r] & 2 \ar[r] & 3} \ko \{x_1, x_2, x_3\})
\]
with respect to the vertices $1$ and $2$ are the seeds
\begin{align} \label{eq:seed2-for-A3}
(\xymatrix{1 & 2 \ar[l] \ar[r] & 3} \ko \{\frac{1+x_2}{x_1}, x_2, x_3\})
\quad\mbox{ and } \\
\label{eq:seed3-for-A3}
(\xymatrix{1 \ar@/^1pc/[rr] & 2 \ar[l] & 3 \ar[l]} \ko \{ x_1, \frac{x_1+x_3}{x_2}, x_3\}).
\end{align}

Let us fix a quiver $Q$. The {\em initial seed}\index{initial seed} of $Q$ is
$(Q,\{x_1, \ldots, x_n\})$. A {\em cluster}\index{cluster} associated with $Q$
is a sequence $u$ which appears in a seed $(R,u)$ obtained from
the initial seed by iterated mutation. The {\em cluster variables}\index{cluster variable}
are the elements of the clusters. The {\em cluster algebra $\ca_Q$}\index{cluster algebra!definition}
is the $\Q$-subalgebra of $\cf$ generated by the cluster variables.
Clearly, if $(R,u)$ is a seed associated with $Q$, the natural
isomorphism
\[
\Q(u_1, \ldots, u_n) \iso \Q(x_1, \ldots, x_n)
\]
induces an isomorphism of $\ca_R$ onto $\ca_Q$ which preserves
the cluster variables and the clusters. Thus, the cluster algebra
$\ca_Q$ is an invariant of the mutation class of $Q$.
It is useful to introduce a combinatorial object which encodes
the recursive construction of the seeds: the {\em exchange graph}\index{exchange graph}.
By definition, its vertices are the isomorphism classes of seeds
(isomorphisms of seeds renumber the vertices and the variables
simultaneously) and
its edges correspond to mutations. For example, the exchange graph
obtained from the quiver
$
Q: \xymatrix{1 \ar[r] & 2 \ar[r] & 3}
$
is the $1$-skeleton of the Stasheff associahedron \cite{Stasheff63}:
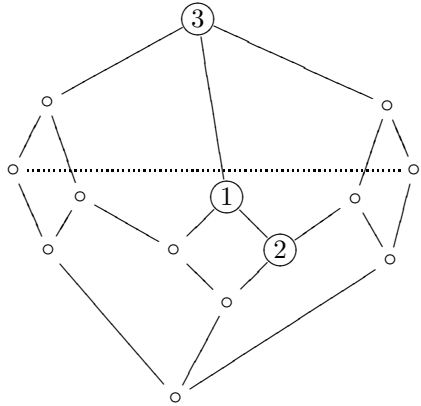
\[
\begin{xy} 0;<0.4pt,0pt>:<0pt,-0.4pt>::
(173,0) *+<8pt>[o][F]{3} ="0",
(0,143) *+{\circ} ="1",
(63,168) *+{\circ} ="2",
(150,218) *+{\circ} ="3",
(250,218) *+<8pt>[o][F]{2} ="4",
(375,143) *+{\circ} ="5",
(350,82) *+{\circ} ="6",
(152,358) *+{\circ} ="7",
(200,168) *+<8pt>[o][F]{1} ="8",
(200,268) *+{\circ} ="9",
(32,79) *+{\circ} ="10",
(33,218) *+{\circ} ="11",
(320,170) *+{\circ} ="12",
(353,228) *+{\circ} ="13",
"0", {\ar@{-}"6"},
"0", {\ar@{-}"8"},
"0", {\ar@{-}"10"},
"1", {\ar@{.}"5"},
"1", {\ar@{-}"10"},
"11", {\ar@{-}"1"},
"2", {\ar@{-}"3"},
"10", {\ar@{-}"2"},
"2", {\ar@{-}"11"},
"3", {\ar@{-}"8"},
"9", {\ar@{-}"3"},
"8", {\ar@{-}"4"},
"4", {\ar@{-}"9"},
"4", {\ar@{-}"12"},
"6", {\ar@{-}"5"},
"5", {\ar@{-}"13"},
"12", {\ar@{-}"6"},
"9", {\ar@{-}"7"},
"11", {\ar@{-}"7"},
"13", {\ar@{-}"7"},
"13", {\ar@{-}"12"},
\end{xy}
\]
Here the vertex $1$ corresponds to the initial seed and the vertices
$2$ and $3$ to the seeds~(\ref{eq:seed2-for-A3}) and~(\ref{eq:seed3-for-A3}).
For analogous polytopes associated with the other Dynkin diagrams, we
refer to \cite{ChapotonFominZelevinsky02}.

Let $Q$ be a connected quiver. If its underlying graph is
a simply laced Dynkin diagram $\Delta$, we say that $Q$ is
a {\em Dynkin quiver of type $\Delta$}\index{Dynkin quiver}.

\begin{theorem}[\cite{FominZelevinsky03}] \label{thm:class-finite-type}
\label{thm:cluster-finite-classification}
\begin{itemize}
\item[a)] Each cluster variable of $\ca_Q$ is a Laurent polynomial with integer
coefficients \cite{FominZelevinsky02}.
\item[b)] The cluster algebra $\ca_Q$ has only a finite number of cluster variables
if and only if $Q$ is mutation equivalent to a Dynkin quiver $R$. In this case, the
underlying graph $\Delta$ of $R$ is unique up to isomorphism and is called the
{\em cluster type\index{cluster type} of $Q$}.
\item[c)] If $Q$ is a Dynkin quiver of type $\Delta$, then the non initial
cluster variables of $\ca_Q$ are in bijection with the positive roots of
the root system $\Phi$ of $\Delta$; more precisely, if $\alpha_1$, \ldots, $\alpha_n$
are the simple roots, then for each positive root $\alpha= d_1 \alpha_1 + \cdots + d_n \alpha_n$,
there is a unique non initial cluster variable $X_\alpha$ whose denominator
is $x_1^{d_1} \ldots x_n^{d_n}$.
\end{itemize}
\end{theorem}

Statement a) is usually refered to as the {\em Laurent phenomenon}\index{Laurent phenomenon}.
A {\em cluster monomial}\index{cluster monomial} is a product of non negative powers of cluster
variables belonging to the same cluster. The construction of a `canonical
basis' of the cluster algebra $\ca_Q$ is an important and largely open
problem, cf. for example
\cite{FominZelevinsky02}
\cite{ShermanZelevinsky04} \cite{Dupont11} \cite{Cerulli11} \cite{GeissLeclercSchroeer12}
\cite{MusikerSchifflerWilliams11} \cite{Lampe10} \cite{Lampe11} \cite{HernandezLeclerc11}.
It is expected that such a basis should contain all cluster monomials.
Whence the following conjecture.

\begin{conjecture}[\cite{FominZelevinsky02}] \label{conj:independence}
The cluster monomials are linearly independent over the field $\Q$.
\end{conjecture}

The conjecture was recently proved in \cite{CerulliKellerLabardiniPlamondon12}
using the additive categorification of \cite{Plamondon11a} and techniques
from \cite{Cerulli11a} \cite{CerulliLabardini11}. It is expected to hold more generally
for cluster algebras associated with valued quivers,
cf.~section~\ref{ss:valued-quivers} below. It is shown
for a certain class of valued quivers by L.~Demonet \cite{Demonet10} \cite{Demonet11}.
For special classes of quivers, a basis containing the cluster monomials is known:
If $Q$ is a Dynkin quiver, one knows \cite{CalderoKeller08} that
the cluster monomials form a basis of $\ca_Q$. If $Q$ is {\em acyclic}\index{quiver!acyclic},
i.e. does not have any oriented cycles, then Geiss-Leclerc-Schr\"oer
\cite{GeissLeclercSchroeer07b} show the existence of a `generic basis'
containing the cluster monomials.

\begin{conjecture}[\cite{FominZelevinsky03}]  \label{conj:positivity}
The cluster variables are Laurent polynomials with non negative
integer coefficients in the variables of each cluster.
\end{conjecture}

For quivers with two vertices, an explicit and manifestly positive
formula for the cluster variables is given in \cite{LeeSchiffler11}.
The technique of monoidal categorification developed by Leclerc
\cite{Leclerc08a} and Her\-nandez-Leclerc \cite{HernandezLeclerc10}
has recently allowed to prove the conjecture first for
the quivers of type $A_n$ and $D_4$, cf. \cite{HernandezLeclerc10},
and then for each bipartite quiver \cite{Nakajima11}, i.e. a
quiver where each vertex is a source or a sink. The positivity
of all cluster variables with respect to the initial seed of
an acyclic quiver is shown by Fan Qin \cite{Qin10} and by
Nakajima \cite[Appendix]{Nakajima11}. This is also
proved by Efimov \cite{Efimov11}, who moreover shows the
positivity of all cluster variables belonging to an acyclic
seed with respect to the initial variables of an arbitrary
quiver. Efimov combines the techniques of \cite{KontsevichSoibelman10} with 
those of \cite{Nagao10}. A proof of the full conjecture for
acyclic quivers using Nakajima quiver varieties is announced by 
Kimura--Qin \cite{KimuraQin12}.
The conjecture has been shown in a combinatorial way by Musiker-Schiffler-Williams
\cite{MusikerSchifflerWilliams11a} for all the quivers associated
with triangulations of surfaces (with boundary and marked points)
and by Di Francesco-Kedem \cite{DiFrancescoKedem09} for the quivers
and the cluster variables
associated with the $T$-system of type $A$, with respect to the
initial cluster.

We refer to \cite{FominZelevinsky03a} and \cite{FominZelevinsky07}
for numerous other conjectures on cluster algebras and to
\cite{DerksenWeymanZelevinsky10}, cf. also \cite{Nagao10} and
\cite{Plamondon11} \cite{Plamondon11a}, for the solution of a large
number of them using additive categorification.

\subsection{Cluster algebras associated with valued quivers}
\label{ss:valued-quivers}
A {\em valued quiver}\index{quiver!valued} is a quiver $Q$ endowed with a function $v: Q_1 \to \N^2$
such that
\begin{itemize}
\item[a)] there are no loops in $Q$,
\item[b)] there is at most one arrow between any two vertices of $Q$ and
\item[c)] there is a function $d: Q_0 \to \N$ such that $d(i)$ is strictly
positive for all vertices $i$ and, for each arrow $\alpha: i\to j$, we have
\[
d(i)\,v(\alpha)_1  = v(\alpha)_2 \, d(j) \ko
\]
where $v(\alpha)=(v(\alpha)_1, v(\alpha)_2)$.
\end{itemize}
For example, we have the valued quivers (we omit the labels $(1,1)$ from our pictures)
\[
\vec{B}_3 : \xymatrix{1 \ar[r] & 2 \ar[r]^{(1,2)} & 3 }
\quad\mbox{and}\quad
\vec{C}_3 : \xymatrix{1 \ar[r] & 2 \ar[r]^{(2,1)} & 3} \ko
\]
where possible functions $d$ are given by $d(1)=d(2)=2$, $d(3)=1$
respectively $d(1)=d(2)=1$, $d(3)=2$.
A valued quiver $(Q,v)$ is {\em equally valued}\index{quiver!equally valued} if we have
$v(\alpha)_1=v(\alpha)_2$ for each arrow $\alpha$. If
$Q$ is an ordinary quiver without loops nor $2$-cycles,
the {\em associated valued quiver} is the equally valued
quiver which has an arrow $\alpha: i \to j$ if there is at least
one arrow $i\to j$ in $Q$ and where $v(\alpha)=(m,m)$, where
$m$ is the number of arrows from $i$ to $j$ in $Q$. For example,
the equally valued quiver
\[
\xymatrix{1 \ar[r]^{(2,2)} & 2 } \quad\mbox{corresponds to the Kronecker quiver}\quad
\xymatrix{1 \ar@<0.5ex>[r]\ar@<-0.5ex>[r]  & 2.}
\]
In this way, the ordinary quivers without loops nor $2$-cycles
correspond bijectively to the equally valued quivers (up to
isomorphism fixing the vertices).
Let $Q$ be a valued quiver with vertex set $I$. We associate
an integer matrix $B=(b_{ij})_{i,j\in I}$ with it as follows
\[
b_{ij} = \left\{ \begin{array}{ll} 0 & \mbox{ if there is no arrow between $i$ and $j$;} \\
v(\alpha)_1 & \mbox{ if there is an arrow } \alpha: i \to j ; \\
-v(\alpha)_2 & \mbox{ if there is an arrow } \alpha: j \to i.
\end{array} \right.
\]
If $D$ is the diagonal $I\times I$-matrix with diagonal entries
$d_{ii}=d(i)$, $i\in I$, then the matrix $DB$ is skew-symmetric.
The existence of such a matrix $D$ means that the matrix
$B$ is {\em skew-symmetrizable}\index{matrix!skew-symmetrizable}. It is easy to
check that in this way, we obtain a bijection between the
skew-symmetrizable $I\times I$-matrices $B$ and the valued
quivers with vertex set $I$ (up to isomorphism fixing the vertices).
Using this bijection, we define
the {\em mutation of valued quivers}\index{mutation!of valued quivers} using Fomin-Zelevinsky's
matrix mutation rule~(\ref{eq:matrix-mutation}). For example,
the mutation at $2$ transforms the valued quiver
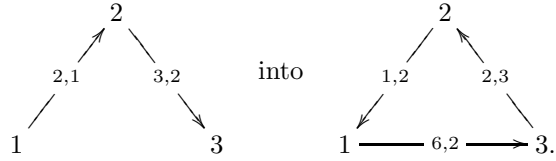
\[
\begin{xy} 0;<0.5pt,0pt>:<0pt,-0.5pt>::
(0,100) *+{1} ="0",
(75,0) *+{2} ="1",
(150,100) *+{3} ="2",
"0", {\ar|*+{\scriptstyle 2,1}"1"},
"1", {\ar|*+{\scriptstyle 3,2}"2"},
\end{xy}
\quad\raisebox{-0.8cm}{\mbox{into}}\quad
\begin{xy} 0;<0.5pt,0pt>:<0pt,-0.5pt>::
(0,100) *+{1} ="0",
(75,0) *+{2} ="1",
(150,100) *+{3.} ="2",
"1", {\ar|*+{\scriptstyle 1,2}"0"},
"0", {\ar|*+{\scriptstyle 6,2}"2"},
"2", {\ar|*+{\scriptstyle 2,3}"1"},
\end{xy}
\]
We extend the notion of an {\em ($X$-)seed}\index{seed!valued} $(R,u)$ by
now allowing the first component $R$ to be any valued quiver
and we extend the construction of seed mutation by
using the rule (\ref{eq:exchange-B}), where $B$ is
the skew-symmetrizable matrix associated with $R$.
For example the mutations of the seed
\[
(\xymatrix{1 \ar[r] & 2 \ar[r]^{(1,2)} & 3}, \{x_1, x_2, x_3\})
\]
at the vertices $1$ and $2$ are the seeds
\[
(\xymatrix{1  & 2 \ar[l] \ar[r]^{(1,2)} & 3}\ko \{\frac{1+x_2}{x_1}, x_2, x_3\})
\]
and
\[
(\xymatrix{1 \ar@/^1pc/[rr]^{(1,2)} & 2 \ar[l] & 3 \ar[l]^{(2,1)}} \ko \{x_1, \frac{x_1+x_3^2}{x_2}, x_3\}).
\]
Given a valued quiver $Q$, we define its associated {\em clusters}\index{cluster!of a valued quiver},
{\em cluster variables}\index{cluster variable!of a valued quiver},
{\em cluster monomials}\index{cluster monomial!of a valued quiver},
the {\em cluster algebra $\ca_Q$}\index{cluster algebra!of a valued quiver} and
the {\em exchange graph}\index{exchange graph!of a valued quiver} in complete analogy with the constructions
in section~\ref{ss:seed-mutation}. For example, the exchange graph
of the above quivers $\vec{B}_3$ and $\vec{C}_3$ is the
$3$rd cyclohedron \cite{ChapotonFominZelevinsky02}, with $4$ quadrilateral, $4$
pentagonal and $4$ hexagonal faces:
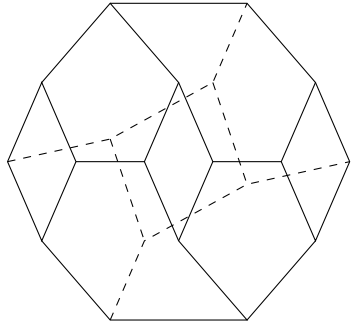
\begin{center}
\begin{tikzpicture}[scale=0.15]
\draw (0,0)--(3,7)--(6,0)--(3,-7)--(0,0);
\draw (6,0)--(12,0)--(15,7)--(9,14)--(3,7);
\draw (12,0)--(15,-7)--(21,-14)--(9,-14)--(3,-7);
\draw (15,-7)--(18,0)--(15,7);
\draw (18,0)--(24,0)--(27,7)--(21,14)--(9,14);
\draw (24,0)--(27,-7)--(30,0)--(27,7);
\draw (27,-7)--(21,-14);
\draw[dashed] (0,0)--(9,2)--(18,7)--(21,14);
\draw[dashed] (9,-14)--(12,-7)--(21,-2)--(30,0);
\draw[dashed] (9,2)--(12,-7);
\draw[dashed] (18,7)--(21,-2);
\end{tikzpicture}
\end{center}

Let $(Q,v)$ be a valued quiver with vertex set $I=Q_0$. Its
{\em associated Cartan matrix}\index{Cartan matrix} is the Cartan companion\index{Cartan companion}
\cite{FominZelevinsky03}
of the skew-symmetrizable matrix $B$ associated with $Q$.
Explicitly, it is the the $I\times I$-matrix $C$
whose coefficient $c_{ij}$ vanishes if there are no
arrows between $i$ and $j$, equals $2$ if $i=j$,
equals $-v(\alpha)_1$ if there is an arrow $\alpha: i \to j$
and equals $-v(\alpha)_2$ if there is an arrow $\alpha: j\to i$. Thus, the Cartan
matrix associated with the above valued quiver $\vec{B}_2$
equals
\[
\left[ \begin{array}{cc} 2 & -2 \\ -1 & 2 \end{array} \right].
\]
Fomin-Zelevinsky have shown in \cite{FominZelevinsky03} that
the the analogue of Theorem~\ref{thm:cluster-finite-classification}
holds for valued quivers. In particular, the Laurent phenomenon
holds and the cluster algebra associated with a valued quiver
$Q$ has only finitely many cluster variables iff $Q$ is
mutation-equivalent to a valued quiver whose associated
Cartan matrix corresponds to a finite root system.

For valued quivers, the independence conjecture~\ref{conj:independence}
is open except for the valued quivers treated by Demonet \cite{Demonet10} \cite{Demonet11}.
The positivity conjecture~\ref{conj:positivity} is open except
in rank two, where it was shown by Dupont in \cite{Dupont09}.

\section{Cluster algebras of geometric type}
\label{s:cluster-algebras-of-geometric-type}
\index{cluster algebra!of geometric type}

We will slightly generalize the definition of
section~\ref{s:cluster-algebras-associated-with-quivers}
in order to obtain the class of `skew-symmetrizable cluster algebras of geometric type'. This
class contains many algebras of geometric origin which are equipped with
`dual semi-canonical bases' \cite{Lusztig00}. The construction of a large part of such a
basis in \cite{GeissLeclercSchroeer06} is one of the most remarkable
applications of cluster algebras so far.

We refer to section~\ref{ss:cluster-alg-semifield}
for the definition of the `skew-symmetrizable cluster algebras with
coefficients in a semi-field', which constitute so
far the most general class considered.

\subsection{Definition} \label{ss:def-cluster-alg-coeff}
Let $1 \leq n \leq m$ be integers. Let $\tilde{Q}$ be an {\em ice quiver\index{quiver!ice}\index{ice quiver}
of type $(n,m)$}, i.e. a quiver with $m$ vertices and which does not
have any arrows between vertices $i$, $j$ which are both strictly
greater than $n$. The {\em principal part}\index{ice quiver!principal part of} of $\tilde{Q}$ is the
full subquiver $Q$ whose vertices are $1$, \ldots, $n$ (a subquiver
is {\em full}\index{subquiver!full} if, with any two vertices, it contains all the
arrows linking them). The vertices $n+1$, \ldots, $m$ are
called the {\em frozen vertices}\index{frozen vertex}. The {\em cluster algebra\index{cluster algebra!of an ice quiver} associated
with the ice quiver $\tilde{Q}$}
\[
\ca_{\tilde{Q}} \subset \Q(x_1, \ldots, x_m)
\]
is defined in the same manner as the cluster algebra associated
with a quiver (section~\ref{s:cluster-algebras-associated-with-quivers})
but
\begin{itemize}
\item only mutations with respect to non frozen vertices are allowed
and no arrows between frozen vertices are added in the mutations;
\item the variables $x_{n+1}$, \ldots, $x_m$, which belong to
all clusters, are called {\em coefficients}\index{coefficients!of a cluster algebra} rather than cluster
variables;
\item the {\em cluster type}\index{cluster type!of an ice quiver} of the ice quiver is that of
its principal part (if it is defined).
\end{itemize}
Notice that the datum of $\tilde{Q}$ is equivalent to that
of the integer $m\times n$-matrix $\tilde{B}$ whose coefficient
$b_{ij}$ is the difference of the number of arrows from $i$ to
$j$ minus the number of arrows from $j$ to $i$ for all
$1\leq i\leq m$ and all $1\leq j\leq n$. The top $n\times n$
part $B$ of $\tilde{B}$ is called its {\em principal part}\index{principal part}.
In complete analogy, one defines the cluster algebra associated
with a {\em valued ice quiver} respectively with an integer
$m\times n$-matrix whose principal part is skew-symmetrizable.

We have the following sharpening of the Laurent phenomenon
proved in Proposition~11.2 of \cite{FominZelevinsky03}.

\begin{theorem}[\cite{FominZelevinsky03}] \label{thm:sharpened-Laurent}
Each cluster variable in $\ca_{\tilde{Q}}$ is a Laurent polynomial
in the initial variables $x_1$, \ldots, $x_n$ with
coefficients in $\Z[x_{n+1}, \ldots, x_m]$.
\end{theorem}

Often one considers localizations of $\ca_{\tilde{Q}}$ obtained
by inverting some or all of the coefficients. If $K$ is an extension
field of $\Q$ and $A$ a commutative $K$-algebra without zero
divisors, a {\em cluster structure\index{cluster structure} of type $\tilde{Q}$ on $A$} is given
by an isomorphism $\phi$ from $\ca_{\tilde{Q}} \ten_\Q K$ onto
$A$. Such an isomorphism is determined by the images of the
coefficients and of the initial cluster variables $\phi(x_i)$,
$1\leq i\leq m$. We call the datum of $\tilde{Q}$ and of
the $\phi(x_i)$ an {\em initial seed\index{initial seed} for $A$}. The following
proposition is a reformulation of Proposition~11.1
of \cite{FominZelevinsky03}, cf. also Proposition~1 of \cite{Scott06}:

\begin{proposition} \label{prop:cluster-structure}
Let $X$ be a rational quasi-affine irreducible algebraic
variety over $\C$. Let $\tilde{Q}$ be an ice quiver of type $(m,n)$.
Assume that we are given a regular function $\phi_c$ on $X$ for each
coefficient $c=x_i$, $n< i\leq m$, and a regular function $\phi_x$ on $X$ for
each cluster variable $x$ of $\ca_{\tilde{Q}}$ such that
\begin{itemize}
\item[a)] the dimension of $X$ equals $m$;
\item[b)] the functions $\phi_x$ and $\phi_c$ generate the coordinate
algebra $\C[X]$;
\item[c)]
the correspondence $x \mapsto \phi_x$, $c\mapsto \phi_c$ takes
each exchange relation of $\ca_{\tilde{Q}}$ to an equality in $\C[X]$.
\end{itemize}
Then the correspondence $x \mapsto \phi_x$, $c \mapsto \phi_c$
extends to an algebra isomorphism $\phi : \ca_{\tilde{Q}}\ten_\Q \C \iso \C[X]$
so that $\C[X]$ carries a cluster algebra structure of type $\tilde{Q}$ with
initial seed $\phi_{x_i}$, $1\leq i\leq m$.
\end{proposition}

\subsection{Example: Planes in a vector space}
Let $n\geq 1$ be an integer. Let $A$ be the algebra of polynomial
functions on the cone over the Grassmannian\index{Grassmannian!of planes} of planes in $\C^{n+3}$.
This algebra is generated by the Pl\"ucker coordinates $x_{ij}$, $1\leq i<j\leq n+3$,
subject to the Pl\"ucker relations: for each quadruple of integers
$i<j<k<l$ between $1$ and $n+3$, we have
\begin{equation} \label{eq:Plucker}
x_{ik} x_{jl} = x_{ij} x_{kl} + x_{jk} x_{il}.
\end{equation}
Notice that the monomials in this relation are naturally associated
with the diagonals and the sides of the square
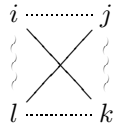
\[
\xymatrix{ i \ar@{.}[r] \ar@{~}[d] \ar@{-}[dr] & j \ar@{-}[dl] \ar@{~}[d] \\
l \ar@{.}[r] & k}
\]
The idea is to interpret this relation as an exchange relation
in a cluster algebra with coefficients. In order to describe
this algebra, let us consider, in the euclidean plane, a regular
polygon $P$ whose vertices are numbered from $1$ to $n+3$.
Consider the variable $x_{ij}$ as associated with the
segment $[ij]$ which links the vertices $i$ and~$j$.

\begin{proposition}[\protect{\cite[Example 12.6]{FominZelevinsky03}}]
\label{prop:plans} The algebra $A$ has a cluster algebra
structure such that
\begin{itemize}
\item[-] the coefficients are the variables $x_{ij}$ associated
with the sides of $P$;
\item[-] the cluster variables are the variables $x_{ij}$ associated
with the diagonals of $P$;
\item[-] the clusters are the $n$-tuples of cluster variables
corresponding to diagonals which form a triangulation of $P$.
\end{itemize}
Moreover, the exchange relations are exactly the Pl\"ucker
relations and the cluster type is $A_n$.
\end{proposition}

A triangulation of $P$ determines an initial seed for the
cluster algebra and the exchange relations satisfied by
the initial cluster variables determine the ice quiver
$\tilde{Q}$. For example, one can check that in the
following picture, the triangulation and the ice
quiver (whose frozen vertices are in boxes) correspond
to each other.
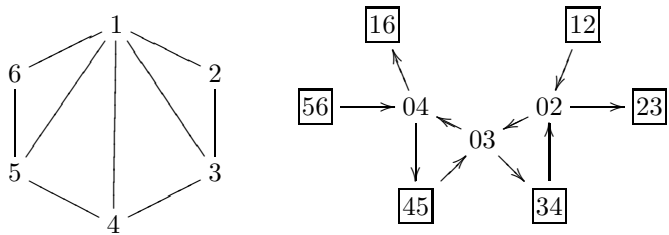
\[
\begin{xy} 0;<0.5pt,0pt>:<0pt,-0.5pt>::
(76,0) *+{1} ="0",
(150,37) *+{2} ="1",
(150,112) *+{3} ="2",
(74,151) *+{4} ="3",
(0,112) *+{5} ="4",
(0,37) *+{6} ="5",
(275,0) *+{\framebox[3ex]{16}} ="6",
(300,62) *+{04} ="7",
(350,87) *+{03} ="8",
(400,62) *+{02} ="9",
(425,0) *+{\framebox[3ex]{12}} ="10",
(225,62) *+{\framebox[3ex]{56}} ="11",
(300,137) *+{\framebox[3ex]{45}} ="12",
(400,137) *+{\framebox[3ex]{34}} ="13",
(475,62) *+{\framebox[3ex]{23}} ="14",
"0", {\ar@{-}"1"},
"0", {\ar@{-}"2"},
"0", {\ar@{-}"3"},
"0", {\ar@{-}"4"},
"5", {\ar@{-}"0"},
"1", {\ar@{-}"2"},
"2", {\ar@{-}"3"},
"3", {\ar@{-}"4"},
"4", {\ar@{-}"5"},
"7", {\ar"6"},
% "6", {\ar"11"},
"8", {\ar"7"},
"11", {\ar"7"},
"7", {\ar"12"},
"9", {\ar"8"},
"12", {\ar"8"},
"8", {\ar"13"},
"10", {\ar"9"},
"13", {\ar"9"},
"9", {\ar"14"},
% "14", {\ar"10"},
\end{xy}
\]
The hypotheses of proposition~\ref{prop:cluster-structure} are
straightforward to check in this example.
Many other (homogeneous) coordinate algebras of classical algebraic
varieties admit cluster algebra structures (or `upper cluster algebra
structures') and in particular the Grassmannians \cite{Scott06},
cf. section~\ref{ss:Gr-3-6} below, and
the double Bruhat cells \cite{BerensteinFominZelevinsky05}. Some
of these algebras have only finitely many cluster variables and
thus a well-defined cluster type. Here is a list of some examples of
varieties and their cluster type extracted from \cite{FominZelevinsky03a},
where $N$ denotes a maximal unipotent subgroup of the
corresponding reductive algebraic group:
\begin{center}
\begin{tabular}{|c|c|c|c|} \hline
$Gr_{2,n+3}$ & $Gr_{3,6}$ & $Gr_{3,7}$ & $Gr_{3,8}$  \\ \hline
$A_n$ & $D_4$ & $E_6$ & $E_8$  \\ \hline
\end{tabular}
\end{center}
\begin{center}
\begin{tabular}{|c|c|c|c|c|c|}\hline
$SL_3/N$ & $SL_4/N$ & $SL_5/N$ & $Sp_4/N$ & $SL_2$ & $SL_3$ \\ \hline
$A_1$ & $A_3$ & $D_6$ & $B_2$ & $A_1$ & $D_4$ \\ \hline
\end{tabular}
\end{center}
\bigskip

A theorem analogous to proposition~\ref{prop:plans} for `reduced
double Bruhat cells' is due to Yang and Zelevinsky \cite{YangZelevinsky08}.
They thus obtain a cluster algebra (with principal coefficients) with
an explicit description of the cluster variables for each Dynkin
diagram.

\subsection{Example: The Grassmannian $Gr(3,6)$}  \label{ss:Gr-3-6}
\index{Grassmannian!of subspaces}
Let us consider the cone $X$ over
the Pl\"ucker embedding of the variety $Gr(3,6)$ of $3$-dimensional
subspaces in $6$-dimensional complex space $\C^6$, considered
as a space of rows. The Pl\"ucker coordinates of the subspace
generated by the rows of a complex $3\times 6$-matrix are the
$3\times 3$-minors of the matrix, \ie
the determinants $D(\mathbf{j})$ of the $3\times 3$-submatrices
formed by the columns with indices in a $3$-element
subset $\mathbf{j}$ of $\{1, \ldots, 6\}$. It is a particular case
of Scott's theorem \cite{Scott06}, cf.~also Example~10.3
of \cite{GeissLeclercSchroeer08b}, that the algebra $\C[X]$ admits
a cluster algebra structure of the type
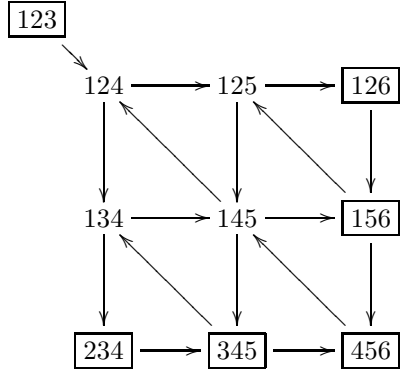
\[
\begin{xy} 0;<0.5pt,0pt>:<0pt,-0.5pt>::
(50,50) *+{124} ="0",
(150,50) *+{125} ="1",
(250,50) *+{\framebox[5ex]{126}} ="2",
(50,150) *+{134} ="3",
(150,150) *+{145} ="4",
(250,150) *+{\framebox[5ex]{156}} ="5",
(50,250) *+{\framebox[5ex]{234}} ="6",
(150,250) *+{\framebox[5ex]{345}} ="7",
(250,250) *+{\framebox[5ex]{456}} ="8",
(0,0) *+{\framebox[5ex]{123}} ="9",
"0", {\ar"1"},
"0", {\ar"3"},
"4", {\ar"0"},
"9", {\ar"0"},
"1", {\ar"2"},
"1", {\ar"4"},
"5", {\ar"1"},
"2", {\ar"5"},
"3", {\ar"4"},
"3", {\ar"6"},
"7", {\ar"3"},
"4", {\ar"5"},
"4", {\ar"7"},
"8", {\ar"4"},
"5", {\ar"8"},
"6", {\ar"7"},
"7", {\ar"8"},
\end{xy}
\]
whose initial seed is given by the minors $D(\mathbf{j})$ associated with
the vertices $\mathbf{j}$ of this quiver (frozen vertices appear in boxes).
If we mutate the principal part of this quiver at the vertex $124$, we obtain
a Dynkin quiver of type $D_4$, which is thus the cluster type of this cluster algebra.
It admits $4+12=16$ cluster variables. As shown in \cite{Scott06}, fourteen
among these are minors and the remaining two are
\[
X_1 = |P_1 \wedge Q_1, P_2 \wedge Q_2, P_3\wedge Q_3| \quad\mbox{and}\quad
X_2 = |Q_1 \wedge P_2, Q_2\wedge P_3, Q_3 \wedge P_1| \ko
\]
where we denote the columns of our matrix by $P_1$, $Q_1$, $P_2$, $Q_2$,
$P_3$, $Q_3$ (in this order) and write $| |$ for the determinant.
In this cluster algebra, we have
the remarkable identity \cite{Zagier11}
\begin{equation} \label{eq:Zagier}
|P_1 P_2 Q_2| | P_2 P_3 Q_3| |P_3 P_1 Q_1| -
|P_1 P_2 Q_1| |P_2 P_3 Q_2| |P_3 P_1 Q_3| = |P_1 P_2 P_3| \, X_1 \ko
\end{equation}
which we can rewrite as
\[
D(134) D(356) D(125) + D(123) D(345) D(156) = D(135) X_1.
\]
This is in fact an exchange relation in our cluster algebra (many
thanks to B.~Leclerc for pointing this out): Indeed, if we successively
mutate the initial seed at the vertices $124$ and $145$, we obtain the cluster
\[
D(135),D(125),D(356),D(134)
\]
(exercise: compute the corresponding quiver!) and if we now mutate at the
variable $D(135)$, we obtain $X_1$ and the exchange relation~(\ref{eq:Zagier}).
This relation appears implicitly in \cite{Goncharov91} and finding a
suitable generalization to higher dimensions would be
of interest in view of Zagier's conjecture \cite{Zagier91}.

\subsection{Example: Rectangular matrices} \label{ss:polynomial-algebras}
Polynomial algebras admit many interesting cluster algebra structures.
As a representative example, let us consider such a
structure on the algebra $A$ of polynomial
functions on the space of complex $4\times 5$-matrices. For
$1\leq i \leq 4$ and $1\leq j\leq 5$, let $D(ij)$ be the determinant
of the largest square submatrix of a $4\times 5$-matrix whose
upper left corner is the $(i,j)$-coefficient. Then the algebra $A$
admits a cluster structure of type $\tilde{Q}$
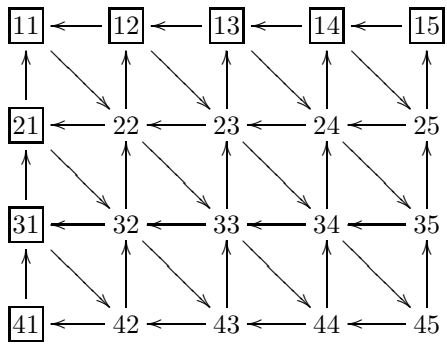
\[
\begin{xy} 0;<0.5pt,0pt>:<0pt,-0.5pt>::
(0,0) *+{\framebox[3ex]{11}} ="0",
(75,0) *+{\framebox[3ex]{12}} ="1",
(150,0) *+{\framebox[3ex]{13}} ="2",
(225,0) *+{\framebox[3ex]{14}} ="3",
(300,0) *+{\framebox[3ex]{15}} ="4",
(0,75) *+{\framebox[3ex]{21}} ="5",
(75,75) *+{22} ="6",
(150,75) *+{23} ="7",
(225,75) *+{24} ="8",
(300,75) *+{25} ="9",
(0,150) *+{\framebox[3ex]{31}} ="10",
(75,150) *+{32} ="11",
(150,150) *+{33} ="12",
(225,150) *+{34} ="13",
(300,150) *+{35} ="14",
(0,225) *+{\framebox[3ex]{41}} ="15",
(75,225) *+{42} ="16",
(150,225) *+{43} ="17",
(225,225) *+{44} ="18",
(300,225) *+{45} ="19",
"1", {\ar"0"},
"5", {\ar"0"},
"0", {\ar"6"},
"2", {\ar"1"},
"6", {\ar"1"},
"1", {\ar"7"},
"3", {\ar"2"},
"7", {\ar"2"},
"2", {\ar"8"},
"4", {\ar"3"},
"8", {\ar"3"},
"3", {\ar"9"},
"9", {\ar"4"},
"6", {\ar"5"},
"10", {\ar"5"},
"5", {\ar"11"},
"7", {\ar"6"},
"11", {\ar"6"},
"6", {\ar"12"},
"8", {\ar"7"},
"12", {\ar"7"},
"7", {\ar"13"},
"9", {\ar"8"},
"13", {\ar"8"},
"8", {\ar"14"},
"14", {\ar"9"},
"11", {\ar"10"},
"15", {\ar"10"},
"10", {\ar"16"},
"12", {\ar"11"},
"16", {\ar"11"},
"11", {\ar"17"},
"13", {\ar"12"},
"17", {\ar"12"},
"12", {\ar"18"},
"14", {\ar"13"},
"18", {\ar"13"},
"13", {\ar"19"},
"19", {\ar"14"},
"16", {\ar"15"},
"17", {\ar"16"},
"18", {\ar"17"},
"19", {\ar"18"},
\end{xy}
\]
whose initial seed is formed by the functions $D(ij)$ associated
with the vertices of the quiver $\tilde{Q}$. This is a particular case
of a theorem of Geiss--Leclerc--Schr\"oer \cite{GeissLeclercSchroeer11b}.
Perhaps the most remarkable fact is that iterated mutations of
the initial seed still produce polynomials in the matrix coefficients (and not fractions).
Geiss--Leclerc--Schr\"oer's proof of this fact in \cite{GeissLeclercSchroeer11b} is ultimately
based on Lusztig's results \cite{Lusztig00}. They sketch a more elementary
approach in section~7.3 of \cite{GeissLeclercSchroeer11a}, cf.~section~\ref{ss:factoriality} below.
It is not obvious either that the cluster variables generate the polynomial ring.
To prove it, we first notice that the variables $x_{25}$, $x_{35}$, $x_{42}$,
$x_{43}$, $x_{44}$, $x_{45}$ already belong to the initial seed.
Now, following \cite{GeissLeclercSchroeer11b}, we
consider the sequence of mutations at the vertices
\[
45,44,43,42,35,25; 34,33,32,24; 23,22; 45,44,43,35; 34,33; 45, 44.
\]
The sequence naturally splits into `hooks', which we have
separated by semicolons. The cluster variables which appear
successively under this sequence of mutations are
\[
x_{34},x_{33}, x_{32}, \ldots\ , x_{24},  \ldots\ , x_{23}, x_{22}, \ldots\  \ko
\]
where we have only indicated those variables associated
with mutations at the vertices of the lower right rim: $25$, $35$,
$42$, $43$, $44$, $45$. So we see that in fact all the
functions $x_{ij}$ are cluster variables.

\subsection{Finite generation} \label{ss:finite-generation}
\index{cluster algebra!finite generation of}
In general, cluster algebras are not finitely generated as algebras.
For example, consider the cluster algebra $\ca_Q$ associated with
the quiver
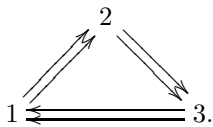
\[
\xymatrix{ & 2 \ar@<0.4ex>[rd] \ar@<-0.4ex>[rd] & \\
1  \ar@<0.4ex>[ru] \ar@<-0.4ex>[ru] & & 3.  \ar@<0.4ex>[ll] \ar@<-0.4ex>[ll]}
\]
Let us show, following \cite{Muller11}, that $\ca_Q$ is not
even Noetherian. Indeed, up to isomorphism, the quiver $Q$ is
invariant under mutations. Hence all exchange relations are of
the form
\[
u_k u_k' = u_i^2 + u_j^2
\]
for three pairwise distinct indices $i$, $j$ and $k$. It follows
that $\ca_Q$ admits a grading such that all cluster variables
have degree $1$. Since $Q$ is not mutation-equivalent to a Dynkin
quiver, by Theorem~\ref{thm:cluster-finite-classification}, there
are infinitely many cluster variables and by Conjecture~\ref{conj:independence},
proved in \cite{CerulliKellerLabardiniPlamondon12}, they are linearly
independent over the field $\Q$, which is the degree $0$ part of $\ca_Q$.
But a positively graded commutative algebra whose degree $1$
part is not a finitely generated module over its degree $0$ part
cannot be Noetherian. Many more examples are provided by the
following theorem

\begin{theorem}[Th.~1.24 of \cite{BerensteinFominZelevinsky05}]
If $Q$ is a valued quiver with three vertices, the cluster algebra $\ca_Q$
is finitely generated over the rationals if and only if $Q$ is
mutation-equivalent to an acyclic valued quiver.
\end{theorem}

For an acyclic valued quiver with $n$ vertices, the cluster algebra
$\ca_Q$ admits a set of $2n$ generators. More precisely, we have
the following theorem.

\begin{theorem}[Cor.~1.21 of \cite{BerensteinFominZelevinsky05}]
If $Q$ is acyclic with $n$ vertices, the cluster algebra $\ca_Q$
is generated over the rationals by the initial variables $x_1$, \ldots, $x_n$
and the cluster variables $x'_j$, $1\leq j\leq n$, obtained by mutating the
initial seed at each vertex~$j$.
\end{theorem}

Moreover, by Cor.~1.21 of \cite{BerensteinFominZelevinsky05}, if
$Q$ is acyclic, the generators $x_1$, \ldots, $x_n$, $x_1'$, \ldots, $x'_n$
together with the exchange relations between $x_j$ and $x_j'$, $1\leq j\leq n$,
form a presentation of $\ca_Q$ and the monomials in the generators not
containing any product $x_j x'_j$ form a $\Q$-basis.

The class of `locally ayclic' cluster algebras is introduced
in \cite{Muller11}. It contains all ayclic cluster algebras.
As shown in \cite{Muller11}, each locally acyclic cluster
algebra is finitely generated, integrally closed and
locally a complete intersection.

\subsection{Factoriality} \label{ss:factoriality}
\index{cluster algebra!factoriality of} In general,
cluster algebras need not be factorial, even when the exchange
matrix is of full rank. The following example, based on
an idea of P.~Lampe, is given in
\cite{GeissLeclercSchroeer11a}.
Let $Q$ be the generalized Kronecker quiver
\[
\xymatrix{ 1 \ar@<1ex>[r] \ar[r] \ar@<-1ex>[r] & 2}
\]
and $x_1'$ the cluster variable obtained by mutating the inital
seed at the vertex $1$. Then we have
\[
x_1 x_1' = 1+x_2^3 = (1+x_2)(1-x_2+x_2^2)
\]
and one can show that these are essentially different factorizations
of the product $x_1 x'_1$ in $\ca_Q$, cf. Prop.~6.3 of \cite{GeissLeclercSchroeer11a}.

Now let $\tilde{Q}$ be a valued ice quiver of type $(n,m)$ and let $n\leq p\leq m$
be an integer. Let $\cp$ be the polynomial ring $\Z[x_{n+1},\ldots, x_m]$ and
$\cl$ its localization at $x_{n+1}$, \ldots, $x_p$.
Let
\[
\ca= \ca_{\tilde{Q}}\ten_\cp \cl
\]
be the localization of the cluster algebra $\ca_{\tilde{Q}}$ at
$x_{n+1}$, \ldots, $x_p$. Notice that the invertible elements of $\cl$ are the
Laurent monomials in $x_{n+1}$, \ldots, $x_p$ multiplied by $\pm 1$.
\begin{theorem}[\cite{GeissLeclercSchroeer11a}]
\begin{itemize}
\item[a)] The invertible elements of $\ca$ are those of $\cl$.
\item[b)] Each cluster variable of $\ca$ is irreducible and two
cluster variables are associate iff they are equal.
\end{itemize}
\end{theorem}
As an application, let us show that the cluster algebra associated
with a Dynkin quiver of type $A_3$ is not factorial. Indeed, consider
the cluster algebra $\ca$ associated with the quiver
\[
Q: \xymatrix{1 \ar[r] & 2 \ar[r] & 3.}
\]
Let $x'_1$ and $x'_3$ be the cluster variables obtained from
the initial seed by mutating respectively at the vertices $1$ and $3$.
We have
\[
x'_1 = \frac{1+x_2}{x_1} \quad\mbox{and}\quad x'_3 = \frac{1+x_2}{x_3}
\]
and therefore
\[
x'_1 x_1 = x'_3 x_3.
\]
Since $x_1$, $x'_1$, $x_3$, $x'_3$ are pairwise distinct cluster
variables, it follows from the theorem
that these are essentially distinct factorizations.

Despite these examples, many cluster algebras appearing `in nature'
are in fact factorial. The following theorem often allows to check
this.

\begin{theorem}[\cite{GeissLeclercSchroeer11a}] As above, let $\ca$ be the cluster
algebra associated with a valued ice quiver of type $(n,m)$ localized at a subset
$x_{n+1}$, \ldots, $x_p$ of the set of coefficients. Let $y$ and $z$
be disjoint clusters and $U\subset \ca$ a subalgebra which is
factorial and contains $y$, $z$ and the localized coefficient
algebra $\cl$. Then $\ca$ equals $U$ and an element $x$ of
the ambient field $\Q(x_1, \ldots, x_m)$ belongs to $\ca$ iff
it is a Laurent polynomial with coefficients in $\cl$
both in $y$ and in $z$.
\end{theorem}

As a prototypical example, consider the ice quiver
\[
\tilde{Q}: \xymatrix{1 \ar[r] & 2 \ar[r] & 3 \ar[r] & *+{\framebox[3ex]{4}}}.
\]
We will parametrise its coefficient and its cluster variables by
the vertices of the following quiver
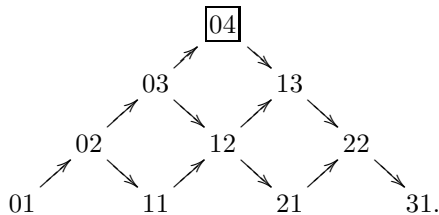
\[
\begin{xy} 0;<0.5pt,0pt>:<0pt,-0.3pt>::
(150,0) *+{\framebox[3ex]{04}} ="0",
(100,75) *+{03} ="1",
(200,75) *+{13} ="2",
(50,150) *+{02} ="3",
(150,150) *+{12} ="4",
(250,150) *+{22} ="5",
(0,225) *+{01} ="6",
(100,225) *+{11} ="7",
(200,225) *+{21} ="8",
(300,225) *+{31.} ="9",
"1", {\ar"0"},
"0", {\ar"2"},
"3", {\ar"1"},
"1", {\ar"4"},
"4", {\ar"2"},
"2", {\ar"5"},
"6", {\ar"3"},
"3", {\ar"7"},
"7", {\ar"4"},
"4", {\ar"8"},
"8", {\ar"5"},
"5", {\ar"9"},
\end{xy}
\]
Namely, to a vertex $ij$, we associate a cluster variable
$x_{i,j}$ in such a way that $x_{0,j}$ equals $x_j$,
$1\leq j \leq 4$, and
each `mesh' gives rise to an exchange relation: We have
\[
x_{i,1}x_{i+1,1}=x_{i,2}+1 \quad \mbox{for } 0\leq i \leq 2
\]
and
\[
x_{i,j}\, x_{i+1,j} = x_{i,j+1} x_{i+1,j-1} +1
\]
for all vertices $ij$ among $02$, $03$, $12$.
Then the set of cluster variables is the set of the
$x_{i,j}$, where $ij$ runs through the vertices other than $04$. The
variables at the bottom are
\[
x_{0,1}=x_1\ko x_{1,1} = \frac{1+x_2}{x_1} \ko x_{2,1}=\frac{x_1+x_3}{x_2} \ko
x_{3,1}=\frac{x_2+x_4}{x_3}
\]
They are algebraically independent and the polynomial ring
\[
U=\Z[x_{0,1}, x_{1,1}, x_{2,1}, x_{3,1}]
\]
contains the disjoint clusters
$y=\{x_1, x_2, x_3\}$ and $z=\{x_{1,3}, x_{2,2}, x_{3,1}\}$ appearing on
the left and the right rim. We see
from the theorem that $U$ equals the cluster algebra
and that an element of the ambient field belongs to the cluster algebra
iff it is a Laurent polynomial with coefficients in $\Z[x_4]$ both
in $y$ and in $z$. We refer to section~7.3 of \cite{GeissLeclercSchroeer11a}
for more elaborate examples arising as coordinate algebras of unipotent cells in
Kac-Moody groups.

\section{General cluster algebras}
\label{s:general-cluster-algebras}

\subsection{Parametrization of seeds by the $n$-regular tree} \label{ss:param-n-reg-tree}
\begin{figure}
\[
\begin{tikzpicture}[scale=2.8]
\draw (0,0) node[right=10pt, above] {$t_0$};
\draw[-] (   0   ,   0   ) node {$\bullet$}--(   1.   ,   0   ) node {$\bullet$} node[midway,fill,fill=white]{{\small    1   }};
\draw[-] (   1.   ,   0   ) node {$\bullet$}--(   .8895   ,   -.3402   ) node {$\bullet$} node[midway,fill,fill=white]{{\small    2   }};
\draw[-] (   .8895   ,   -.3402   ) --(   .7859   ,   -.2650   );
\draw[-] (   .7859   ,   -.2650   ) --(   .8230   ,   -.2381   );
\draw[-] (   .7859   ,   -.2650   ) --(   .7717   ,   -.2214   );
\draw[-] (   .7859   ,   -.2650   ) --(   .7401   ,   -.2650   );
\draw[-] (   .7859   ,   -.2650   ) --(   .7718   ,   -.3085   );
\draw[-] (   .8895   ,   -.3402   ) --(   .7859   ,   -.4155   );
\draw[-] (   .7859   ,   -.4155   ) --(   .7718   ,   -.3720   );
\draw[-] (   .7859   ,   -.4155   ) --(   .7401   ,   -.4155   );
\draw[-] (   .7859   ,   -.4155   ) --(   .7717   ,   -.4591   );
\draw[-] (   .7859   ,   -.4155   ) --(   .8230   ,   -.4424   );
\draw[-] (   .8895   ,   -.3402   ) --(   .9291   ,   -.4620   );
\draw[-] (   .9291   ,   -.4620   ) --(   .8833   ,   -.4621   );
\draw[-] (   .9291   ,   -.4620   ) --(   .9150   ,   -.5055   );
\draw[-] (   .9291   ,   -.4620   ) --(   .9662   ,   -.4889   );
\draw[-] (   .9291   ,   -.4620   ) --(   .9661   ,   -.4351   );
\draw[-] (   .8895   ,   -.3402   ) --(   1.017   ,   -.3402   );
\draw[-] (   1.017   ,   -.3402   ) --(   1.003   ,   -.3838   );
\draw[-] (   1.017   ,   -.3402   ) --(   1.054   ,   -.3671   );
\draw[-] (   1.017   ,   -.3402   ) --(   1.054   ,   -.3133   );
\draw[-] (   1.017   ,   -.3402   ) --(   1.003   ,   -.2966   );
\draw[-] (   1.   ,   0   ) node {$\bullet$}--(   1.289   ,   -.2103   ) node {$\bullet$} node[midway,fill,fill=white]{{\small    3   }};
\draw[-] (   1.289   ,   -.2103   ) --(   1.185   ,   -.2856   );
\draw[-] (   1.185   ,   -.2856   ) --(   1.170   ,   -.2421   );
\draw[-] (   1.185   ,   -.2856   ) --(   1.139   ,   -.2856   );
\draw[-] (   1.185   ,   -.2856   ) --(   1.171   ,   -.3292   );
\draw[-] (   1.185   ,   -.2856   ) --(   1.222   ,   -.3125   );
\draw[-] (   1.289   ,   -.2103   ) --(   1.329   ,   -.3321   );
\draw[-] (   1.329   ,   -.3321   ) --(   1.284   ,   -.3321   );
\draw[-] (   1.329   ,   -.3321   ) --(   1.314   ,   -.3757   );
\draw[-] (   1.329   ,   -.3321   ) --(   1.366   ,   -.3591   );
\draw[-] (   1.329   ,   -.3321   ) --(   1.366   ,   -.3051   );
\draw[-] (   1.289   ,   -.2103   ) --(   1.417   ,   -.2103   );
\draw[-] (   1.417   ,   -.2103   ) --(   1.403   ,   -.2539   );
\draw[-] (   1.417   ,   -.2103   ) --(   1.454   ,   -.2372   );
\draw[-] (   1.417   ,   -.2103   ) --(   1.454   ,   -.1834   );
\draw[-] (   1.417   ,   -.2103   ) --(   1.403   ,   -.1667   );
\draw[-] (   1.289   ,   -.2103   ) --(   1.329   ,   -.08860   );
\draw[-] (   1.329   ,   -.08860   ) --(   1.366   ,   -.1155   );
\draw[-] (   1.329   ,   -.08860   ) --(   1.366   ,   -.06160   );
\draw[-] (   1.329   ,   -.08860   ) --(   1.314   ,   -.04500   );
\draw[-] (   1.329   ,   -.08860   ) --(   1.284   ,   -.08850   );
\draw[-] (   1.   ,   0   ) node {$\bullet$}--(   1.289   ,   .2103   ) node {$\bullet$} node[midway,fill,fill=white]{{\small    4   }};
\draw[-] (   1.289   ,   .2103   ) --(   1.329   ,   .08860   );
\draw[-] (   1.329   ,   .08860   ) --(   1.284   ,   .08850   );
\draw[-] (   1.329   ,   .08860   ) --(   1.314   ,   .04500   );
\draw[-] (   1.329   ,   .08860   ) --(   1.366   ,   .06160   );
\draw[-] (   1.329   ,   .08860   ) --(   1.366   ,   .1155   );
\draw[-] (   1.289   ,   .2103   ) --(   1.417   ,   .2103   );
\draw[-] (   1.417   ,   .2103   ) --(   1.403   ,   .1667   );
\draw[-] (   1.417   ,   .2103   ) --(   1.454   ,   .1834   );
\draw[-] (   1.417   ,   .2103   ) --(   1.454   ,   .2372   );
\draw[-] (   1.417   ,   .2103   ) --(   1.403   ,   .2539   );
\draw[-] (   1.289   ,   .2103   ) --(   1.329   ,   .3321   );
\draw[-] (   1.329   ,   .3321   ) --(   1.366   ,   .3051   );
\draw[-] (   1.329   ,   .3321   ) --(   1.366   ,   .3591   );
\draw[-] (   1.329   ,   .3321   ) --(   1.314   ,   .3757   );
\draw[-] (   1.329   ,   .3321   ) --(   1.284   ,   .3321   );
\draw[-] (   1.289   ,   .2103   ) --(   1.185   ,   .2856   );
\draw[-] (   1.185   ,   .2856   ) --(   1.222   ,   .3125   );
\draw[-] (   1.185   ,   .2856   ) --(   1.171   ,   .3292   );
\draw[-] (   1.185   ,   .2856   ) --(   1.139   ,   .2856   );
\draw[-] (   1.185   ,   .2856   ) --(   1.170   ,   .2421   );
\draw[-] (   1.   ,   0   ) node {$\bullet$}--(   .8895   ,   .3402   ) node {$\bullet$} node[midway,fill,fill=white]{{\small    5   }};
\draw[-] (   .8895   ,   .3402   ) --(   1.017   ,   .3402   );
\draw[-] (   1.017   ,   .3402   ) --(   1.003   ,   .2966   );
\draw[-] (   1.017   ,   .3402   ) --(   1.054   ,   .3133   );
\draw[-] (   1.017   ,   .3402   ) --(   1.054   ,   .3671   );
\draw[-] (   1.017   ,   .3402   ) --(   1.003   ,   .3838   );
\draw[-] (   .8895   ,   .3402   ) --(   .9291   ,   .4620   );
\draw[-] (   .9291   ,   .4620   ) --(   .9661   ,   .4351   );
\draw[-] (   .9291   ,   .4620   ) --(   .9662   ,   .4889   );
\draw[-] (   .9291   ,   .4620   ) --(   .9150   ,   .5055   );
\draw[-] (   .9291   ,   .4620   ) --(   .8833   ,   .4621   );
\draw[-] (   .8895   ,   .3402   ) --(   .7859   ,   .4155   );
\draw[-] (   .7859   ,   .4155   ) --(   .8230   ,   .4424   );
\draw[-] (   .7859   ,   .4155   ) --(   .7717   ,   .4591   );
\draw[-] (   .7859   ,   .4155   ) --(   .7401   ,   .4155   );
\draw[-] (   .7859   ,   .4155   ) --(   .7718   ,   .3720   );
\draw[-] (   .8895   ,   .3402   ) --(   .7859   ,   .2650   );
\draw[-] (   .7859   ,   .2650   ) --(   .7718   ,   .3085   );
\draw[-] (   .7859   ,   .2650   ) --(   .7401   ,   .2650   );
\draw[-] (   .7859   ,   .2650   ) --(   .7717   ,   .2214   );
\draw[-] (   .7859   ,   .2650   ) --(   .8230   ,   .2381   );
\draw[-] (   0   ,   0   ) node {$\bullet$}--(   .3090   ,   .9510   ) node {$\bullet$} node[midway,fill,fill=white]{{\small    2   }};
\draw[-] (   .3090   ,   .9510   ) node {$\bullet$}--(   .5984   ,   .7408   ) node {$\bullet$} node[midway,fill,fill=white]{{\small    3   }};
\draw[-] (   .5984   ,   .7408   ) --(   .4949   ,   .6656   );
\draw[-] (   .4949   ,   .6656   ) --(   .4808   ,   .7091   );
\draw[-] (   .4949   ,   .6656   ) --(   .4491   ,   .6656   );
\draw[-] (   .4949   ,   .6656   ) --(   .4807   ,   .6220   );
\draw[-] (   .4949   ,   .6656   ) --(   .5320   ,   .6387   );
\draw[-] (   .5984   ,   .7408   ) --(   .6380   ,   .6191   );
\draw[-] (   .6380   ,   .6191   ) --(   .5922   ,   .6192   );
\draw[-] (   .6380   ,   .6191   ) --(   .6238   ,   .5756   );
\draw[-] (   .6380   ,   .6191   ) --(   .6750   ,   .5922   );
\draw[-] (   .6380   ,   .6191   ) --(   .6750   ,   .6461   );
\draw[-] (   .5984   ,   .7408   ) --(   .7264   ,   .7409   );
\draw[-] (   .7264   ,   .7409   ) --(   .7123   ,   .6973   );
\draw[-] (   .7264   ,   .7409   ) --(   .7634   ,   .7140   );
\draw[-] (   .7264   ,   .7409   ) --(   .7634   ,   .7678   );
\draw[-] (   .7264   ,   .7409   ) --(   .7123   ,   .7845   );
\draw[-] (   .5984   ,   .7408   ) --(   .6380   ,   .8624   );
\draw[-] (   .6380   ,   .8624   ) --(   .6750   ,   .8356   );
\draw[-] (   .6380   ,   .8624   ) --(   .6751   ,   .8893   );
\draw[-] (   .6380   ,   .8624   ) --(   .6239   ,   .9059   );
\draw[-] (   .6380   ,   .8624   ) --(   .5922   ,   .8625   );
\draw[-] (   .3090   ,   .9510   ) node {$\bullet$}--(   .5984   ,   1.161   ) node {$\bullet$} node[midway,fill,fill=white]{{\small    4   }};
\draw[-] (   .5984   ,   1.161   ) --(   .6379   ,   1.040   );
\draw[-] (   .6379   ,   1.040   ) --(   .5921   ,   1.040   );
\draw[-] (   .6379   ,   1.040   ) --(   .6237   ,   .9967   );
\draw[-] (   .6379   ,   1.040   ) --(   .6749   ,   1.013   );
\draw[-] (   .6379   ,   1.040   ) --(   .6749   ,   1.066   );
\draw[-] (   .5984   ,   1.161   ) --(   .7264   ,   1.161   );
\draw[-] (   .7264   ,   1.161   ) --(   .7123   ,   1.117   );
\draw[-] (   .7264   ,   1.161   ) --(   .7634   ,   1.134   );
\draw[-] (   .7264   ,   1.161   ) --(   .7634   ,   1.188   );
\draw[-] (   .7264   ,   1.161   ) --(   .7123   ,   1.205   );
\draw[-] (   .5984   ,   1.161   ) --(   .6380   ,   1.283   );
\draw[-] (   .6380   ,   1.283   ) --(   .6750   ,   1.257   );
\draw[-] (   .6380   ,   1.283   ) --(   .6751   ,   1.310   );
\draw[-] (   .6380   ,   1.283   ) --(   .6239   ,   1.326   );
\draw[-] (   .6380   ,   1.283   ) --(   .5922   ,   1.283   );
\draw[-] (   .5984   ,   1.161   ) --(   .4949   ,   1.236   );
\draw[-] (   .4949   ,   1.236   ) --(   .5319   ,   1.263   );
\draw[-] (   .4949   ,   1.236   ) --(   .4808   ,   1.280   );
\draw[-] (   .4949   ,   1.236   ) --(   .4492   ,   1.236   );
\draw[-] (   .4949   ,   1.236   ) --(   .4807   ,   1.193   );
\draw[-] (   .3090   ,   .9510   ) node {$\bullet$}--(   .1984   ,   1.291   ) node {$\bullet$} node[midway,fill,fill=white]{{\small    5   }};
\draw[-] (   .1984   ,   1.291   ) --(   .3264   ,   1.291   );
\draw[-] (   .3264   ,   1.291   ) --(   .3123   ,   1.247   );
\draw[-] (   .3264   ,   1.291   ) --(   .3634   ,   1.264   );
\draw[-] (   .3264   ,   1.291   ) --(   .3634   ,   1.318   );
\draw[-] (   .3264   ,   1.291   ) --(   .3123   ,   1.335   );
\draw[-] (   .1984   ,   1.291   ) --(   .2380   ,   1.412   );
\draw[-] (   .2380   ,   1.412   ) --(   .2750   ,   1.386   );
\draw[-] (   .2380   ,   1.412   ) --(   .2750   ,   1.439   );
\draw[-] (   .2380   ,   1.412   ) --(   .2238   ,   1.455   );
\draw[-] (   .2380   ,   1.412   ) --(   .1922   ,   1.412   );
\draw[-] (   .1984   ,   1.291   ) --(   .09480   ,   1.366   );
\draw[-] (   .09480   ,   1.366   ) --(   .1319   ,   1.393   );
\draw[-] (   .09480   ,   1.366   ) --(   .08070   ,   1.410   );
\draw[-] (   .09480   ,   1.366   ) --(   .04910   ,   1.366   );
\draw[-] (   .09480   ,   1.366   ) --(   .08070   ,   1.323   );
\draw[-] (   .1984   ,   1.291   ) --(   .09480   ,   1.215   );
\draw[-] (   .09480   ,   1.215   ) --(   .08070   ,   1.258   );
\draw[-] (   .09480   ,   1.215   ) --(   .04900   ,   1.215   );
\draw[-] (   .09480   ,   1.215   ) --(   .08070   ,   1.171   );
\draw[-] (   .09480   ,   1.215   ) --(   .1319   ,   1.188   );
\draw[-] (   .3090   ,   .9510   ) node {$\bullet$}--(   -.04880   ,   .9510   ) node {$\bullet$} node[midway,fill,fill=white]{{\small    1   }};
\draw[-] (   -.04880   ,   .9510   ) --(   -.009300   ,   1.073   );
\draw[-] (   -.009300   ,   1.073   ) --(   .02780   ,   1.046   );
\draw[-] (   -.009300   ,   1.073   ) --(   .02780   ,   1.100   );
\draw[-] (   -.009300   ,   1.073   ) --(   -.02340   ,   1.116   );
\draw[-] (   -.009300   ,   1.073   ) --(   -.05500   ,   1.073   );
\draw[-] (   -.04880   ,   .9510   ) --(   -.1524   ,   1.026   );
\draw[-] (   -.1524   ,   1.026   ) --(   -.1154   ,   1.053   );
\draw[-] (   -.1524   ,   1.026   ) --(   -.1666   ,   1.070   );
\draw[-] (   -.1524   ,   1.026   ) --(   -.1982   ,   1.026   );
\draw[-] (   -.1524   ,   1.026   ) --(   -.1666   ,   .9828   );
\draw[-] (   -.04880   ,   .9510   ) --(   -.1524   ,   .8758   );
\draw[-] (   -.1524   ,   .8758   ) --(   -.1666   ,   .9193   );
\draw[-] (   -.1524   ,   .8758   ) --(   -.1982   ,   .8758   );
\draw[-] (   -.1524   ,   .8758   ) --(   -.1666   ,   .8322   );
\draw[-] (   -.1524   ,   .8758   ) --(   -.1154   ,   .8489   );
\draw[-] (   -.04880   ,   .9510   ) --(   -.009300   ,   .8293   );
\draw[-] (   -.009300   ,   .8293   ) --(   -.05500   ,   .8293   );
\draw[-] (   -.009300   ,   .8293   ) --(   -.02340   ,   .7858   );
\draw[-] (   -.009300   ,   .8293   ) --(   .02780   ,   .8024   );
\draw[-] (   -.009300   ,   .8293   ) --(   .02780   ,   .8562   );
\draw[-] (   0   ,   0   ) node {$\bullet$}--(   -.8090   ,   .5878   ) node {$\bullet$} node[midway,fill,fill=white]{{\small    3   }};
\draw[-] (   -.8090   ,   .5878   ) node {$\bullet$}--(   -.5196   ,   .7980   ) node {$\bullet$} node[midway,fill,fill=white]{{\small    4   }};
\draw[-] (   -.5196   ,   .7980   ) --(   -.4801   ,   .6764   );
\draw[-] (   -.4801   ,   .6764   ) --(   -.5259   ,   .6765   );
\draw[-] (   -.4801   ,   .6764   ) --(   -.4943   ,   .6329   );
\draw[-] (   -.4801   ,   .6764   ) --(   -.4431   ,   .6495   );
\draw[-] (   -.4801   ,   .6764   ) --(   -.4431   ,   .7034   );
\draw[-] (   -.5196   ,   .7980   ) --(   -.3916   ,   .7979   );
\draw[-] (   -.3916   ,   .7979   ) --(   -.4058   ,   .7543   );
\draw[-] (   -.3916   ,   .7979   ) --(   -.3546   ,   .7710   );
\draw[-] (   -.3916   ,   .7979   ) --(   -.3546   ,   .8248   );
\draw[-] (   -.3916   ,   .7979   ) --(   -.4058   ,   .8415   );
\draw[-] (   -.5196   ,   .7980   ) --(   -.4800   ,   .9197   );
\draw[-] (   -.4800   ,   .9197   ) --(   -.4430   ,   .8928   );
\draw[-] (   -.4800   ,   .9197   ) --(   -.4430   ,   .9466   );
\draw[-] (   -.4800   ,   .9197   ) --(   -.4942   ,   .9632   );
\draw[-] (   -.4800   ,   .9197   ) --(   -.5258   ,   .9198   );
\draw[-] (   -.5196   ,   .7980   ) --(   -.6231   ,   .8732   );
\draw[-] (   -.6231   ,   .8732   ) --(   -.5861   ,   .9001   );
\draw[-] (   -.6231   ,   .8732   ) --(   -.6373   ,   .9168   );
\draw[-] (   -.6231   ,   .8732   ) --(   -.6689   ,   .8732   );
\draw[-] (   -.6231   ,   .8732   ) --(   -.6373   ,   .8297   );
\draw[-] (   -.8090   ,   .5878   ) node {$\bullet$}--(   -.9194   ,   .9280   ) node {$\bullet$} node[midway,fill,fill=white]{{\small    5   }};
\draw[-] (   -.9194   ,   .9280   ) --(   -.7914   ,   .9280   );
\draw[-] (   -.7914   ,   .9280   ) --(   -.8056   ,   .8844   );
\draw[-] (   -.7914   ,   .9280   ) --(   -.7544   ,   .9011   );
\draw[-] (   -.7914   ,   .9280   ) --(   -.7544   ,   .9549   );
\draw[-] (   -.7914   ,   .9280   ) --(   -.8056   ,   .9716   );
\draw[-] (   -.9194   ,   .9280   ) --(   -.8799   ,   1.050   );
\draw[-] (   -.8799   ,   1.050   ) --(   -.8429   ,   1.024   );
\draw[-] (   -.8799   ,   1.050   ) --(   -.8429   ,   1.077   );
\draw[-] (   -.8799   ,   1.050   ) --(   -.8941   ,   1.093   );
\draw[-] (   -.8799   ,   1.050   ) --(   -.9257   ,   1.050   );
\draw[-] (   -.9194   ,   .9280   ) --(   -1.023   ,   1.004   );
\draw[-] (   -1.023   ,   1.004   ) --(   -.9864   ,   1.031   );
\draw[-] (   -1.023   ,   1.004   ) --(   -1.037   ,   1.048   );
\draw[-] (   -1.023   ,   1.004   ) --(   -1.069   ,   1.004   );
\draw[-] (   -1.023   ,   1.004   ) --(   -1.038   ,   .9605   );
\draw[-] (   -.9194   ,   .9280   ) --(   -1.023   ,   .8528   );
\draw[-] (   -1.023   ,   .8528   ) --(   -1.038   ,   .8963   );
\draw[-] (   -1.023   ,   .8528   ) --(   -1.069   ,   .8528   );
\draw[-] (   -1.023   ,   .8528   ) --(   -1.037   ,   .8092   );
\draw[-] (   -1.023   ,   .8528   ) --(   -.9864   ,   .8259   );
\draw[-] (   -.8090   ,   .5878   ) node {$\bullet$}--(   -1.167   ,   .5878   ) node {$\bullet$} node[midway,fill,fill=white]{{\small    1   }};
\draw[-] (   -1.167   ,   .5878   ) --(   -1.127   ,   .7095   );
\draw[-] (   -1.127   ,   .7095   ) --(   -1.090   ,   .6826   );
\draw[-] (   -1.127   ,   .7095   ) --(   -1.090   ,   .7364   );
\draw[-] (   -1.127   ,   .7095   ) --(   -1.142   ,   .7530   );
\draw[-] (   -1.127   ,   .7095   ) --(   -1.172   ,   .7095   );
\draw[-] (   -1.167   ,   .5878   ) --(   -1.271   ,   .6630   );
\draw[-] (   -1.271   ,   .6630   ) --(   -1.234   ,   .6899   );
\draw[-] (   -1.271   ,   .6630   ) --(   -1.285   ,   .7066   );
\draw[-] (   -1.271   ,   .6630   ) --(   -1.317   ,   .6630   );
\draw[-] (   -1.271   ,   .6630   ) --(   -1.286   ,   .6195   );
\draw[-] (   -1.167   ,   .5878   ) --(   -1.271   ,   .5126   );
\draw[-] (   -1.271   ,   .5126   ) --(   -1.286   ,   .5561   );
\draw[-] (   -1.271   ,   .5126   ) --(   -1.317   ,   .5126   );
\draw[-] (   -1.271   ,   .5126   ) --(   -1.285   ,   .4691   );
\draw[-] (   -1.271   ,   .5126   ) --(   -1.234   ,   .4857   );
\draw[-] (   -1.167   ,   .5878   ) --(   -1.127   ,   .4661   );
\draw[-] (   -1.127   ,   .4661   ) --(   -1.172   ,   .4661   );
\draw[-] (   -1.127   ,   .4661   ) --(   -1.142   ,   .4226   );
\draw[-] (   -1.127   ,   .4661   ) --(   -1.090   ,   .4392   );
\draw[-] (   -1.127   ,   .4661   ) --(   -1.090   ,   .4930   );
\draw[-] (   -.8090   ,   .5878   ) node {$\bullet$}--(   -.9196   ,   .2476   ) node {$\bullet$} node[midway,fill,fill=white]{{\small    2   }};
\draw[-] (   -.9196   ,   .2476   ) --(   -1.023   ,   .3228   );
\draw[-] (   -1.023   ,   .3228   ) --(   -.9864   ,   .3497   );
\draw[-] (   -1.023   ,   .3228   ) --(   -1.037   ,   .3664   );
\draw[-] (   -1.023   ,   .3228   ) --(   -1.069   ,   .3228   );
\draw[-] (   -1.023   ,   .3228   ) --(   -1.038   ,   .2793   );
\draw[-] (   -.9196   ,   .2476   ) --(   -1.023   ,   .1724   );
\draw[-] (   -1.023   ,   .1724   ) --(   -1.038   ,   .2159   );
\draw[-] (   -1.023   ,   .1724   ) --(   -1.069   ,   .1724   );
\draw[-] (   -1.023   ,   .1724   ) --(   -1.037   ,   .1288   );
\draw[-] (   -1.023   ,   .1724   ) --(   -.9864   ,   .1455   );
\draw[-] (   -.9196   ,   .2476   ) --(   -.8800   ,   .1258   );
\draw[-] (   -.8800   ,   .1258   ) --(   -.9258   ,   .1259   );
\draw[-] (   -.8800   ,   .1258   ) --(   -.8942   ,   .08230   );
\draw[-] (   -.8800   ,   .1258   ) --(   -.8430   ,   .09890   );
\draw[-] (   -.8800   ,   .1258   ) --(   -.8430   ,   .1528   );
\draw[-] (   -.9196   ,   .2476   ) --(   -.7916   ,   .2476   );
\draw[-] (   -.7916   ,   .2476   ) --(   -.8058   ,   .2040   );
\draw[-] (   -.7916   ,   .2476   ) --(   -.7546   ,   .2207   );
\draw[-] (   -.7916   ,   .2476   ) --(   -.7546   ,   .2745   );
\draw[-] (   -.7916   ,   .2476   ) --(   -.8058   ,   .2912   );
\draw[-] (   0   ,   0   ) node {$\bullet$}--(   -.8090   ,   -.5878   ) node {$\bullet$} node[midway,fill,fill=white]{{\small    4   }};
\draw[-] (   -.8090   ,   -.5878   ) node {$\bullet$}--(   -.9196   ,   -.2476   ) node {$\bullet$} node[midway,fill,fill=white]{{\small    5   }};
\draw[-] (   -.9196   ,   -.2476   ) --(   -.7916   ,   -.2476   );
\draw[-] (   -.7916   ,   -.2476   ) --(   -.8058   ,   -.2912   );
\draw[-] (   -.7916   ,   -.2476   ) --(   -.7546   ,   -.2745   );
\draw[-] (   -.7916   ,   -.2476   ) --(   -.7546   ,   -.2207   );
\draw[-] (   -.7916   ,   -.2476   ) --(   -.8058   ,   -.2040   );
\draw[-] (   -.9196   ,   -.2476   ) --(   -.8800   ,   -.1258   );
\draw[-] (   -.8800   ,   -.1258   ) --(   -.8430   ,   -.1528   );
\draw[-] (   -.8800   ,   -.1258   ) --(   -.8430   ,   -.09890   );
\draw[-] (   -.8800   ,   -.1258   ) --(   -.8942   ,   -.08230   );
\draw[-] (   -.8800   ,   -.1258   ) --(   -.9258   ,   -.1259   );
\draw[-] (   -.9196   ,   -.2476   ) --(   -1.023   ,   -.1724   );
\draw[-] (   -1.023   ,   -.1724   ) --(   -.9864   ,   -.1455   );
\draw[-] (   -1.023   ,   -.1724   ) --(   -1.037   ,   -.1288   );
\draw[-] (   -1.023   ,   -.1724   ) --(   -1.069   ,   -.1724   );
\draw[-] (   -1.023   ,   -.1724   ) --(   -1.038   ,   -.2159   );
\draw[-] (   -.9196   ,   -.2476   ) --(   -1.023   ,   -.3228   );
\draw[-] (   -1.023   ,   -.3228   ) --(   -1.038   ,   -.2793   );
\draw[-] (   -1.023   ,   -.3228   ) --(   -1.069   ,   -.3228   );
\draw[-] (   -1.023   ,   -.3228   ) --(   -1.037   ,   -.3664   );
\draw[-] (   -1.023   ,   -.3228   ) --(   -.9864   ,   -.3497   );
\draw[-] (   -.8090   ,   -.5878   ) node {$\bullet$}--(   -1.167   ,   -.5878   ) node {$\bullet$} node[midway,fill,fill=white]{{\small    1   }};
\draw[-] (   -1.167   ,   -.5878   ) --(   -1.127   ,   -.4661   );
\draw[-] (   -1.127   ,   -.4661   ) --(   -1.090   ,   -.4930   );
\draw[-] (   -1.127   ,   -.4661   ) --(   -1.090   ,   -.4392   );
\draw[-] (   -1.127   ,   -.4661   ) --(   -1.142   ,   -.4226   );
\draw[-] (   -1.127   ,   -.4661   ) --(   -1.172   ,   -.4661   );
\draw[-] (   -1.167   ,   -.5878   ) --(   -1.271   ,   -.5126   );
\draw[-] (   -1.271   ,   -.5126   ) --(   -1.234   ,   -.4857   );
\draw[-] (   -1.271   ,   -.5126   ) --(   -1.285   ,   -.4691   );
\draw[-] (   -1.271   ,   -.5126   ) --(   -1.317   ,   -.5126   );
\draw[-] (   -1.271   ,   -.5126   ) --(   -1.286   ,   -.5561   );
\draw[-] (   -1.167   ,   -.5878   ) --(   -1.271   ,   -.6630   );
\draw[-] (   -1.271   ,   -.6630   ) --(   -1.286   ,   -.6195   );
\draw[-] (   -1.271   ,   -.6630   ) --(   -1.317   ,   -.6630   );
\draw[-] (   -1.271   ,   -.6630   ) --(   -1.285   ,   -.7066   );
\draw[-] (   -1.271   ,   -.6630   ) --(   -1.234   ,   -.6899   );
\draw[-] (   -1.167   ,   -.5878   ) --(   -1.127   ,   -.7095   );
\draw[-] (   -1.127   ,   -.7095   ) --(   -1.172   ,   -.7095   );
\draw[-] (   -1.127   ,   -.7095   ) --(   -1.142   ,   -.7530   );
\draw[-] (   -1.127   ,   -.7095   ) --(   -1.090   ,   -.7364   );
\draw[-] (   -1.127   ,   -.7095   ) --(   -1.090   ,   -.6826   );
\draw[-] (   -.8090   ,   -.5878   ) node {$\bullet$}--(   -.9194   ,   -.9280   ) node {$\bullet$} node[midway,fill,fill=white]{{\small    2   }};
\draw[-] (   -.9194   ,   -.9280   ) --(   -1.023   ,   -.8528   );
\draw[-] (   -1.023   ,   -.8528   ) --(   -.9864   ,   -.8259   );
\draw[-] (   -1.023   ,   -.8528   ) --(   -1.037   ,   -.8092   );
\draw[-] (   -1.023   ,   -.8528   ) --(   -1.069   ,   -.8528   );
\draw[-] (   -1.023   ,   -.8528   ) --(   -1.038   ,   -.8963   );
\draw[-] (   -.9194   ,   -.9280   ) --(   -1.023   ,   -1.004   );
\draw[-] (   -1.023   ,   -1.004   ) --(   -1.038   ,   -.9605   );
\draw[-] (   -1.023   ,   -1.004   ) --(   -1.069   ,   -1.004   );
\draw[-] (   -1.023   ,   -1.004   ) --(   -1.037   ,   -1.048   );
\draw[-] (   -1.023   ,   -1.004   ) --(   -.9864   ,   -1.031   );
\draw[-] (   -.9194   ,   -.9280   ) --(   -.8799   ,   -1.050   );
\draw[-] (   -.8799   ,   -1.050   ) --(   -.9257   ,   -1.050   );
\draw[-] (   -.8799   ,   -1.050   ) --(   -.8941   ,   -1.093   );
\draw[-] (   -.8799   ,   -1.050   ) --(   -.8429   ,   -1.077   );
\draw[-] (   -.8799   ,   -1.050   ) --(   -.8429   ,   -1.024   );
\draw[-] (   -.9194   ,   -.9280   ) --(   -.7914   ,   -.9280   );
\draw[-] (   -.7914   ,   -.9280   ) --(   -.8056   ,   -.9716   );
\draw[-] (   -.7914   ,   -.9280   ) --(   -.7544   ,   -.9549   );
\draw[-] (   -.7914   ,   -.9280   ) --(   -.7544   ,   -.9011   );
\draw[-] (   -.7914   ,   -.9280   ) --(   -.8056   ,   -.8844   );
\draw[-] (   -.8090   ,   -.5878   ) node {$\bullet$}--(   -.5196   ,   -.7980   ) node {$\bullet$} node[midway,fill,fill=white]{{\small    3   }};
\draw[-] (   -.5196   ,   -.7980   ) --(   -.6231   ,   -.8732   );
\draw[-] (   -.6231   ,   -.8732   ) --(   -.6373   ,   -.8297   );
\draw[-] (   -.6231   ,   -.8732   ) --(   -.6689   ,   -.8732   );
\draw[-] (   -.6231   ,   -.8732   ) --(   -.6373   ,   -.9168   );
\draw[-] (   -.6231   ,   -.8732   ) --(   -.5861   ,   -.9001   );
\draw[-] (   -.5196   ,   -.7980   ) --(   -.4800   ,   -.9197   );
\draw[-] (   -.4800   ,   -.9197   ) --(   -.5258   ,   -.9198   );
\draw[-] (   -.4800   ,   -.9197   ) --(   -.4942   ,   -.9632   );
\draw[-] (   -.4800   ,   -.9197   ) --(   -.4430   ,   -.9466   );
\draw[-] (   -.4800   ,   -.9197   ) --(   -.4430   ,   -.8928   );
\draw[-] (   -.5196   ,   -.7980   ) --(   -.3916   ,   -.7979   );
\draw[-] (   -.3916   ,   -.7979   ) --(   -.4058   ,   -.8415   );
\draw[-] (   -.3916   ,   -.7979   ) --(   -.3546   ,   -.8248   );
\draw[-] (   -.3916   ,   -.7979   ) --(   -.3546   ,   -.7710   );
\draw[-] (   -.3916   ,   -.7979   ) --(   -.4058   ,   -.7543   );
\draw[-] (   -.5196   ,   -.7980   ) --(   -.4801   ,   -.6764   );
\draw[-] (   -.4801   ,   -.6764   ) --(   -.4431   ,   -.7034   );
\draw[-] (   -.4801   ,   -.6764   ) --(   -.4431   ,   -.6495   );
\draw[-] (   -.4801   ,   -.6764   ) --(   -.4943   ,   -.6329   );
\draw[-] (   -.4801   ,   -.6764   ) --(   -.5259   ,   -.6765   );
\draw[-] (   0   ,   0   ) node {$\bullet$}--(   .3090   ,   -.9510   ) node {$\bullet$} node[midway,fill,fill=white]{{\small    5   }};
\draw[-] (   .3090   ,   -.9510   ) node {$\bullet$}--(   -.04880   ,   -.9510   ) node {$\bullet$} node[midway,fill,fill=white]{{\small    1   }};
\draw[-] (   -.04880   ,   -.9510   ) --(   -.009300   ,   -.8293   );
\draw[-] (   -.009300   ,   -.8293   ) --(   .02780   ,   -.8562   );
\draw[-] (   -.009300   ,   -.8293   ) --(   .02780   ,   -.8024   );
\draw[-] (   -.009300   ,   -.8293   ) --(   -.02340   ,   -.7858   );
\draw[-] (   -.009300   ,   -.8293   ) --(   -.05500   ,   -.8293   );
\draw[-] (   -.04880   ,   -.9510   ) --(   -.1524   ,   -.8758   );
\draw[-] (   -.1524   ,   -.8758   ) --(   -.1154   ,   -.8489   );
\draw[-] (   -.1524   ,   -.8758   ) --(   -.1666   ,   -.8322   );
\draw[-] (   -.1524   ,   -.8758   ) --(   -.1982   ,   -.8758   );
\draw[-] (   -.1524   ,   -.8758   ) --(   -.1666   ,   -.9193   );
\draw[-] (   -.04880   ,   -.9510   ) --(   -.1524   ,   -1.026   );
\draw[-] (   -.1524   ,   -1.026   ) --(   -.1666   ,   -.9825   );
\draw[-] (   -.1524   ,   -1.026   ) --(   -.1982   ,   -1.026   );
\draw[-] (   -.1524   ,   -1.026   ) --(   -.1666   ,   -1.070   );
\draw[-] (   -.1524   ,   -1.026   ) --(   -.1154   ,   -1.053   );
\draw[-] (   -.04880   ,   -.9510   ) --(   -.009300   ,   -1.073   );
\draw[-] (   -.009300   ,   -1.073   ) --(   -.05500   ,   -1.073   );
\draw[-] (   -.009300   ,   -1.073   ) --(   -.02340   ,   -1.116   );
\draw[-] (   -.009300   ,   -1.073   ) --(   .02780   ,   -1.100   );
\draw[-] (   -.009300   ,   -1.073   ) --(   .02780   ,   -1.046   );
\draw[-] (   .3090   ,   -.9510   ) node {$\bullet$}--(   .1984   ,   -1.291   ) node {$\bullet$} node[midway,fill,fill=white]{{\small    2   }};
\draw[-] (   .1984   ,   -1.291   ) --(   .09480   ,   -1.215   );
\draw[-] (   .09480   ,   -1.215   ) --(   .1319   ,   -1.188   );
\draw[-] (   .09480   ,   -1.215   ) --(   .08070   ,   -1.171   );
\draw[-] (   .09480   ,   -1.215   ) --(   .04900   ,   -1.215   );
\draw[-] (   .09480   ,   -1.215   ) --(   .08070   ,   -1.258   );
\draw[-] (   .1984   ,   -1.291   ) --(   .09480   ,   -1.366   );
\draw[-] (   .09480   ,   -1.366   ) --(   .08070   ,   -1.323   );
\draw[-] (   .09480   ,   -1.366   ) --(   .04910   ,   -1.366   );
\draw[-] (   .09480   ,   -1.366   ) --(   .08070   ,   -1.410   );
\draw[-] (   .09480   ,   -1.366   ) --(   .01319   ,   -1.393   );
\draw[-] (   .1984   ,   -1.291   ) --(   .2380   ,   -1.412   );
\draw[-] (   .2380   ,   -1.412   ) --(   .1922   ,   -1.412   );
\draw[-] (   .2380   ,   -1.412   ) --(   .2238   ,   -1.455   );
\draw[-] (   .2380   ,   -1.412   ) --(   .2750   ,   -1.439   );
\draw[-] (   .2380   ,   -1.412   ) --(   .2750   ,   -1.386   );
\draw[-] (   .1984   ,   -1.291   ) --(   .3264   ,   -1.291   );
\draw[-] (   .3264   ,   -1.291   ) --(   .3123   ,   -1.335   );
\draw[-] (   .3264   ,   -1.291   ) --(   .3634   ,   -1.318   );
\draw[-] (   .3264   ,   -1.291   ) --(   .3634   ,   -1.264   );
\draw[-] (   .3264   ,   -1.291   ) --(   .3123   ,   -1.247   );
\draw[-] (   .3090   ,   -.9510   ) node {$\bullet$}--(   .5984   ,   -1.161   ) node {$\bullet$} node[midway,fill,fill=white]{{\small    3   }};
\draw[-] (   .5984   ,   -1.161   ) --(   .4949   ,   -1.236   );
\draw[-] (   .4949   ,   -1.236   ) --(   .4807   ,   -1.193   );
\draw[-] (   .4949   ,   -1.236   ) --(   .4492   ,   -1.236   );
\draw[-] (   .4949   ,   -1.236   ) --(   .4808   ,   -1.280   );
\draw[-] (   .4949   ,   -1.236   ) --(   .5319   ,   -1.263   );
\draw[-] (   .5984   ,   -1.161   ) --(   .6380   ,   -1.283   );
\draw[-] (   .6380   ,   -1.283   ) --(   .5922   ,   -1.283   );
\draw[-] (   .6380   ,   -1.283   ) --(   .6239   ,   -1.326   );
\draw[-] (   .6380   ,   -1.283   ) --(   .6751   ,   -1.310   );
\draw[-] (   .6380   ,   -1.283   ) --(   .6750   ,   -1.257   );
\draw[-] (   .5984   ,   -1.161   ) --(   .7264   ,   -1.161   );
\draw[-] (   .7264   ,   -1.161   ) --(   .7123   ,   -1.205   );
\draw[-] (   .7264   ,   -1.161   ) --(   .7634   ,   -1.188   );
\draw[-] (   .7264   ,   -1.161   ) --(   .7634   ,   -1.134   );
\draw[-] (   .7264   ,   -1.161   ) --(   .7123   ,   -1.117   );
\draw[-] (   .5984   ,   -1.161   ) --(   .6379   ,   -1.040   );
\draw[-] (   .6379   ,   -1.040   ) --(   .6749   ,   -1.066   );
\draw[-] (   .6379   ,   -1.040   ) --(   .6749   ,   -1.013   );
\draw[-] (   .6379   ,   -1.040   ) --(   .6237   ,   -.9967   );
\draw[-] (   .6379   ,   -1.040   ) --(   .5921   ,   -1.040   );
\draw[-] (   .3090   ,   -.9510   ) node {$\bullet$}--(   .5984   ,   -.7408   ) node {$\bullet$} node[midway,fill,fill=white]{{\small    4   }};
\draw[-] (   .5984   ,   -.7408   ) --(   .6380   ,   -.8624   );
\draw[-] (   .6380   ,   -.8624   ) --(   .5922   ,   -.8625   );
\draw[-] (   .6380   ,   -.8624   ) --(   .6239   ,   -.9059   );
\draw[-] (   .6380   ,   -.8624   ) --(   .6751   ,   -.8893   );
\draw[-] (   .6380   ,   -.8624   ) --(   .6750   ,   -.8356   );
\draw[-] (   .5984   ,   -.7408   ) --(   .7264   ,   -.7409   );
\draw[-] (   .7264   ,   -.7409   ) --(   .7123   ,   -.7845   );
\draw[-] (   .7264   ,   -.7409   ) --(   .7634   ,   -.7678   );
\draw[-] (   .7264   ,   -.7409   ) --(   .7634   ,   -.7140   );
\draw[-] (   .7264   ,   -.7409   ) --(   .7123   ,   -.6973   );
\draw[-] (   .5984   ,   -.7408   ) --(   .6380   ,   -.6191   );
\draw[-] (   .6380   ,   -.6191   ) --(   .6750   ,   -.6461   );
\draw[-] (   .6380   ,   -.6191   ) --(   .6750   ,   -.5922   );
\draw[-] (   .6380   ,   -.6191   ) --(   .6238   ,   -.5756   );
\draw[-] (   .6380   ,   -.6191   ) --(   .5922   ,   -.6192   );
\draw[-] (   .5984   ,   -.7408   ) --(   .4949   ,   -.6656   );
\draw[-] (   .4949   ,   -.6656   ) --(   .5320   ,   -.6387   );
\draw[-] (   .4949   ,   -.6656   ) --(   .4807   ,   -.6220   );
\draw[-] (   .4949   ,   -.6656   ) --(   .4491   ,   -.6656   );
\draw[-] (   .4949   ,   -.6656   ) --(   .4808   ,   -.7091   );
\end{tikzpicture}
\]
\caption{A picture, up to depth $4$, of the $5$-regular tree} \label{fig:n-regular-tree}
\end{figure}
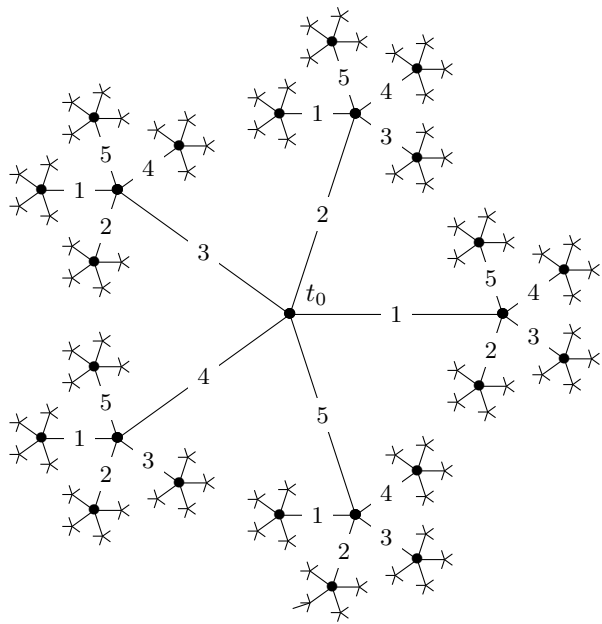

Let us introduce a convenient parametrization of the seeds in the
mutation class of a given initial seed.
Let $1\leq n\leq m$ be integers and $\tilde{Q}$ a valued ice quiver of type $(n,m)$.
Let $X=\{x_1, \ldots, x_m\}$ be the initial cluster and
$(\tilde{Q}, X)$ the initial seed.
Let $\mathbb{T}_n$ be the {\em $n$-regular tree:}\index{tree!$n$-regular} Its edges
are labeled by the integers $1$, \ldots, $n$ such that
the $n$ edges emanating from each vertex carry different
labels, cf. figure~\ref{fig:n-regular-tree}. Let $t_0$ be a vertex of $\mathbb{T}_n$.
To each vertex $t$ of $\mathbb{T}_n$ we associate
a seed $(\tilde{Q}(t), X(t))$
such that at $t=t_0$, we have the initial seed
and whenever $t$ is linked to $t'$ by an edge labeled $k$,
the seeds associated with $t$ and $t'$ are related by
the mutation at $k$. We write $x_i(t)$, $1\leq i\leq n$, for the cluster
variables in the seed $X(t)$. If $\tilde{B}$ is the
$m\times n$-matrix associated with $\tilde{Q}$, we write
$\tilde{B}(t)$ for the matrix associated with $\tilde{Q}(t)$.

\subsection{Principal coefficients} \label{ss:principal-coefficients}
Let $n\geq 1$ be an integer and $Q$ a valued quiver with $n$ vertices.
Let $B$ be the associated skew-symmetri\-zable integer $n\times n$-matrix.
In the next subsections, following \cite{FominZelevinsky07},
we will define data associated with $Q$ which are relevant for all cluster
algebras with coefficients associated with valued ice quivers
$\tilde{Q}$ whose principal part is $Q$. This will become
apparent from a general formula expressing the cluster
variables in terms of these data, cf.~section~\ref{ss:cluster-alg-semifield}.

\subsection{Principal coefficients: $c$-vectors}
\label{ss:principal-coefficients-c-vectors}
Let $Q_{pr}$ be the {\em principal extension\index{principal extension!of a quiver} of $Q$}, i.e.
the valued quiver obtained from $Q$ by adding new vertices
$n+1$, \ldots, $2n$ and new arrows $i+n \to i$, $1\leq i\leq n$,
for each vertex $i$ of $Q$. For example, if we have
\[
Q: \xymatrix{1 \ar[r] & 2}\ko \quad\mbox{then} \quad
Q_{pr}: \quad\raisebox{4ex}{\xymatrix{ 1' \ar[d] & 2' \ar[d] \\ 1 \ar[r] & 2,}}
\]
where we write $i'$ for $i+n$.
The {\em cluster algebra with principal coefficients\index{cluster algebra!with principal coefficients} associated with $Q$}
is the cluster algebra associated with $Q_{pr}$. We write $B_{pr}$
for the corresponding integer $2n\times n$-matrix. It is obtained
from $B$ by appending an $n\times n$ identity matrix at the bottom:
\[
B_{pr} = \left[ \begin{array}{cc} B \\ I_n \end{array} \right].
\]

For a vertex $t$ of the $n$-regular tree, the {\em matrix
of $c$-vectors $C(t)$}\index{$c$-vector} is by definition the $n\times n$-matrix
appearing in the bottom part of $B_{pr}(t)$, so that we have
\[
B_{pr}(t) = \left[ \begin{array}{cc} B(t) \\ C(t) \end{array} \right].
\]
Its columns are the {$c$-vectors at $t$}.
When necessary, we will denote the matrix $C(t)$ by $C(B,t_0,t)$ to clarify
its dependence on $B$ and the sequence of mutations linking $t_0$ to $t$.
For example, if we successively mutate the quiver $Q_{pr}$ associated
with $Q: 1\to 2$ at the vertices $1$, $2$, $1$, \ldots, we obtain
the sequence
\begin{equation} \label{eq:A2-mut-seq}
\xymatrix@R=0.4cm@C=0.4cm{1' \ar[d] & 2' \ar[d] \\ *+<6pt>[o][F]{1} \ar[r] & 2} \raisebox{-3ex}{$\longmapsto$}
\xymatrix@R=0.4cm@C=0.4cm{1' \ar[rd] & 2' \ar[d] \\ 1 \ar[u] & *+<6pt>[o][F]{2}\ar[l]} \raisebox{-3ex}{$\longmapsto$}
\xymatrix@R=0.4cm@C=0.4cm{1' & 2' \ar[ld] \\ *+<6pt>[o][F]{1} \ar[r] & 2 \ar[ul] \ar[u]} \raisebox{-3ex}{$\longmapsto$}
\xymatrix@R=0.4cm@C=0.4cm{1' & 2' \\ 1 \ar[ru] & *+<6pt>[o][F]{2} \ar[l] \ar[lu] }\raisebox{-3ex}{$\longmapsto$}
\xymatrix@R=0.4cm@C=0.4cm{1' \ar[rd] & 2' \\ *+<6pt>[o][F]{1} \ar[ru] \ar[r] & 2}\raisebox{-3ex}{$\longmapsto$}
\xymatrix@R=0.4cm@C=0.4cm{1' \ar[rd] & 2' \ar[ld] \\ 1 & 2, \ar[l]}
\end{equation}
which yields the sequence of matrices of $c$-vectors
\begin{equation} \label{eq:A2-Cmat-seq}
\left[ \begin{array}{cc} 1 & 0 \\ 0 & 1 \end{array} \right] \ko
\left[ \begin{array}{cc} -1 & 1 \\ 0 & 1 \end{array} \right] \ko
\left[ \begin{array}{cc} 0 & -1 \\ 1 & -1 \end{array} \right] \ko
\left[ \begin{array}{cc} 0 & -1 \\ -1 & 0 \end{array} \right] \ko
\left[ \begin{array}{cc} 0 & 1 \\ -1 & 0  \end{array} \right] \ko
\left[ \begin{array}{cc} 0 & 1 \\ 1 & 0 \end{array} \right].
\end{equation}
Notice that in total, we find $6$ distinct $c$-vectors
and that these are in natural bijection with the (positive and negative)
roots of the root system corresponding to the underlying
graph $A_2$ of the quiver $Q$: We simply map a $c$-vector
with components $c_1$ and $c_2$ to the root $c_1 \alpha_1+c_2 \alpha_2$,
where $\alpha_1$ and $\alpha_2$ are the simple roots.
\[
\begin{xy}
 0;<3pt,0pt>:<0pt,3pt>::
(0,0) *{\circ} = "0",
(10,0) *+!L{\alpha_1} ="A",
(-5,8.6) *+!D{\alpha_2}="B",
"0", {\ar "A"} ,
"0", {\ar (5,8.6)},
"0", {\ar "B"},
"0", {\ar (-10,0)},
"0", {\ar (-5,-8.6)},
"0", {\ar (5,-8.6)}
\end{xy}
\]
As shown in \cite{ReadingSpeyer11}, cf.~also \cite{SpeyerThomas12a},
this bijection generalizes to all cluster-finite cluster
algebras.
In particular, we see that in these examples, each
$c$-vector is non zero and has all its
components of the same sign. This is conjectured to
be true in full generality:

\begin{main-conjecture}[\cite{FominZelevinsky07}] \label{conj:main}
Each $c$-vector
associated with a valued quiver is non zero and has either all
components non negative or all components non positive.
\end{main-conjecture}

For equally valued quivers, this conjecture follows from the
results of \cite{DerksenWeymanZelevinsky10}, which are based
on categorification using decorated representations of
quivers with potential, cf. below. Two different proofs
were given in \cite{Plamondon11a} and, up to a technical
extra hypothesis which is most probably superfluous, in \cite{Nagao10}.
In the case of valued quivers, the conjecture
is open in general, but known to be true in many important
cases thanks to the work of Demonet \cite{Demonet10}.
The determination of the $c$-vectors for general quivers
seems to be an open problem. A non acyclic example
is computed in \cite{Najera11}.

\subsection{Principal coefficients: $F$-polynomials and $g$-vectors}
\label{ss:principal-coefficients-g-vectors}
We keep the above notations $Q$, $B$, $Q_{pr}$ and $B_{pr}$. By the sharpened
Laurent phenomenon (Theorem~\ref{thm:sharpened-Laurent}), each
cluster variable of the cluster algebra $\ca(Q_{pr})$ associated
with $Q_{pr}$ is a Laurent polynomial in $x_1$, \ldots, $x_n$ with
coefficients in $\Z[x_{n+1}, \ldots, x_{2n}]$. In other words, for each
vertex $t$ of the $n$-regular tree and each $1\leq j\leq n$,
the cluster variable $x_j(t)$ belongs to the ring
\[
\Z[x_1^{\pm 1}, \ldots, x_n^{\pm 1}, x_{n+1}, \ldots, x_{2n}].
\]
The {\em $F$-polynomial}\index{$F$-polynomial}
\[
F_j(t) \in \Z[x_{n+1}, \ldots, x_{2n}]
\]
is by definition the specialization of $x_j(t)$ at
$x_1=1$, $x_2=1$, \ldots, $x_n=1$.

To define the $g$-vectors, let us endow the ring
\[
\Z[x_1^{\pm 1}, \ldots, x_n^{\pm 1}, x_{n+1}, \ldots, x_{2n}]
\]
with the $\Z^n$-grading such that
\[
\deg(x_j)=e_j  \mbox{ and }
\deg(x_{n+j})=-Be_j \mbox{ for } 1\leq j\leq n.
\]
For each vertex $t$ of the $n$-regular tree and each $1\leq j\leq n$,
the cluster variable $x_j(t)$ of $\ca(Q_{pr})$ is in fact homogeneous for this
grading (Prop.~6.1 of \cite{FominZelevinsky07}). Its degree is by definition
the {\em $g$-vector $g_j(t)$}\index{$g$-vector}. The {\em matrix of $g$-vectors
$G(t)$} has as its columns the vectors $g_j(t)$. When necessary, we
will denote this matrix by $G(B,t_0,t)$ to clarify its dependence
on $B$ and the sequence of mutations linking $t_0$ to $t$.

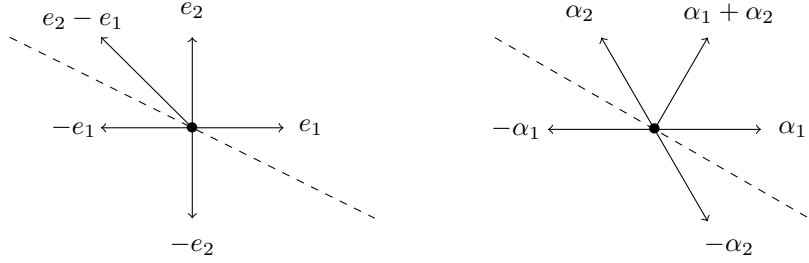
\begin{figure}
\[
\begin{tikzpicture}[scale=1.2]
\draw (0,0) node {$\bullet$};
\draw[->] (0,0)--(1,0);
\draw[->] (0,0)--(0,1);
\draw[->] (0,0)--(-1,1);
\draw[->] (0,0)--(-1,0);
\draw[->] (0,0)--(0,-1);
\draw (1.3,0) node {$e_1$};
\draw (0,1.3) node {$e_2$};
\draw (-1.2,1.2) node {$e_2-e_1$};
\draw (-1.3,0) node {$-e_1$};
\draw (0,-1.3) node {$-e_2$};
\draw[dashed] (-2,1)--(2,-1);
\end{tikzpicture}
\quad\quad\quad\quad
\begin{tikzpicture}[scale=1.4]
\draw (0,0) node {$\bullet$};
\draw[->] (0,0)--(1,0);
\draw[->] (0,0)--(0.5,0.87);
\draw[->] (0,0)--(-0.5,0.87);
\draw[->] (0,0)--(-1,0);
%\draw[->,dashed] (0,0)--(-0.5,-0.87);
\draw[->] (0,0)--(0.5, -0.87);
\draw (1.3,0) node {$\alpha_1$};
\draw (0.7,1.1) node {$\alpha_1+\alpha_2$};
\draw (-0.7,1.1) node {$\alpha_2$};
\draw (-1.3,0) node {$-\alpha_1$};
\draw (0.7,-1.1) node {$-\alpha_2$};
\draw[dashed] (-1.5,0.87)--(1.5,-0.87);
\end{tikzpicture}
\]
\caption{$g$-vectors and almost positive roots for $A_2$}
\label{fig:g-vectors-A2}
\end{figure}

For example, if $B$ is associated with $Q:1\to 2$ and we
mutate along the path
\[
\xymatrix{t_0 \ar@{-}[r]^1 & t_1 \ar@{-}[r]^2 & t_2 \ar@{-}[r]^1 & t_3 \ar@{-}[r]^2 & t_4 \ar@{-}[r]^1 & t_5}
\]
in the $2$-regular tree, then, in addition to the $g$-vectors $g_1(t_0)=e_1$
and $g_2(t_0)=e_2$ and the $F$-polynomials $F_1(t_0)=F_2(t_0)=1$
associated with the initial variables, we successively
find the following cluster variables in $\ca(B_{pr})$ and the corresponding
$F$-polynomials and $g$-vectors:
\begin{align*}
x_1(t_1)= \frac{x_2+x_3}{x_1}\ko &\quad F_1(t_1)=1+x_3  \ko\quad  g_1(t_1) = e_2-e_1\ko \\
x_2(t_2)= \frac{x_2+x_3+x_1 x_3 x_4}{x_1x_2}\ko &\quad F_2(t_2)= 1+x_3+x_3 x_4\ko\quad g_2(t_2) = -e_1 \ko\\
x_1(t_3)= \frac{1+x_1 x_4}{x_2}\ko &\quad F_1(t_3)= 1+x_4\ko \quad g_1(t_3) = -e_2\\
x_2(t_4)= x_1 \ko &\quad F_2(t_4)=1 \ko \quad g_2(t_4)= e_1 \\
x_1(t_5)= x_2 \ko &\quad F_1(t_5)=1 \ko \quad g_1(t_5)= e_2.
\end{align*}
The associated $G$-matrices are
\begin{equation} \label{eq:A2-Gmat-seq}
\left[ \begin{array}{cc} 1 & 0 \\ 0 & 1 \end{array} \right] \ko
\left[ \begin{array}{cc} -1 & 0 \\ 1 & 1 \end{array} \right] \ko
\left[ \begin{array}{cc} -1 & -1 \\ 1 & 0 \end{array} \right] \ko
\left[ \begin{array}{cc} 0 & -1 \\ -1 & 0 \end{array} \right] \ko
\left[ \begin{array}{cc} 0 & 1 \\ -1 & 0 \end{array} \right] \ko
\left[ \begin{array}{cc} 0 & 1 \\ 1 & 0 \end{array} \right].
\end{equation}
If we let $\alpha_1$ and $\alpha_2$ be the simple roots of the
root system $A_2$, then clearly the linear map which takes
$e_1$ to $\alpha_1$ and $e_2$ to $\alpha_1+\alpha_2$
yields a bijection from the set of the $g$-vectors to the
{\em set of  almost positive roots}\index{root!almost positive}, i.e. the union of the set of
positive roots with the set of opposites of the simple roots,
cf. Figure~\ref{fig:g-vectors-A2}.

This statement generalizes to an acyclic (equally valued)
quiver $Q$ as follows: For two vertices $i,j$ of $Q$, let
$p_{ij}$ be the number of paths from $i$ to $j$ (i.e. formal
compositions of $\geq 0$ arrows). Let $\alpha_1$, \ldots, $\alpha_n$
be the simple roots of the root system corresponding to the
underlying graph of $Q$. The following theorem is a consequence
of the results of \cite{CalderoKeller06}.

\begin{theorem}
The linear map taking $e_j$ to $\sum_{i=1}^n p_{ij} \alpha_i$, $1\leq j\leq n$, is
a bijection from the set of $g$-vectors of $Q$ to the
union of the set of real (positive) Schur roots with the set
of negative simple roots.
\end{theorem}

\subsection{Tropical duality}\label{ss:tropical-duality} \index{duality!tropical} Let
$Q$ be a valued quiver, $B$ the associated skew-symmetri\-zable
$n\times n$-matrix and $D$
a diagonal integer $n\times n$-matrix with strictly positive
diagonal entries such that the transpose $(DB)^T$ of $DB$
equals $-DB$. The {\em opposite valued quiver $Q^{op}$}\index{quiver!opposite valued} corresponds
to the matrix $-B$. For example, the opposite valued quiver of
\[
\vec{B}_3: \xymatrix{1 \ar[r] & 2 \ar[r]^{(1,2)} & 3} \quad\mbox{is}\quad
\vec{B}_3^{op}: \xymatrix{1 & 2 \ar[l] & 3 \ar[l]_{(2,1)} }\ko
\]
which is in fact mutation equivalent to $\vec{B}_3$
(we mutate at $1$ and $3$).

\begin{theorem}[\cite{NakanishiZelevinsky11}] \label{thm:tropical-duality}
Suppose that
the main conjecture~\ref{conj:main} holds for $Q$. Then for
each vertex $t$ of the $n$-regular tree, we have
\begin{equation} \label{eq:tropical-duality}
G(t)^T D\, C(t) = D \ko
\end{equation}
\begin{equation} \label{eq:tropical-inversion}
C(t)^{-1} = C(Q(t)^{op},t,t_0) \quad\mbox{and}\quad G(t)^{-1}=G(Q(t)^{op},t,t_0).
\end{equation}
\end{theorem}

To check the equality (\ref{eq:tropical-duality}) in the example
of the quiver $1 \to 2$, the reader may inspect the
$C$- and $G$-matrices given in (\ref{eq:A2-Cmat-seq}) and
(\ref{eq:A2-Gmat-seq}).
The equalities~(\ref{eq:tropical-inversion}) are
given in Theorem~1.2 of \cite{NakanishiZelevinsky11}.
The equality~(\ref{eq:tropical-duality}) is equation (3.11) from
\cite{NakanishiZelevinsky11}, cf. also Prop.~3.2 of \cite{Nakanishi11a}.
For skew-symmetric matrices $B$, it was first proved using
Plamondon's results \cite{Plamondon11a} in Prop.~4.1 of
\cite{Nakanishi11} by T.~Nakanishi, who had discovered the
statement by combining in Cor.~6.10 and 6.11 of \cite{Keller10a}.

Let $v: Q_1 \to \N^2$ denote the valuation of the valued quiver $Q$,
cf.~section~\ref{ss:valued-quivers}. Following \cite{FockGoncharov09},
we define
the {\em Langlands dual $Q^\vee$}\index{quiver!Langlands dual of} as the valued quiver whose underlying
oriented graph equals that of $Q$ and whose valuation $v^\vee$ is
defined by reversing the valuation of $Q$: For each arrow $\alpha$, we put
\begin{equation} \label{eq:Landlands-dual-valuation}
v^\vee(\alpha) = (v(\alpha)_2, v(\alpha)_1).
\end{equation}
The corresponding skew-symmetrizable matrix $B^\vee$ equals $-B^T$.
For example, if $Q$ is the valued quiver
\[
\vec{B}_3: \xymatrix{1 \ar[r] & 2 \ar[r]^{(1,2)} & 3}\ko
\]
then $Q^\vee$ is
\[
\vec{C}_3: \xymatrix{1 \ar[r] & 2 \ar[r]^{(2,1)} & 3}.
\]

\begin{theorem}[Th.~1.2 of \cite{NakanishiZelevinsky11}] \label{thm:tropical-Langlands-duality}
Suppose that
the main conjecture~\ref{conj:main} holds for $Q$. Then for
each vertex $t$ of the $n$-regular tree, we have
\[
G(Q,t_0,t)^T = C(Q^\vee,t_0,t)^{-1}.
\]
\end{theorem}

For example, if we successively mutate the principal extension of
the above valued quiver $\vec{C}_3$ at the vertices $1,2,3,1,2,3$, we find
the valued quiver
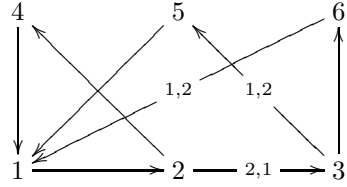
\[
\begin{xy} 0;<0.6pt,0pt>:<0pt,-0.6pt>::
(0,100) *+{1} ="0",
(100,100) *+{2} ="1",
(200,100) *+{3} ="2",
(0,0) *+{4} ="3",
(100,0) *+{5} ="4",
(200,0) *+{6} ="5",
"0", {\ar"1"},
"3", {\ar"0"},
"4", {\ar"0"},
"5", {\ar|*+{\scriptstyle 1,2}"0"},
"1", {\ar|*+{\scriptstyle 2,1}"2"},
"1", {\ar"3"},
"2", {\ar|*+{\scriptstyle 1,2}"4"},
"2", {\ar"5"},
\end{xy}
\]
and hence the $C$-matrix
\[
C(\vec{C}_3, t_0, t)= \left[ \begin{array}{ccc} 1 & -1 & 0 \\ 1 & 0 & -2 \\ 1 & 0 & -1 \end{array} \right].
\]
On the other hand, if successively mutate the initial seed of
the principal extension of $\vec{B}_3$ at $1,2,3,1,2,3$, we find the cluster
\begin{align*}
x_1(t) &= \frac{1}{x_2}(x_3^2+x_1 x^5) \ko \\
x_2(t) &= \frac{1}{x_1 x_2^2 x_3^2}(x_1^2 x_2^2 x_4 x_5^2+ 2 x_1^2 x_2 x_4 x_5^2 x_6 + \cdots + x_3^4 x_4) \ko \\
x_3(t) &= \frac{1}{x_2 x_3}(x_3^2 + x_1 x_5 + x_1 x_2 x_5 x_6)
\end{align*}
and thus the $G$-matrix
\[
G(\vec{B}_3,t_0,t) = \left[ \begin{array}{ccc} 0 & -1 & 0 \\ -1 & -1 & -1 \\ 2 & 2 & 1 \end{array} \right].
\]
This is indeed the inverse transpose of $C(\vec{B}_3,t_0,t)$.
This was to be expected by theorem~\ref{thm:tropical-Langlands-duality} since
the main conjecture holds for $\vec{C}_3$ by
Demonet's work \cite{Demonet11} \cite{Demonet10}.

\subsection{Product formulas for $c$-matrices and $g$-matrices} \label{ss:product-formulas-for-c-matrices}
\index{product formulas}
We will give a key ingredient for the proof of theorem~\ref{thm:tropical-Langlands-duality}
which is also useful in the investigation of quantum cluster algebras
(section~\ref{ss:quantum-mutations}).
Let $Q$ be a valued quiver, $B$ the associated skew-symmetri\-zable
$n\times n$-matrix and $D$ a diagonal integer $n\times n$-matrix with
strictly positive diagonal entries such that the transpose $(DB)^T$ of $DB$
equals $-DB$.
$1\leq k\leq n$ be an integer. Choose a sign $\eps$ equal to $1$ or $-1$.
Let $F_\eps=F_{k,\eps}(Q)$ be the $n\times n$-matrix which differs
from the identity matrix only in its $k$th row, whose
coefficients are given by
\[
(F_\eps)_{kj} = \left\{
\begin{array}{cc} -1               & \mbox{if }j=k; \\
                 {[\eps b_{kj}]_+} & \mbox{if } j\neq k. \end{array} \right.
\]
Let $E_\eps=E_{k,\eps}(Q)$ be the $n\times n$-matrix which differs from
the identity matrix only in its $k$th column, whose
coefficients are given by
\[
(E_\eps)_{ik} = \left\{
\begin{array}{cc} -1                & \mbox{if } i=k; \\
                 {[-\eps b_{ik}]_+} & \mbox{if } i\neq k. \end{array} \right.
\]
Notice that both $E_\eps$ and $F_\eps$ square to the identity
matrix. Parts a) to d) of the following lemma become natural in the categorical picture
to be developed in section~\ref{s:categorification}, cf. Corollary~\ref{cor:categorification-of-E-eps}.
Part e) seems harder to interpret.
\begin{lemma}
\begin{itemize}
\item[a)] We have $E_\eps \mu_k(B) = B F_\eps$ and $E_\eps^T D F_\eps = D$.
\item[b)] We have $E_{k,-\eps}(\mu_k(Q)) = E_{k,\eps}(Q)^{-1}$ and $F_{k,-\eps}(\mu_k(Q)) = F_{k,\eps}(Q)^{-1}$.
\item[c)] For $1\leq k\leq n$, let $T_k=E_{k,\eps}(\mu_k(Q)) E_{k,\eps}(Q)$. Then for
two vertices $i$, $j$, the matrices $T_i$ and $T_j$ satisfy the braid relation
associated with the full valued subquiver whose vertices are $i$ and $j$, i.e. we
have
\begin{equation} \label{eq:braid-relation}
\underbrace{T_i T_j T_i \ldots}_{m \mbox{ {\scriptsize factors}}} =
\underbrace{T_j T_i T_j \ldots}_{m \mbox{ {\scriptsize factors}}} \ko
\end{equation}
where the number of factors $m$ equals $2$, $3$, $4$ or $6$
depending on whether $|b_{ij}b_{ji}|$ equals $0$, $1$, $2$ or $3$.
\item[d)] We have $E_{k,\eps}(Q^{op})= E_{k,-\eps}(Q)$, $F_{k,\eps}(Q^{op})= F_{k,-\eps}(Q)$.
\item[e)] We have $E_{k,\eps}(Q^\vee)^T = F_{k,\eps}(Q)$.
\end{itemize}
\end{lemma}

Now let
\[
\xymatrix{t_0 \ar@{-}[r]^{i_1} & t_1 \ar@{-}[r]^{i_2} & t_2 \ar@{-}[r]^{i_3}  & \ldots \ar@{-}[r]^{i_N} & t_N} .
\]
be a path in the $n$-regular tree, let $\eps_s$ be the sign
of the $c$-vector $C(t_{s-1})e_s$ and
let $E_{i_s, \eps_s}(t_s)$ resp.~$F_{i_s,\eps_s}(t_s)$ be the matrix
$E_{\eps_s}$ resp. $F_{\eps_s}$ associated with the quiver
$Q(t_{s-1})$ and the vertex $i_s$, $1\leq s\leq N$.

\begin{theorem}[\cite{NakanishiZelevinsky11}] \label{thm:product-formula-for-c-matrices}
If the main conjecture~\ref{conj:main}
holds for $Q$, we have
\[
G(t_N) = E_{i_1,\eps_1}(t_1) \ldots\, E_{i_N,\eps_N}(t_N) \mbox{ and }
C(t_N) = F_{i_1,\eps_1}(t_1) \ldots\, F_{i_N,\eps_N}(t_N).
\]
\end{theorem}

\subsection{Cluster algebras with coefficients in a semifield} \label{ss:cluster-alg-semifield}
A {\em semifield}\index{semifield} is an abelian group $\P$ endowed with an additional binary
operation $\oplus: \P\times \P\to \P$ which is commutative, associative
and distributive with respect to the group law of $\P$. For example,
the {\em tropical semifield}\index{semifield!tropical} $\mbox{Trop}(u_1, \ldots, u_n)$ is the
free (multiplicative) abelian group generated by the indeterminates $u_i$
endowed with the operation $\oplus$ defined by
\[
(\prod u_i^{l_i}) \oplus (\prod u_i^{m_i}) = \prod u_i^{\min(l_i,m_i)}.
\]
Clearly, it is isomorphic to $\Z_{\mbox{\scriptsize trop}}^n$, where
$\Z_{\mbox{\scriptsize trop}}$ is the abelian group $\Z$ endowed with the
operation $\oplus$ defined by $x\oplus y=\min(x,y)$.
It is shown in Lemma~2.1.6 of \cite{BerensteinFominZelevinsky96} that
the {\em universal semifield}\index{semifield!universal} $\Q_{sf}(x_1, \ldots, x_n)$ on given indeterminates
$x_1$, \ldots, $x_n$ is the closure, in $\Q(x_1, \ldots, x_n)$, of the
set $\{x_1, \ldots, x_n\}$ under multiplication, division and addition.
Notice that this closure contains polynomials whose coefficients are not
all positive; for example, the polynomial
\[
x^2 - x + 1 = \frac{x^3+1}{x+1}
\]
belongs to $\Q_{sf}(x)$. The abelian group underlying a semifield $\P$
is torsion-free. Indeed, if an element $x$ satisfies $x^m=1$, then
\[
x=\frac{x^m\oplus x^{m-1} \oplus \cdots \oplus x}{x^{m-1}\oplus x^{m-2} \oplus \cdots \oplus 1}
 =\frac{1 \oplus x^{m-1} \oplus \cdots \oplus x}{x^{m-1}\oplus x^{m-2} \oplus \cdots \oplus 1}=1.
\]
Thus, the group ring $\Z\P$ is integral.

Let us fix a semifield $\P$ and an integer $n\geq 1$.
A {\em $Y$-seed}\index{$Y$-seed} of rank $n$ with values in $\P$
is a pair $(Q,Y)$ formed by a valued quiver $Q$ with $n$ vertices
and by a sequence $Y=(y_1, \ldots, y_n)$ of elements of $\P$. Let $B$
be the skew-symmetrizable matrix corresponding to $Q$.
If $(Q,Y)$ is a $Y$-seed and $k$ a vertex of $Q$, the {\em mutated $Y$-seed}\index{mutation!of $Y$-seeds} $\mu_k(Q,Y)$
is the $Y$-seed $(Q',Y')$ where $Q'=\mu_k(Q)$ and, for $1\leq j\leq n$, we have
\begin{equation} \label{eq:Y-var-mutation}
y'_j = \left\{ \begin{array}{ll} y_k^{-1} & \mbox{if $j=k$;} \\
       y_j y_k^{[b_{kj}]_+} (1\oplus y_k)^{-b_{kj}} & \mbox{if $j\neq k$.}
              \end{array} \right.
\end{equation}
One checks that $\mu_k^2(Q,Y)=(Q,Y)$.
For example, the following $Y$-seeds are related by a mutation at
the vertex $1$
\[
\xymatrix{ y_1 \ar[r] & y_2 \ar[ld] \\ y_3 \ar[u] \ar[r] & y_4 \ar[u]} \quad\quad\quad
\xymatrix{ 1/y_1 \ar[d] & y_2/(1\oplus y_1^{-1}) \ar[l] \\ y_3(1\oplus y_1) \ar[r] & y_4\ar[u]} \ko
\]
where we write the element $y_i$ in place of the vertex $i$.

Let $\Q\P$ be the fraction field of the group ring $\Z\P$ and
$\cf$ any field obtained from $\Q\P$ by adjoining $n$ indeterminates.
A {\em seed}\index{seed!semifield-valued} with coefficients in $\P$ is a triple $(Q,Y,X)$, where
$(Q,Y)$ is a $Y$-seed of rank $n$ with values in $\P$ and $X$
is a sequence $(x_1, \ldots, x_n)$ of elements of $\cf$ which
freely generate the field $\cf$. If $(Q,Y,X)$ is a seed and
$k$ a vertex of $Q$, the mutation $\mu_k(Q,Y,X)$ is the
seed formed by the mutation $\mu_k(Q,Y)$ and the
sequence $X'$ with $x'_j=x_j$ for $j\neq k$ and
$x'_k$ defined by the exchange relation
\begin{equation} \label{eq:exchange-rel-semifield}
 x'_k x_k (1\oplus y_k) = y_k  \prod_i x_i^{[b_{ik}]_+} + \prod_i x_i^{[-b_{ik}]_+}.
\end{equation}

A {\em seed pattern}\index{seed pattern}\index{pattern!of seeds} is the datum, for each vertex $t$
of the $n$-regular tree, of a seed $(Q(t), Y(t), X(t))$ such
that if $t$ and $t'$ are linked by an edge labeled $k$, then
the seeds corresponding to $t$ and $t'$ are linked by the
mutation at $k$. The {\em cluster algebra}\index{cluster algebra!general} is the
$\Z\P$-subalgebra of the field $\cf$ generated by the
cluster variables.

We recover the cluster algebra of geometric type
associated with an $m \times n$-matrix $\tilde{B}$ as follows:
We let $B$ be the principal part of $\tilde{B}$; we define the
semifield $\P$ to be the tropical semifield $\mbox{Trop}(x_{n+1}, \ldots, x_{m})$
and the initial $Y$-variables to be
\[
y_j = \prod_{i={n+1}}^m x_i^{b_{ij}} \ko 1 \leq j\leq n.
\]
As a simple example of a cluster algebra of `non geometric' type,
consider the case where $n=1$, $\P=\Q_{sf}(y_1,y_2)$ and $Q: 1 \to 2$.
Then the sequence of mutations
\[
\xymatrix{t_0 \ar@{-}[r]^1 & t_1 \ar@{-}[r]^2 & t_2 \ar@{-}[r]^1 & t_3 \ar@{-}[r]^2 & t_4 \ar@{-}[r]^1 & t_5}
\]
starting from the initial seed $(1 \to 2, \{x_1,x_2\}, \{y_1,y_2\})$ yields
\begin{align*}
y_1(t_1) &=\frac{1}{y_1} \ko              & y_2(t_1) &=\frac{y_1 y_2}{1+y_1}\ko   & x_1(t_1) &= \frac{y_1+x_2}{x_1(1+y_1)} \\
y_1(t_2) &=\frac{y_2}{1+y_1+y_1 y_2}  \ko & y_2(t_2) &=\frac{1+y_1}{y_1 y_2}\ko   & x_2(t_2) &= \frac{x_1 y_1 y_2 + x_2+y_1}{x_1 x_2 (1+y_1 + y_2 y_1)} \\
y_1(t_3) &=\frac{1+y_1 + y_1 y_2}{y_2}\ko & y_2(t_3) &=\frac{1}{y_1(1 + y_2)} \ko & x_1(t_3) &=\frac{x_1 y_2 + 1}{x_2 (1+y_2)} \\
y_1(t_4) &=\frac{1}{y_2} \ko  & y_2(t_4)  &=y_1(1+y_2) \ko  & x_2(t_4) &=x_1 \\
y_1(t_5) &= y_2 \ko                       & y_2(t_5) &= y_1 \ko                    & x_1(t_5) &=x_2.
\end{align*}

\subsection{The separation formulas}\index{separation formula!classical} \label{ss:separation-formulas}
Let a seed pattern be given and let
us write $(Q,Y,X)$ for the initial seed $(Q(t_0),Y(t_0),X(t_0))$
associated with the chosen root $t_0$ of the $n$-regular tree.
Let us write $c_{ij}(t)$ for the coefficients of the
$c$-matrix $C(t)$ and $g_{ij}(t)$ for those of the
$g$-matrix $G(t)$ associated with a vertex $t$ of the $n$-regular tree.
Recall that
\[
F_j(t) \in \Z[x_{n+1}, \ldots, x_{2n}] \ko 1\leq j\leq n\ko
\]
are the $F$-polynomials at the vertex $t$. By construction,
they belong to the universal semifield $\Q_{sf}(x_{n+1}, \ldots, x_{2n})$
and thus it makes sense to consider their evaluations
\[
F_j(t)(y_1, \ldots, y_n)
\]
at the elements $y_1, \ldots, y_n$ of $\P$ and more generally,
at an $n$-tuple of elements of any semifield.

\begin{theorem}[Prop.~3.13 and Cor.~6.3 of \cite{FominZelevinsky07}]
\label{thm:separation-formula}
For each vertex $t$ of the $n$-regular tree and each $1\leq j\leq n$,
we have
\begin{align}
y_j(t) &= y_1^{c_{1j}(t)} \ldots\, y_n^{c_{nj}(t)} \prod_{i} F_i(t)(y_1, \ldots, y_n)^{b_{ij}(t)} \\
x_j(t) &= x_1^{g_{1j}(t)} \ldots\, x_n^{g_{nj}(t)} \frac{F_{j}(t)(\hat{y}_1, \ldots, \hat{y}_n)}{F_j(t)(y_1, \ldots, y_n)} \ko \label{eq:separation-classical}
\end{align}
where $\hat{y}_l = y_j \prod_i x_i^{b_{il}}$, $1 \leq l \leq n$.
\end{theorem}

\section{Quantum cluster algebras and quantum dilogarithms}
\label{s:quantum-cluster-algebras-and-quantum-dilogarithms}

\subsection{The quantum dilogarithm} \label{ss:quantum-dilog}
Let $q^{1/2}$ be an indeterminate. We will denote its square by $q$.
The (exponential of) the {\em quantum dilogarithm series}\index{quantum dilogarithm}\index{dilogarithm!quantum} is
\[
\E(y)=\E_q(y) = 1 + \frac{q^{1/2}}{q-1} y + \ldots + \frac{q^{n^2/2}}{(q^n-1)(q^{n-1}-1) \ldots (q-1)} + \ldots .
\]
It is a series in the indeterminate $y$ with coefficients in the field $\Q(q^{1/2})$.
It is related to the classical dilogarithm\index{dilogarithm!classical}
\[
\mbox{Li}_2(x) = \sum_{n=1}^\infty \frac{x^n}{n^2} = - \int_0^x \frac{\log(1-y)}{y} \, dy\ko |x|<1 \ko
\]
by the asymptotic expansion
\[
\E_q(y) \sim \exp (-\frac{\mbox{Li}_2(-y)}{\log(q)})
\]
when $q$ goes to $1^-$. An easy computation shows that we have
the functional equation
\begin{equation} \label{eq:functional-eq}
(1+q^{1/2} y)\, \E(y) = \E(qy).
\end{equation}
The quantum dilogarithm is related to the classical $q$-exponential
function by the substitution $y\mapsto \frac{q^{1/2}}{q-1}y$. Therefore,
as discovered by Sch\"utzenberger \cite{Schuetzenberger53}, if $y_1$ and
$y_2$ are two indeterminates which $q$-commute, i.e. $y_1 y_2 = q y_2 y_1$,
then we have
\begin{equation} \label{eq:q-exp-eq}
\E(y_1+y_2)= \E(y_2)\E(y_1).
\end{equation}
In 1993, Faddeev, Kashaev and Volkov \cite{FaddeevKashaev94} \cite{FaddeevVolkov93}
discovered that (\ref{eq:functional-eq}) and (\ref{eq:q-exp-eq}) together
imply the {\em pentagon identity}\index{pentagon identity}
\begin{equation} \label{eq:pentagon}
y_1 y_2 = q y_2 y_1 \Longrightarrow \E(y_1) \E(y_2) = \E(y_2) \E(q^{-1/2} y_1 y_2) \E(y_1) \ko
\end{equation}
cf.~\cite{Volkov11} for a recent account. Their main result
states that this identity implies the classical five-term identity
\[
L(x)+L(y)-L(xy) = L(\frac{x-xy}{1-xy}) + L(\frac{y-xy}{1-xy})
\]
for the Rogers dilogarithm\index{dilogarithm!Rogers}
\[
L(x) = \mbox{Li}_2(x) + \log(1-x)\log(x)/2.
\]
We refer to \cite{Nakanishi11} \cite{KashaevNakanishi11} for more information
on the many recent developments around this subject and to
\cite{Zagier07} for more information on the dilogarithm function.

\subsection{Quantum mutations and quantum cluster algebras}
\label{ss:quantum-mutations} \index{mutation!quantum}
We will construct quantum cluster algebras following
Berenstein--Zelevinsky \cite{BerensteinZelevinsky05}. Quantum
cluster algebras are certain non commutative deformations of
cluster algebras of geometric type. Let $1\leq n\leq m$ be
integers, $\tilde{B}$ an integer $m\times n$-matrix with
skew-symmetrizable principal part $B$ and $\Lambda$ a
skew-symmetric integer $m\times m$-matrix. Let $\tilde{Q}$
and $Q$ be the associated valued ice quivers.
Recall from section~\ref{s:cluster-algebras-of-geometric-type} that the
datum of $\tilde{B}$ gives rise to a cluster algebra of
geometric type. Let us assume that $(\Lambda, \tilde{B})$
is a {\em compatible pair}\index{compatible pair}, i.e. we have
\begin{equation} \label{eq:compatible-pair}
\tilde{B}^T \Lambda = [D \, 0]\ko
\end{equation}
where $D$ is a diagonal $n\times n$-matrix whose diagonal
coefficients are strictly positive integers. This will ensure
that $\Lambda$ gives rise to a (non commutative) deformation
of the cluster algebra associated with $\tilde{B}$. We first
need to define the mutation of compatible pairs: Let
$1\leq k\leq n$ be an integer and choose a sign $\eps$ equal to $1$ or $-1$.
In the notations of section~\ref{ss:product-formulas-for-c-matrices},
let $F_\eps$ be the $n\times n$-matrix $F_{k,\eps}(Q)$ and
$E_\eps$ the $m\times m$-matrix $E_{k,\eps}(\tilde{Q})$.
The {\em mutation $\mu_k(\tilde{B},\Lambda)$}\index{mutation!of compatible pairs} is defined
to be the compatible pair $(\tilde{B}', \Lambda')$ with
\[
\tilde{B}'= E_\eps \tilde{B} F_\eps\quad\mbox{and}\quad
\Lambda' = E_\eps^T \Lambda E_\eps.
\]
One checks that $\tilde{B}'$ equals $\mu_k(\tilde{B})$ and
that $(\Lambda', \tilde{B}')$ does not depend on the choice of~$\eps$
and is again a compatible pair (with the same matrix $D$).
One checks that mutation of compatible pairs is an involution.
Thus, given a compatible pair $(\tilde{B},\Lambda)$, we can
assign a compatible pair $(\tilde{B}(t), \Lambda(t))$ to each
vertex $t$ of the $n$-regular tree such that the given
pair is assigned to $t_0$ and, whenever $t$ and $t'$ are linked
by an edge labeled $k$, the corresponding pairs are related
by the mutation at $k$.

The {\em quantum affine space $\mathbb{A}_\Lambda$}\index{quantum affine space} associated with $\Lambda$
is by definition the $\Z[q^{\pm 1/2}]$-algebra generated by all
symbols $x^\alpha$, $\alpha\in \N^m$, subject to the relations
\[
x^\alpha x^\beta = q^{\frac{1}{2} \alpha^T \Lambda \beta} x^{\alpha+\beta}.
\]
The {\em quantum torus $\mathbb{T}_\Lambda$}\index{quantum torus} is defined similarly
on generators $x^\alpha$, $\alpha\in\Z^m$. One checks that the underlying 
$\Z[q^{\pm 1/2}]$-module
of $\mathbb{A}_\Lambda$ resp. $\mathbb{T}_\Lambda$ is free on the basis formed by the
$x^\alpha$, $\alpha\in\N^m$ resp. $\alpha\in\Z^m$.
The {\em completed quantum affine space
$\hat{\mathbb{A}}_\Lambda$} is the completion of $\mathbb{A}_\Lambda$ with respect
to the kernel of the projection $\mathbb{A}_Q\to \Z[q^{\pm 1/2}]$.
The algebras $\mathbb{A}_\Lambda$ and $\mathbb{T}_\Lambda$ are Ore domains (cf. the Appendix to \cite{BerensteinZelevinsky05})
and so have a field of fractions $\mathbb{F}_\Lambda$ whose elements
are given by right fractions (or left fractions).

The {\em initial quantum seed}\index{seed!quantum} is $(\tilde{B}, \Lambda, X)$,
where $X$ is the sequence of the $x_i = x^{e_i}$. Its {\em mutation
at $k$}\index{mutation!of quantum seeds}, where $1\leq k\leq n$, is $(\tilde{B}', \Lambda', X')$, where the sequence
$X'$ is formed by the $x_i$, $i\neq k$, and by the element
$x'_k$ defined by the {\em quantum exchange relation}\index{exchange relation!quantum}
\begin{equation} \label{eq:quantum-exchange-relation}
x'_k = x^{E_+ e_k} + x^{E_- e_k}.
\end{equation}
By part (3) of Prop.~4.7 of \cite{BerensteinZelevinsky05}, there
is a unique morphism of $\Z[q^{\pm 1/2}]$-algebras
\[
\mu_k^{\#}: \mathbb{A}_{\Lambda'} \to \mathbb{T}_{\Lambda}
\]
taking $x_i$ to $x'_i$, $1\leq i\leq m$; moreover it is injective and induces an
isomorphism
\[
\mu_k^{\#}: \mathbb{F}_{\Lambda'} \iso \mathbb{F}_{\Lambda}.
\]
One checks that mutation of quantum seeds is an involution.
Thus, with each vertex $t$ of the $n$-regular tree, one can
associate a quantum seed $(\tilde{B}(t), \Lambda(t), X(t))$
such that the initial quantum seed is associated with $t_0$
and seeds with vertices $t$ and $t'$ linked by an
edge labeled $k$ are related by a quantum mutation. The
{\em quantum cluster variables}\index{cluster variable!quantum} are the $x_j(t)$, $1\leq j\leq n$,
associated with the vertices $t$ of the $n$-regular tree.
The {\em quantum cluster algebra}\index{cluster algebra!quantum} is the $\Z[q^{1/2}]$-subalgebra
of $\mathbb{F}_\Lambda$ generated by the quantum cluster variables.
We have the quantum Laurent phenomenon:

\begin{theorem}[Cor.~5.2 of \cite{BerensteinZelevinsky05}]
The quantum cluster variables are contained in the quantum
torus $\mathbb{T}_\Lambda$.
\end{theorem}

We refer to \cite{BerensteinZelevinsky05} \cite{GrabowskiLaunois11} \cite{GeissLeclercSchroeer11}
for examples of quantum cluster algebras.
The {\em exchange graph}\index{exchange graph!quantum} of quantum seeds associated
with $(\tilde{B}, \Lambda)$ is defined in analogy with
the exchange graph of (classical) seeds associated with
$\tilde{B}$, cf. section~\ref{ss:seed-mutation}. The
specialization map
\[
\Z[q^{\pm 1/2}] \to \Z
\]
taking $q^{1/2}$ to $1$ yields a morphism of
$\Z[q^{\pm 1/2}]$-modules
\[
\mathbb{T}_\Lambda \to \Z[x_1^{\pm 1}, \ldots, x_m^{\pm 1}]
\]
which takes quantum cluster variables to classical ones
and induces a map from the quantum exchange graph to
the classical exchange graph.

\begin{theorem}[Th.~6.1 of \cite{BerensteinZelevinsky05}] \label{thm:quantum-exchange-graph}
The specialization at $q^{1/2}=1$ yields an isomorphism from the
quantum to the classical exchange graph.
\end{theorem}

\subsection{Fock-Goncharov's separation formula}\index{separation formula!Fock-Goncharov's} \label{ss:fock-goncharov-separation}
Recall that the numbers $d_i$ are the coefficients of the diagonal matrix $D$
appearing in the compatibility condition~(\ref{eq:compatible-pair}).
We consider the mutation at $k$ of a given initial quantum
seed $(\tilde{B},\Lambda, X)$.

\begin{lemma}[\cite{FockGoncharov09}] We have the separation
formulas
\begin{equation} \label{eq:separation-one}
\mu_k^{\#} = \Ad'(\E_{q^{d_k}}(y_k)) \circ \phi_{k,+} = \Ad'(\E_{q^{d_k}}(y_k^{-1}))^{-1} \circ \phi_{k,-} \ko
\end{equation}
where the right adjoint action $\Ad'(u)$ takes an element $v$ to $u^{-1}v u$, we put
\[
y_k = x^{\tilde{B} e_k}
\]
and $\phi_{k,\eps} : \mathbb{T}_{\Lambda'} \to \mathbb{T}_{\Lambda}$ is the unique morphism of $\Z[q^{\pm 1/2}]$-algebras taking $x^\alpha$ to $x^{E_\eps \alpha}$.
\end{lemma}

Thus, we have separated the mutation isomorphism into a `tropical' part
and a `transcendental' part.
Notice that in order to give meaning to the formulas~(\ref{eq:separation-one}),
we need to embed the quantum tori into suitable localizations of
completions of quantum affine space. Using
formula~(\ref{eq:functional-eq}) one then checks the claim. Of course,
one would like to iterate this formula. The iteration should be
meaningful in (at least) two ways:
\begin{itemize}
\item[(1)] the product of the appearing power series should have
a meaning, i.e. all the series should live in a common completion of
quantum affine space;
\item[(2)] the composition of the `tropical parts' should have
a meaning from the point of view of `tropical' cluster theory,
as we have seen it in sections~\ref{ss:principal-coefficients-c-vectors}
and \ref{ss:principal-coefficients-g-vectors}.
\end{itemize}
In order to obtain both, it is essential to choose the sign $\pm$
in each factor carefully. This can be achieved using the
main conjecture~\ref{conj:main}.

\subsection{The quantum separation formula}\index{separation formula!quantum} \label{ss:quantum-separation}
To simplify the notations, let us assume from now on that
$(\tilde{B},\Lambda)$ is {\em unitally compatible}\index{compatible pair!unitally}, i.e.
equation~(\ref{eq:compatible-pair}) holds with $D$ the $n\times n$-identity
matrix. Let $\mathbf{i}=(i_1, \ldots, i_N)$ be a sequence of
vertices in $\{1, \ldots, n\}$. Consider the corresponding
path in the $n$-regular tree
\[
\xymatrix{t_0 \ar@{-}[r]^{i_1} & t_1 \ar@{-}[r]^{i_2} & t_2 \ar@{-}[r]^{i_3}  & \ldots \ar@{-}[r]^{i_N} & t_N} .
\]
It yields a chain of mutation isomorphisms between the
associated quantum tori:
\[
\xymatrix{\mathbb{T}_\Lambda & \mathbb{T}_{\Lambda(t_1)} \ar[l]_{\mu^{\#}_{i_1}} & \mathbb{T}_{\Lambda(t_2)} \ar[l]_{\mu^{\#}_{i_2}}
 & \ldots \ar[l]_{\mu^{\#}_{i_3}} & \mathbb{T}_{\Lambda(t_N)} \ar[l]_{\mu^{\#}_{i_N}}} .
\]
Let us write $\Phi(\mathbf{i})$ for the composition of these isomorphisms.
We would like to write down a separation formula for $\Phi(\mathbf{i})$
which generalizes (\ref{eq:separation-one}). We need some more notation:
For $1\leq s\leq N$, let $\beta_s$ be the $c$-vector $C(t_{s-1})e_{i_s}$
and let $\eps_s$ be the common sign of the components of $\beta_s$
(cf. section~\ref{ss:principal-coefficients-g-vectors}).
For a vector $\alpha$ in $\Z^n$, let us write $\E(\alpha)$ for
$\E(y^\alpha)$, where $y^\alpha=x^{\tilde{B} \alpha}$.

\begin{theorem}[\cite{Nagao10}] \label{thm:quantum-separation}
Put
\begin{align}
\E(\mathbf{i}) &= \E(\eps_N\beta_N)^{\eps_N} \ldots\  \E(\eps_1\beta_1)^{\eps_1} \\
\phi(\mathbf{i}) &= \phi_{i_1, \eps_1} \circ \ldots \circ \phi_{i_N, \eps_N}.
\end{align}
Then we have
\begin{equation} \label{eq:separation-many}
\Phi(\mathbf{i})= \Ad'(\E(\mathbf{i})) \circ \phi \ko
\end{equation}
the isomorphism $\phi$ sends $x^\alpha$ to $x^{G(t_N) \alpha}$,
where $G(t_N)$ is the $g$-matrix at $t_N$ (section~\ref{ss:principal-coefficients-g-vectors}),
and $\Ad'(\E(\mathbf{i}))$ acts on $x^{e_i}$ by multiplication with
the quantum $F$-poly\-nomial of \cite{Tran11}.
\end{theorem}

Notice that by construction all the vectors $\eps_s \beta_s$ have
non negative components so that all the series $\E(\eps_s \beta_s)^{\eps_s}$
belong to the same completion of quantum affine space.
If we replace the right adjoint action of $\E(\mathbf{i})$ by the multiplication
with the quantum $F$-polynomials, we obtain Tran's formula
(Theorem~6.1 of \cite{Tran11}), which is the quantum analogue
of Fomin-Zelevinsky's \cite{FominZelevinsky07}
separation formula~(\ref{eq:separation-classical}).
The theorem is due, in a different language, to Nagao \cite{Nagao10}
(cf. also Theorem~5.1 in \cite{Efimov11}). Alternatively, using
Theorem~\ref{thm:product-formula-for-c-matrices} and Tran's formula, it is not
hard to prove the analogous theorem for arbitrary valued quivers $Q$
for which the main conjecture~\ref{conj:main} holds.

Let us keep the notations from Theorem~\ref{thm:quantum-separation}. It is
not hard to check that there is a unique $\Q(q^{1/2})$-algebra
embedding
\[
\mathbb{A}_{B} \to \mathbb{T}_\Lambda
\]
taking an element $x^\alpha$ to $x^{\tilde{B}\alpha}$ (if $D$ is
not the identity matrix, it is an embedding $\mathbb{A}_{DB} \to \mathbb{T}_\Lambda$).
Thus, by construction, the product $\E(\mathbf{i})$ lies in a
completed quantum affine subspace isomorphic to $\hat{\mathbb{A}}_B$ and independent of
the choice of the non principal part in $\tilde{B}$. For example, we can always choose
$\tilde{B}=B_{pr}$, cf. section~\ref{ss:principal-coefficients-c-vectors}
and
\[
\Lambda= \left[ \begin{array}{cc} 0 & -I \\ I & B^T \end{array}\right] .
\]

\begin{theorem}[\cite{Keller11a} \cite{Nagao11}] \label{thm:quantum-dilog-identity}
\begin{itemize}
\item[a)] If $C(t_N)$ is a permutation matrix, then $\E(\mathbf{i})=1$.
\item[b)] If the opposite matrix $-C(t_N)$ is a permutation matrix, then 
$\E(\mathbf{i}) \in \hat{\mathbb{A}}_B$ is
Kontsevich-Soibelman's non commutative Donaldson--Thomas invariant
\cite{KontsevichSoibelman08} associated with the quiver corresponding
to $B$ (when this invariant is defined, cf. section~~\ref{ss:proof-quantum-dilog}).
\end{itemize}
\end{theorem}

\begin{remark} \label{rem:combinatorial-DT-invariant} One can
sharpen part a) as follows: Let
$\mathbf{i}$ and $\mathbf{i}'$ be two sequences of vertices in
$\{1, \ldots, n\}$ and let $t$ and $t'$ be the end points of the
corresponding paths in the $n$-regular tree. Suppose that
we have $P C(t)=C(t')$ for a permutation matrix $P$. We will
show in section~\ref{ss:proof-quantum-dilog} that
we then have $\E(\mathbf{i})=\E(\mathbf{i}')$. Thus, if $Q$ admits
some sequence $\mathbf{i}$ such that $-C(t)$ is a permutation
matrix, then the series $\E(\mathbf{i}) \in \hat{\mathbb{A}}_B$ is
independent of the choice of the sequence $\mathbf{i}$ with
this property. We then call this series the
{\em combinatorial DT-invariant} associated with $Q$.
\end{remark}

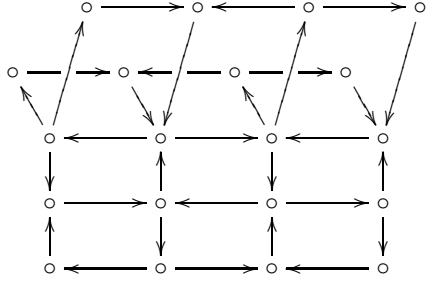
\begin{figure}
\[
\xymatrix@C=0.1cm@R=0.5cm{ & & \circ \ar[rrr] &  &  & \circ \ar[ldd]
& & &
\circ \ar[lll] \ar[rrr] & & & \circ \ar[ldd] \\
     \circ \ar[rrr]|!{[dr];[urr]}\hole & &  & \circ \ar[rd] & & &
\circ \ar[lll]|!{[lu];[dll]}\hole \ar[rrr]|!{[rd];[rru]}\hole & & & \circ \ar[rd] & &  \\
           & \circ \ar[lu] \ar[ruu] \ar[d] & & & \circ \ar[lll] \ar[rrr] & & &
\circ \ar[ruu] \ar[lu] \ar[d] & & & \circ \ar[lll] & \\
           & \circ \ar[rrr] & & & \circ \ar[d] \ar[u] & & &
\circ \ar[lll] \ar[rrr] & & & \circ \ar[u] \ar[d] & \\
           & \circ \ar[u] & & & \circ \ar[lll] \ar[rrr] &  & &
\circ \ar[u] & & & \circ \ar[lll] & }
\]
\caption{The quiver $\vec{A}_4 \square \vec{D}_5$}
\label{fig:a4-square-d5}
\end{figure}

We will give a proof of the theorem and the remark 
in section~\ref{ss:proof-quantum-dilog}, 
cf. also Theorem~3.5 in \cite{KashaevNakanishi11}. 
Let us illustrate the theorem on the example of the mutation sequence
$\mathbf{i}=(1,2,1,2,1)$ of the quiver $\vec{A}_2: 1\to 2$. We have computed
the sequence of $c$-matrices $C(t_s)$, $1\leq s\leq 5$, in equation
(\ref{eq:A2-Cmat-seq}). We obtain
\[
\beta_1 = e_1\ko \beta_2 = e_1+e_2\ko \beta_3 = e_2 \ko \beta_4 = -e_1 \ko
\beta_5 = -e_2.
\]
Since $C(t_5)$ is the matrix of the transposition $e_1 \leftrightarrow e_2$,
part a) of the theorem yields the identity
\[
\E(e_2)^{-1} \E(e_1)^{-1} \E(e_2) \E(e_1+e_2) \E(e_1) = 1 \ko
\]
which is of course equivalent to the pentagon identity (\ref{eq:pentagon}).
Since $C(t_3)$ is the opposite of the transposition matrix, we find
that Kontsevich-Soibelman's DT invariant equals
\[
\E(e_2) \E(e_1+e_2) \E(e_1)
\]
for the quiver $\vec{A}_2$, as is well-known, cf. Example~2)
in section~6.4 of \cite{KontsevichSoibelman08}. This example can
be generalized to any Dynkin quiver, which yields a family
of quantum dilogarithm identities due to Reineke \cite{Reineke10},
cf. also Cor.~1.7 in \cite{Keller11a} and \cite{Qiu11a} \cite{Qiu11b}. 
Namely, let $\Delta$ be a simply laced Dynkin diagram and let $Q$ 
be an alternating quiver (i.e. each vertex is a source or a sink) whose 
underlying graph is $\Delta$. Let $i_+$ be the sequence of sources of $Q$
and $i_-$ its sequence of sinks (in any order). Let
\[
\mathbf{i}=\underbrace{i_+ i_- i_+ \ldots}_{h \mbox{ {\scriptsize factors}}} \ko
\]
where $h$ is the Coxter number of $\Delta$
and let $\mathbf{i}'={i_- i_+}$. Let $\mu_{\mathbf{i}}(t_0)$ be the
final vertex in the path in the regular tree which starts at $t_0$
and runs through the sequence of mutations $\mathbf{i}$ starting
at the leftmost vertex in the sequence. Then
one can show that both $-C(\mu_{\mathbf{i}}(t_0))$ and  $-C(\mu_{\mathbf{i}'}(t_0))$
are permutation matrices and so the Kontsevich-Soibelman invariant
associated with $Q$ equals
\[
\E(\mathbf{i})=\E(\mathbf{i}') \ko
\]
which is Reineke's identity associated with $Q$. One can further
generalize this class as follows: Let $\Delta$ and $\Delta'$ be
two simply laced Dynkin diagrams and $\vec{\Delta}$ and $\vec{\Delta}'$ alternating
quivers with underlying graphs $\Delta$ and $\Delta'$. Let
$Q$ be the square product $\vec{\Delta}\square \vec{\Delta}'$ as
defined in section~3.3 of \cite{Keller10a}. For example, the
square product of the quivers
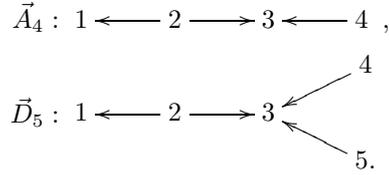
\begin{align*}
\vec{A}_4 &: \xymatrix{ 1 & 2 \ar[l] \ar[r] & 3 & 4 \ar[l] } \ko\\
\vec{D}_5 &: \raisebox{0.75cm}{\xymatrix@R=0.2cm{   &  &  & 4 \ar[dl] \\ 1 & 2 \ar[l] \ar[r] & 3 & \\
& & & 5. \ar[ul] }}
\end{align*}
is depicted in Figure~\ref{fig:a4-square-d5}.
Let $i_+$ be the sequence of all source-sinks of $\vec{\Delta}\square \vec{\Delta}'$
(i.e. vertices $(u,v)$ such that $u$ is a source in the full
subquiver $p_2^{-1}(v)$ and $v$ a sink the full subquiver $p_1^{-1}(u)$,
where the $p_i$ are the projections) and let $i_-$ be
the sequence of all sink-sources. Let
\[
\mathbf{i}=\underbrace{i_+ i_- i_+ \ldots}_{h \mbox{ {\scriptsize factors}}}\quad\mbox{and}\quad
\mathbf{i'}=\underbrace{i_- i_+ i_- \ldots}_{h' \mbox{ {\scriptsize factors}}} \ko
\]
where $h$ is the Coxeter number of $\Delta$ and $h'$ that of $\Delta'$.
Again one can show that both $-C(\mu_{\mathbf{i}}(t_0))$ and  $-C(\mu_{\mathbf{i}'}(t_0))$
are permutation matrices and so the Kontsevich-Soibelman invariant
associated with $\vec{\Delta}\square \vec{\Delta}'$ equals
\[
\E(\mathbf{i})=\E(\mathbf{i}').
\]
In physics, a related method for computing this invariant is
the {\em mutation method}\index{mutation!method} developed and applied in \cite{CecottiEtAl11}.

\section{Categorification}
\label{s:categorification}
\index{categorification!of cluster algebras}

The setup we will describe uses triangulated $3$-Calabi-Yau categories
(derived categories of Ginzburg dg algebras).
It is due to Kontsevich-Soibelman \cite{KontsevichSoibelman10} and Nagao \cite{Nagao10}.
It is closely related to that of Plamondon \cite{Plamondon11}, who
uses triangulated $2$-Calabi-Yau categories (cluster categories).
Both build on work by Derksen-Weyman-Zelevinsky on quivers with potentials
\cite{DerksenWeymanZelevinsky08}, who first proved a statement
equivalent to the main theorem~\ref{thm:decategorification} using decorated
representations of quivers with potentials \cite{DerksenWeymanZelevinsky10}.

\subsection{Mutation of quivers with potential}
\label{ss:quivers-with-potential} We follow Derksen-Weyman-Zele\-vinsky's
fundamental article \cite{DerksenWeymanZelevinsky08}. Let $Q$ be a finite
quiver. Let $\hat{\C Q}$ be the {\em completed path algebra}\index{path algebra!completed}, \ie the
completion of the path algebra at the ideal generated by the arrows
of $Q$. Thus, $\hat{\C Q}$ is a topological algebra and the paths of
$Q$ form a topologial basis so that the underlying vector space of
$\hat{\C Q}$ is
\[
\prod_{p \mbox{ \tiny path}} kp
\]
and the multiplication is induced from the composition of paths
(we compose paths in the same way as we compose morphisms).
The {\em continuous zeroth Hochschild homology}\index{Hochschild homology} of $\hat{\C Q}$ is the vector
space $HH_0(\hat{\C Q})$ obtained as the quotient of $\hat{\C Q}$ by the closure of
the subspace generated by all commutators. It admits a topological
basis formed by the {\em cycles} of $Q$, \ie the orbits
of paths $p=(i|\alpha_m| \ldots |\alpha_1| i)$ of any length $m\geq 0$ with
identical source and target under the action of the cyclic group of order $m$.
In particular, the space $HH_0(\hat{\C Q})$ is a product of copies of $\C$ indexed
by the vertices if $Q$ does not have oriented cycles.
For each arrow $a$ of $Q$, the {\em cyclic derivative\index{cyclic derivative} with respect to $a$}
is the unique continuous $\C$-linear map
\[
\del_a : HH_0(\hat{\C Q}) \to \hat{\C Q}
\]
which takes the class of a path $p$ to the sum
\[
\sum_{p=uav} vu
\]
taken over all decompositions of $p$ as a concatenation
of paths $u$, $a$, $v$, where $u$ and $v$ are of length $\geq 0$.
A {\em potential}\index{potential} on $Q$ is an element $W$ of $HH_0(\hat{\C Q})$ whose expansion
in the basis of cycles does not involve cycles of length $\leq 1$.
A potential is {\em reduced}\index{potential!reduced} if it does not involve cycles of length
$\leq 2$. The {\em Jacobian algebra}\index{Jacobian algebra} $J(Q,W)$ associated to a quiver
$Q$ with potential $W$ is the quotient of the completed path algebra
by the closure of the $2$-sided ideal generated by the cyclic derivatives
of the elements of $W$. If the potential $W$ is reduced and the Jacobian
algebra $J(Q,W)$ is finite-dimensional, its quiver is isomorphic to $Q$.

As typical examples, we may consider the quiver $Q$
\begin{equation} \label{eq:cyclic-quiver}
\xymatrix{ & 2 \ar[dr]^a & \\
1 \ar[ur]^b & & 3 \ar[ll]^c }
\end{equation}
with the potential $W=abc$ or with the potential $W=(abc)^2$.

In order to define the mutation of a quiver with potential $(Q,W)$
at a vertex $k$, we need to recall the construction of a reduced
quiver with potential from an arbitrary quiver with potential.

Two quivers with potential $(Q,W)$ and $(Q',W')$
are {\em right equivalent}\index{right equivalent} if $Q_0=Q_0'$ and there exists a $\C$-algebra
isomorphism $\phi: \hat{kQ} \to \hat{kQ'}$ such that $\phi$ induces
the identity on the subalgebra $\prod_{Q_0} \C$ and the induced
map in topological Hochschild homology takes $W$ to $W'$.
A quiver with potential $(Q,W)$ is {\em trivial}\index{potential!trivial} if $W$ is
a (possibly infinite) linear combination of $2$-cycles
and $J(Q,W)$ is isomorphic to $\prod_{Q_0} \C$.
If $(Q,W)$ and $(Q',W')$ are two quivers with potential
such that the sets of vertices of $Q$ and $Q'$ coincide,
their {\em direct sum}\index{potential!direct sum of} $(Q,W)\oplus (Q',W')$ is defined
as the pair consisting of the quiver with the same
vertex set, with set of arrows the disjoint union of
those of $Q$ and $Q'$, and with the potential equal
to the sum $W \oplus W'$.

\begin{theorem}[\cite{DerksenWeymanZelevinsky08}, Theorem~4.6 and Proposition~4.5]
Any quiver with potential $(Q,W)$ is right equivalent to the direct
sum of a reduced one $(Q_{red},W_{red})$ and a trivial one
$(Q_{triv},W_{triv})$, both unique up to right equivalence. Moreover,
the inclusion induces an isomorphism from $J(Q_{red},W_{red})$ onto $J(Q,W)$.
\end{theorem}
The quiver with potential $(Q_{red}, W_{red})$ is
the {\em reduced part}\index{potential!reduced part of} of $(Q,W)$.

We can now define the mutation of a quiver with potential.
Let $(Q,W)$ be a quiver with potential such that $Q$ does not
have loops. Let $k$ be a vertex of $Q$
not lying on a $2$-cycle. The {\em mutation}\index{mutation!of a quiver with potential} $\mu_i(Q,W)$ is
defined as the reduced part of the quiver with potential
$\tilde{\mu}_i(Q,W)=(Q',W')$, which is defined as follows:
\begin{itemize}
\item[a)]
\begin{itemize}
\item[(i)] To obtain $Q'$ from $Q$, add a new arrow $[\alpha\beta]$
for each pair of arrows $\alpha:k\to j$ and $\beta: i \to k$ of $Q$ and
\item[(ii)] replace each arrow $\gamma$ with source or target $i$ by
a new arrow $\gamma^*$ with $s(\gamma^*)=t(\gamma)$ and $t(\gamma^*)=s(\gamma)$.
\end{itemize}
\item[b)] Put $W'=[W]+\Delta$, where
\begin{itemize}
\item[(i)] $[W]$ is obtained from $W$ by
replacing, in a representative of $W$ without cycles passing through
$k$, each occurrence of $\alpha\beta$ by $[\alpha\beta]$, for each pair
of arrows $\alpha:i\to k$ and $\beta: i \to k$ of $Q$;
\item[(ii)] $\Delta$ is the sum of the cycles $[\alpha\beta]\beta^* \alpha^*$
taken over all pairs of arrows $\alpha:k\to j$ and $\beta: i \to k$ of $Q$.
\end{itemize}
\end{itemize}
Then $k$ is not contained in a $2$-cycle of $\mu_k(Q,W)$ and
$\mu_k(\mu_k(Q,W))$ is right equivalent to $(Q,W)$ (Theorem~5.7 of \cite{DerksenWeymanZelevinsky08}).
As examples, consider the mutation at $2$ of the cyclic quiver (\ref{eq:cyclic-quiver})
endowed with the potential $W=abc$ and with $W'=(abc)^2$. For $W=abc$,
the mutated quiver with potential is the acyclic quiver
\begin{equation}
\xymatrix{ & 2 \ar[dl]_{b^*} & \\
1  & & 3 \ar[ul]_{a^*} }
\end{equation}
with the zero potential. But for $W=(abc)^2$, the mutated quiver with potential
is
\begin{equation} \label{eq:2-cycle-after-mutation}
\xymatrix{ & 2 \ar[dl]_{b^*} & \\
1 \ar@<0.5ex>[rr]^e & & 3 \ar[ul]_{a^*} \ar@<0.5ex>[ll]^c }
\end{equation}
with the potential $ecec+eb^* a^*$.

The general construction implies that if neither $Q$ nor the quiver $Q'$ in
$(Q',W')=\mu_k(Q,W)$ have loops or $2$-cycles, then $Q$ and $Q'$ are linked
by the quiver mutation rule (cf. Prop.~7.1 of \cite{DerksenWeymanZelevinsky08}) .
Thus, if we want to `extend' this rule to quivers with potentials, it is important
to ensure that no $2$-cycles appear during the mutation process.

Let $Q$ be a finite quiver.
A {\em continuous quotient} of $HH_0(\hat{\C Q})$ is
linear surjection $q: HH_0(\hat{\C Q}) \to V$ such that
for some $N\gg 0$, all potentials involving only cycles
of length $>N$ lie in the kernel of $q$. A {\em polynomial
function $HH_0(\hat{\C Q}) \to \C$} is the composition
of a continuous quotient $HH_0(\hat{\C Q}) \to V$ with
a polynomial map $V \to C$. A {\em hypersurface} in
$HH_0(\hat{\C Q})$ is the set of zeroes of a non zero
polynomial function.

\begin{theorem}[\cite{DerksenWeymanZelevinsky08}, Cor.~7.4] Let $Q$ be a
finite quiver without loops nor $2$-cycles. There is a countable
union of hypersurfaces $C\subset HH_0(\hat{\C Q})$ such that
for each $W$ not belonging to $C$, no $2$-cycles appear
in any iterated mutation of $(Q,W)$.
\end{theorem}

A potential $W$ not belonging to $C$ is called {\em generic}\index{potential!generic}.
So if $Q$ is a quiver without loops nor $2$-cycles and
$W$ a generic potential, we can indefinitely mutate the
quiver with potential $(Q,W)$ and the mutation of the
underlying quivers is given by the quiver mutation rule.
Notice that the potential $W=(abc)^2$ on the quiver~(\ref{eq:cyclic-quiver})
is not generic, which is compatible with the appearance of
a $2$-cycle in (\ref{eq:2-cycle-after-mutation}).

\subsection{Ginzburg algebras} \label{ss:Ginzburg-algebras}
Let $Q$ be a finite quiver and $W$ a potential on $Q$ (cf.~section~\ref{ss:quivers-with-potential}).
Let $\Gamma$ be the {\em Ginzburg \cite{Ginzburg06} dg algebra}\index{dg algebra!Ginzburg}\index{Ginzburg dg algebra}
of $(Q,W)$. It is constructed as follows: Let $\overline{Q}$ be the
graded quiver with the same vertices as $Q$ and whose arrows are
\begin{itemize}
\item the arrows of $Q$ (they all have degree~$0$),
\item an arrow $a^*: j \to i$ of degree $-1$ for each arrow $a:i\to j$ of $Q$,
\item a loop $t_i : i \to i$ of degree $-2$ for each vertex $i$
of $Q$.
\end{itemize}
The underlying graded algebra of $\Gamma(Q,W)$ is the
completion of the graded path algebra $\C\overline{Q}$ in the category
of graded vector spaces with respect to the ideal generated by the
arrows of $\overline{Q}$. Thus, the $n$-th component of
$\Gamma(Q,W)$ consists of elements of the form
$\sum_{p}\lambda_p p$, where $p$ runs over all paths of degree $n$.
The differential of $\Gamma(Q,W)$ is the unique continuous
linear endomorphism homogeneous of degree~$1$ which satisfies the
Leibniz rule
\[
d(uv)= (du) v + (-1)^p u dv \ko
\]
for all homogeneous $u$ of degree $p$ and all $v$, and takes the
following values on the arrows of $\overline{Q}$:
\begin{itemize}
\item $da=0$ for each arrow $a$ of $Q$,
\item $d(a^*) = \del_a W$ for each arrow $a$ of $Q$,
\item $d(t_i) = e_i (\sum_{a} [a,a^*]) e_i$ for each vertex $i$ of $Q$, where
$e_i$ is the lazy path at $i$ and the sum runs over the set of
arrows of $Q$.
\end{itemize}
One checks that $d^2=0$. For example, for the
the cyclic quiver~\ref{eq:cyclic-quiver}
with the potential $W=abc$, the  graded quiver $\overline{Q}$ is
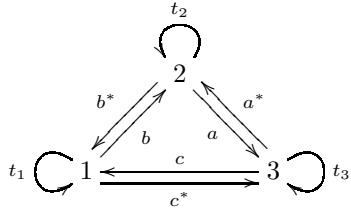
\[
\xymatrix{ & 2 \ar[dr]_a  \ar@<-1ex>[ld]_{b^*} \ar@(ur,ul)[]_{t_2} & \\
1 \ar@(ul,dl)[]_{t_1} \ar[ur]_b \ar@<-1ex>[rr]_{c^*}& &
3 \ar[ll]_c \ar@<-1ex>[ul]_{a^*} \ar@(ur,dr)[]^{t_3} }
\]
and the differential is given by
\[
d(a^*)=bc\ko d(b^*)=ca \ko d(c^*)=ab \ko d(t_1)= cc^* - b^* b \ko \ldots \ .
\]

The Ginzburg algebra should be viewed as a refined version of the
Jacobian algebra $J(Q,W)$. It is concentrated in (cohomological) degrees $\leq 0$
and $H^0(\Gamma)$ is isomorphic to $J(Q,W)$. Two right equivalent
quivers with potential have isomorphic Ginzburg algebras (Lemma~2.9 of \cite{KellerYang11}).
If $(Q,W)=(Q_{triv},W_{triv})\oplus (Q_{red},W_{red})$ is the direct
sum of a trivial and a reduced quiver with potential, then the
projection from $Q$ onto $Q_{red}$ induces a quasi-isomorphism
$\Gamma(Q,W) \to \Gamma(Q_{red},W_{red})$ (Lemma~2.10 of \cite{KellerYang11}).

\subsection{Derived categories of dg algebras} \label{ss:der-cat-dg-alg}
Let us recall the construction of the {\em derived category $\cd(A)$}\index{derived category!of a dg algebra}
of a dg (=differential graded) algebra $A$:
A \emph{(right) dg module}\index{dg module} $M$ over $A$ is a graded $A$-module
equipped with a differential $d$ such that
\[d(ma)=d(m)a+(-1)^{|m|}m d(a)\]
where $m$ in $M$ is homogeneous of degree $|m|$, and $a\in A$.

Given two dg $A$-modules $M$ and $N$, we define the \emph{morphism
complex} to be the graded $\C$-vector space $\mathcal{H}om_A(M,N)$
whose $i$-th component $\mathcal{H}om_A^i(M,N)$ is the subspace of
the product $\prod_{j\in\mathbb{Z}}\Hom_\C(M^j,N^{j+i})$ consisting
of the morphisms $f$ such that \[f(ma)=f(m)a,\] for all $m$ in $M$ and
all $a$ in $A$, together with the differential $d$ given by
\[d(f)=f\circ d_M-(-1)^{|f|}d_N\circ f\]
for a homogeneous morphism $f$ of degree $|f|$.

The {\em category $\mathcal{C}(A)$} of dg $A$-modules is the category
whose objects are the dg $A$-modules, and whose morphisms are the
0-cycles of the morphism complexes. This is an abelian category and
a Frobenius category for the conflations which are split exact as
sequences of graded $A$-modules. Its stable category
$\mathcal{H}(A)$ is called the \emph{homotopy category} of dg
$A$-modules, which is equivalently defined as the category whose
objects are the dg $A$-modules and whose morphism spaces are the
$0$-th homology groups of the morphism complexes. The homotopy
category $\mathcal{H}(A)$ is a triangulated category whose suspension
functor $\Sigma$ is the shift of dg modules $M \mapsto \Sigma M = M[1]$. The \emph{derived
category} $\mathcal{D}(A)$ of dg $A$-modules is the localization of
$\mathcal{H}(A)$ at the full subcategory of acyclic dg $A$-modules.
A short exact sequence
\[\xymatrix{0\ar[r]&M\ar[r]&N\ar[r]&L\ar[r]&0}\]
in $\mathcal{C}(A)$ yields a triangle
\[\xymatrix{M\ar[r]&N\ar[r]&L\ar[r]&\Sigma M}\]
in $\mathcal{D}(A)$.
A dg $A$-module $P$ is \emph{cofibrant} if
\[\Hom_{\mathcal{C}(A)}(P,L)\stackrel{s_*}{\longrightarrow}\Hom_{\mathcal{C}(A)}(P,M)\]
is surjective for each quasi-isomorphism $s:L\rightarrow M$ which is
surjective in each component. We use the term ``cofibrant'' since
these are actually the objects which are cofibrant for a certain
structure of Quillen model category on the category $\cc(A)$,
cf. \cite[Theorem 3.2]{Keller06d}. For an explicit description
of the cofibrant dg $A$-modules, cf. Prop.~2.13 of \cite{KellerYang11}.

The derived category $\cd(A)$ admits arbitrary (set-indexed) coproducts.
An object $P$ of $\cd(A)$ is {\em compact}\index{object!compact} or if the functor
$\Hom_{\cd(A)}(P,?):\cd(A) \to \cd(\C)$ commutes with arbitrary
coproducts. For example, the functor
\[
\Hom_{\cd(A)}(A,?) = H^0(?)
\]
commutes with arbitrary sums and so the free $A$-module of rank $1$
is compact. An arbitrary object of $\cd(A)$ is compact iff it
is {\em perfect}\index{object!perfect}, i.e if it belongs to the closure of $A$ under
left and right shifts, extensions and passage to direct factors.
The {\em perfect derived category $\per(A)\subset \cd(A)$}\index{derived category!perfect} is the full subcategory
on the perfect objects.

\subsection{The derived category of the Ginzburg algebra} \label{ss:der-cat-Ginzburg-alg}
As in section~\ref{ss:Ginzburg-algebras}, let $Q$ be a finite quiver. Assume
that the vertex set of $Q$ is $\{1, \ldots, n\}$. Let
$W$ be a reduced potential on $Q$. Let $\Gamma$ be the Ginzburg dg algebra
of the {\em opposite} quiver with potential $(Q^{op},W^{op})$.
Let $\cd(\Gamma)$\index{derived category!of a Ginzburg dg algebra} be the derived category of $\Gamma$
and $\per(\Gamma)$ the perfect derived category. By Lemma~2.17 of
\cite{KellerYang11}, the category $\per(\Gamma)$ is
a Krull-Schmidt category, i.e. each object decomposes into a
finite direct sum of indecomposable objects and each indecomposable
object has a local endomorphism algebra. In particular, the
free module of rank one $\Gamma \in \per(\Gamma)$ decomposes
into the indecomposable summands $P_i = e_i \Gamma$ associated
with the vertices $i$ of $Q$. The Grothendieck group
$K_0(\per(\Gamma))$ is free on the basis formed by the
classes $[P_i]$, $1\leq i\leq n$.

Now let $\cd_{fd}(\Gamma)$ the {\em finite-dimensional derived category
of $\Gamma$}\index{derived category!finite-dimensional}, i.e. the full subcategory of $\cd(\Gamma)$ formed by the dg modules
whose homology is of finite total dimension. An object $M$
belongs to $\cd_{fd}(\Gamma)$ if and only if for each object
$P$ of $\per(\Gamma)$, the space $\Hom_{\cd(\Gamma)}(P,\Sigma^i M)$
vanishes for almost all $i\in \Z$ and is finite-dimensional
for all $i\in\Z$.
The category $\cd_{fd}(\Gamma)$ is in fact contained in $\per(\Gamma)$.
It is triangulated and has finite-dimensional morphism spaces. More
precisely, for $L$ and $M$ in $\cd_{fd}(\Gamma)$, the spaces
$\Hom_{\cd(\Gamma)}(L, \Sigma^i M)$ are finite-dimensional for all
$i\in \Z$ and vanish for all but finitely many $i\in \Z$. Thus, the
Grothendieck group $K_0(\cd_{fd}(\Gamma))$ carries a well-defined
Euler form:
\[
\langle L, M \rangle = \sum_{p\in\Z} (-1)^p \dim \Hom_{\cd(\Gamma)}(L, \Sigma^p M).
\]
For a vertex $i$ of $Q$, the simple $\overline{Q}^{op}$-representation
$S_i$ concentrated at the vertex $i$ yields a simple dg $\Gamma$-module
still denoted by $S_i$. The $S_i$ generate the triangulated
category $\cd_{fd}(\Gamma)$. The Grothendieck group $K_0(\cd_{fd}(\Gamma))$
is free on the basis given by the classes $[S_i]$, $1\leq i\leq n$.
We have a well-defined non degenerate pairing
\[
\langle\,, \rangle : K_0(\per(\Gamma)) \times K_0(\cd_{fd}(\Gamma)) \to \Z
\]
given again as a Euler form
\[
\langle P, M \rangle = \sum_{p\in\Z} (-1)^p \dim \Hom_{\cd(\Gamma)}(P, \Sigma^p M).
\]
We have
\[
\langle P_i, S_j \rangle = \delta_{ij} \ko 1\leq i,j \leq n\ko
\]
so that the basis of the $[S_j]$, $1\leq j\leq n$, is dual to that of
the $[P_i]$, $1\leq i\leq n$.

Let $j$ be a vertex of $Q$. It follows from the cofibrant resolution of $S_j$ given at the
beginning of the proof of Lemma~3.12 in \cite{KellerYang11} that
the image of the class of $S_j$ in $K_0(\per(\Gamma))$ equals
\[
[P_j] - \sum_{\alpha: t(\alpha)=j} [P_{s(\alpha)}] + \sum_{\beta: s(\beta)=j} [P_{t(\alpha)}] - [P_j]
= \sum_{i} b_{ij}\, [P_i] \ko
\]
where the source and target maps refer to $Q^{op}$ and $B=(b_{ij})$ is the
antisymmetric matrix associated with the quiver $Q$. Thus, the matrix of the
linear map
\[
K_0(\cd_{fd}(\Gamma)) \to K_0(\per(\Gamma))
\]
in the bases of the $[S_j]$ and the $[P_i]$ is $B$.
It follows that we have
\[
\langle S_i, S_j \rangle = -b_{ij}
\]
so that $-B$ is the matrix of the Euler form $\langle\, , \rangle$ on
$K_0(\cd_{fd}(\Gamma))$ in the basis of the $[S_i]$, $1\leq i\leq n$.

The category $\cd_{fd}(\Gamma)$ is $3$-Calabi-Yau, by which
we mean that we have bifunctorial isomorphisms
\[
D\Hom(X,Y) \iso \Hom(Y,\Sigma^3 X) \ko
\]
where $D$ is the duality functor $\Hom_\C(?,\C)$ and $\Sigma$ the
shift functor. The simple modules $S_i$ are $3$-spherical objects in $\cd_{fd}(\Gamma)$,
i.e. we have an isomorphism
\[
\Ext_\Gamma^*(S_i, S_i) \cong H^*(S^3,\C)
\]
where the left hand side denotes the graded vector space
whose $p$th component is $\Hom_{\cd(\Gamma)}(S_i, \Sigma^p S_i)$ and
the right hand side is the (singular) cohomology of the $3$-sphere with complex
coefficients. The spherical objects $S_i$ yield the Seidel-Thomas \cite{SeidelThomas01} twist functors
$\tw_{S_i}$. These are autoequivalences
of $\cd(\Gamma)$ such that each object $X$ fits into a triangle
\[
\RHom(S_i,X) \ten_k S_i \to X \to \tw_{S_i}(X) \to \Sigma \RHom(S_i,X) \ten_k S_i\;.
\]
By \cite{SeidelThomas01}, the twist functors give rise to a (weak) action on $\cd(\Gamma)$ of the
braid group associated with $Q$, i.e.~the group with generators
$\sigma_i$, $i\in Q_0$, and relations $\sigma_i \sigma_j = \sigma_j \sigma_i$
if $i$ and $j$ are not linked by an arrow and
\[
\sigma_i \sigma_j \sigma_i = \sigma_j \sigma_i \sigma_j
\]
if there is exactly one arrow between $i$ and $j$ (no relation
if there are two or more arrows).

The category $\cd(\Gamma)$ admits a natural $t$-structure
whose truncation functors are those of the natural $t$-structure
on the category of complexes of vector spaces (because $\Gamma$
is concentrated in degrees $\leq 0$). Thus, we have an
induced natural $t$-structure on $\cd_{fd}(\Gamma)$. Its
heart $\ca$ is canonically equivalent to the category
$\nil(J(Q,W))$ of finite-dimensional right modules over
$J(Q,W)$ where all sufficiently long paths act by $0$.
In particular, the inclusion of $\ca$ into $\cd_{fd}(\Gamma)$
induces an isomorphism in the Grothendieck groups
\[
K_0(\ca) \iso K_0(\cd_{fd}(\Gamma)).
\]

The {\em cluster category}\index{cluster category} is the triangle
quotient
\begin{equation} \label{eq:def-cluster-category}
\cc_\Gamma = \per(\Gamma)/\cd_{fd}(\Gamma)
\end{equation}
For acyclic quivers $Q$, Amiot \cite{Amiot09} has shown that
it is equivalent to the cluster category $\cc_{Q^{op}}$ (we pass
to the opposite because $\Gamma$ is associated with $Q^{op}$) in
the sense of \cite{BuanMarshReinekeReitenTodorov06}. For
arbitrary quivers, there is also a close link between $\per(\Gamma)$
and $\cc_{\Gamma}$: For a triangulated category $\ct$ and an object $X$ of $\ct$,
let us denote by $\pr_\ct(X)$ the subcategory of {\em $X$-presentable objects
of $\ct$}\index{object!$X$-presentable}, i.e. the objects $Y$ which occur
in a triangle
\begin{equation} \label{eq:presentation-triangle}
X'' \to X' \to Y \to \Sigma X'' \ko
\end{equation}
where $X''$ and $X'$ belong to the closure $\add(X)$ of $X$ under
taking (finite) direct sums and direct summands.

\begin{proposition}[Prop.~2.10 of \cite{Plamondon11}] \label{prop:equiv-presentable}
The projection
$\per(\Gamma) \to \cc_{\Gamma}$ induces a $\C$-linear
equivalence
\[
\pr_{\per(\Gamma)}(\Gamma) \iso \pr_{\cc_{\Gamma}}(\Gamma) \ko
\]
\end{proposition}

Plamondon also relates the extension groups computed in the
two categories (Prop.~2.19 of \cite{Plamondon11}).

\subsection{Derived equivalences from mutations}\index{derived equivalence}
\label{ss:der-eq}
As in the preceding section, let $Q$ be a finite quiver
without loops nor $2$-cycles with vertex set $\{1, \ldots, n\}$
and let $W$ be a generic potential on $Q$. Let $\Gamma$
denote the Ginzburg algebra associated with the
{\em opposite} quiver with potential $(Q^{op}, W^{op})$.
Let $k$ be a vertex of $Q$ and $\Gamma'$ the Ginzburg
algebra associated with the opposite of the mutated
quiver with potential $\mu_k(Q,W)$. Put
$P_i=e_i \Gamma$ and $P_i' = e_i \Gamma'$, $1\leq i\leq n$.

\begin{theorem}[\cite{KellerYang11}] \label{thm:equiv-Ginzburg}
There are two canonical equivalences
\[
\Phi_{\pm}: \cd(\Gamma') \to \cd(\Gamma)
\]
related by an isomorphism
\[
\tw_{S_k} \circ \Phi_- \iso \Phi_+ \;.
\]
Both $\Phi_+$ and $\Phi_-$ send $P'_i$ to $P_i$ for $i\neq k$ and
the images of $P'_k$ under the two functors fit into triangles
\begin{equation} \label{eq:tria-plus}
\xymatrix{P_k \ar[r] & \bigoplus_{k\to i} P_i \ar[r] &  \Phi_-(P'_k) \ar[r] &  \Sigma P_k}
\end{equation}
and
\begin{equation} \label{eq:tria-minus}
\xymatrix{\Sigma^{-1} P_k \ar[r] & \Phi_+(P'_k) \ar[r] & \bigoplus_{j\to k} P_j \ar[r] &  P_k} \ko
\end{equation}
where the sums are taken over the arrows in $Q^{op}$.
\end{theorem}

The intrinsic characterizations of the subcategories $\per(\Gamma)$ and
$\cd_{fd}(\Gamma)$ show that the equivalences $\Phi_{\pm}$ induce
equivalences
\[
\per(\Gamma') \to \per(\Gamma) \quad\mbox{and}\quad \cd_{fd}(\Gamma')\to \cd_{fd}(\Gamma).
\]
and thus isomorphisms in the associated Grothendieck groups. By the
triangles~(\ref{eq:tria-plus}) and (\ref{eq:tria-minus}), we get the
first statement of the following corollary; the second one follows
by passage to the duals. We use the matrices $E_\eps$ and $F_\eps$
associated with $Q$ in section~\ref{ss:quantum-mutations}
(remember however that $\Gamma$ is the Ginzburg algebra
of the opposite of $(Q,W)$).

\begin{corollary} \label{cor:categorification-of-E-eps}
Let $\eps$ be equal to $1$ or $-1$.
Under the assumptions of the theorem, the matrix of
the induced isomorphism
\[
K_0(\Phi_\eps): K_0(\per(\Gamma')) \to K_0(\per(\Gamma))
\]
in the bases $[P'_j]$ and $[P_i]$ is $E_\eps$ and the
matrix of the induced isomorphism
\[
K_0(\Phi_\eps): K_0(\cd_{fd}(\Gamma')) \to K_0(\per(\Gamma))
\]
in the bases $[S'_j]$ and $[S_i]$ is $F_\eps$.
\end{corollary}

\begin{figure}
\begin{center}
\begin{tikzpicture}[scale=0.45]
\draw (0,4)--(13,4);
\draw (0,9)--(13,9);
\draw (1.5,4)--(1.5,9);
\draw (4,4)--(4,9);
\draw (9.25,4)--(9.25,9);
\draw (11.75,4)--(11.75,9);
\draw (2.75,6.5) node {$S_k$};
\draw (6.5,6.5) node {$S_k^\perp$};
\draw (10.5,6.5) node {$\Sigma^{-1}S_k$};
\draw [decorate, decoration={brace,amplitude=5pt}] (1.5,9.5)--(9.25,9.5)
node [midway, above=6pt] {$\ca$};
\draw [decorate, decoration={brace,amplitude=5pt}] (11.75,3.5)--(4, 3.5)
node [midway, below=6pt] {$\mu_k^+(\ca)$};
\end{tikzpicture}
\quad\quad
\begin{tikzpicture}[scale=0.45]
\draw (0,4)--(13,4);
\draw (0,9)--(13,9);
\draw (1.5,4)--(1.5,9);
\draw (4,4)--(4,9);
\draw (9.25,4)--(9.25,9);
\draw (11.75,4)--(11.75,9);
\draw (2.75,6.5) node {$\Sigma S_k$};
\draw (6.5,6.5) node {$^\perp\!S_k$};
\draw (10.5,6.5) node {$S_k$};
\draw [decorate, decoration={brace,amplitude=5pt}] (4,9.5)--(11.75,9.5)
node [midway, above=6pt] {$\ca$};
\draw [decorate, decoration={brace,amplitude=5pt}] (9.25,3.5)--(1.5, 3.5)
node [midway, below=6pt] {$\mu_k^-(\ca)$};
\end{tikzpicture}
\caption{Right and left mutation of a heart}\label{fig:mut-hearts}
\end{center}
\end{figure}
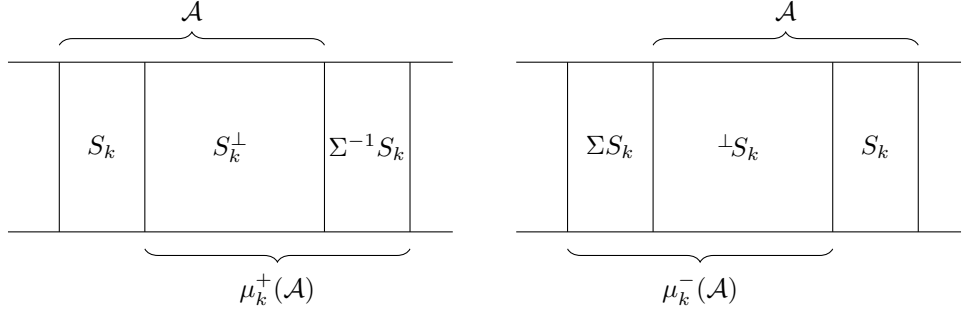

Let $\ca'$ be the heart of the canonical $t$-structure on $\cd_{fd}(\Gamma')$.
The equivalences $\Phi_{\pm}$ send $\ca'$ onto the hearts $\mu_k^{\pm}(\ca)$ of
two new t-structures. These can be described in terms of
$\ca$ and the subcategory $\add S_k$ as follows
(cf. figure~\ref{fig:mut-hearts}): Let
$S_k^\perp$ be the right orthogonal subcategory of $S_k$ in $\ca$,
whose objects are the $M$ with $\Hom(S_k,M)=0$. Then
$\mu_k^+(\ca)$ is formed by the objects $X$ of $\cd_{fd}(\Gamma)$
such that the object $H^0(X)$ belongs to $S_k^\perp$, the object $H^1(X)$ belongs
to $\add S_k$ and $H^p(X)$ vanishes for all
$p\neq 0,1$. Similarly, the subcategory $\mu_k^-(\ca)$ is
formed by the objects $X$ such that the object $H^0(X)$ belongs to
the left orthogonal subcategory $^\perp\!S_k$, the object $H^{-1}(X)$
belongs to $\add(S_k)$ and $H^p(X)$ vanishes for all
$p \neq -1,0$. The subcategory $\mu_k^{+}(\ca)$ is
the {\em right mutation of $\ca$}\index{mutation!of a heart}\index{mutation!right} and $\mu_k^-(\ca)$ is
its {\em left mutation}\index{mutation!left}.

% \begin{figure}
% \begin{center}
% \includegraphics[height=5cm]{mut-hearts.png}
% \end{center}
% \end{figure}
By construction, we have
\[
\tw_{S_k}(\mu_k^-(\ca)) = \mu_k^+(\ca).
\]
Since the categories $\ca$ and $\mu^\pm(\ca)$ are hearts of
bounded, non degenerate $t$-structures on $\cd_{fd}(\Gamma)$, their
Grothendieck groups identify canonically with that of
$\cd_{fd}(\Gamma)$. They are endowed with canonical bases
given by the simples. Those of $\ca$ identify with
the simples $S_i$, $i\in Q_0$, of $\nil(J(Q,W))$.
The simples of $\mu_k^+(\ca)$ are $\Sigma^{-1} S_k$,
the simples $S_i$ of $\ca$ such that $\Ext^1(S_k, S_i)$
vanishes and the objects $\tw_{S_k}(S_i)$ where
$\Ext^1(S_k, S_i)$ is of dimension $\geq 1$. By applying
$\tw_{S_k}^{-1}$ to these objects we obtain the
simples of $\mu_k^-(\ca)$.

\subsection{Torsion subcategories and intermediate $t$-structures}
\label{ss:torsion-subcategories}
In order to investigate the effect on hearts of suitable compositions
of the equivalences $\Phi_\pm$ of Theorem~\ref{thm:equiv-Ginzburg},
let us recall the construction of `tilted hearts', a variation
on a construction of \cite{HappelReitenSmaloe96}. Let $\cd$
be a triangulated category (for example the category $\cd_{fd}(\Gamma)$).
Let $(\cd_{\leq 0}, \cd_{\geq 0})$ be a bounded non degenerate
$t$-structure on $\cd$ and $\ca$ its heart. A {\em torsion pair\index{torsion pair} on $\ca$}
is a pair $(\ct, \cf)$ of full subcategories such that
\begin{itemize}
\item[a)] we have $\Hom(T,F)=0$ for all $T\in \ct$ and $F\in\cf$ and
\item[b)] for each object $M$ of $\ca$, there is a short exact sequence
\[
\xymatrix{ 0 \ar[r] & M_\ct \ar[r] & M \ar[r] & M^\cf \ar[r] & 0}
\]
with $M_\ct$ in $\ct$ and $M^\cf$ in $\cf$.
\end{itemize}
In this case, the subcategories $\ct$ and $\cf$ determine each other:
we have $\cf=\ct^\perp$ and $\ct = \mbox{ }^\perp\!\ct$, where the
orthogonal subcategories are taken in $\ca$.

For two full subcategories $\cu$ and $\cv$ of $\cd$, let us write
$\cu\star\cv$ for the full subcategory whose objects $X$ occur in
triangles
\[
U \to X \to V \to \Sigma X
\]
with $U$ in $\cu$ and $V$ in $\cv$. Let $(\ct,\cf)$ be a torsion
par in $\ca$. Then the subcategory $\cd_{\leq 0} \star (\Sigma^{-1} \cf)$ is the
left aisle of a new $t$-structure, whose heart $\ca(\cf, \Sigma^{-1}\ct)$
equals $\cf\star\Sigma^{-1}\ct$. It is called the {\em right tilt}\index{tilt!right} of $\ca$
at $(\ct,\cf)$. Dually, the subcategory $(\Sigma \cf)\star\cd_{\geq 0}$ is
the right aisle of a new $t$-structure on $\cd$, whose heart $\ca(\Sigma \cf, \ct)$
equals $(\Sigma\cf)\star\ct$. It is called the {\em left tilt}\index{tilt!left} of $\ca$
at $(\ct,\cf)$. The right tilt $\ca(\cf, \Sigma^{-1}\ct)$ admits the
torsion pair $(\cf,\Sigma^{-1}\ct)$ and its left tilt with respect
to this pair equals the original category $\ca = \ct\star\cf$. Similarly,
the left tilt $\ca (\Sigma \cf, \ct)$ admits the torsion pair $(\Sigma \cf, \ct)$
and we recover $\ca$ as its right tilt with respect to this pair.

Clearly, the left aisle $\cd_{\leq 0}\star(\Sigma^{-1}\cf)$ is an
{\em intermediate left aisle}\index{aisle!intermediate}, i.e.~we have
\[
\cd_{\leq 0} \subset \cd_{\leq 0}\star(\Sigma^{-1}\cf) \subset \cd_{\leq 1}.
\]
It is not hard to check that each intermediate left aisle is of this form.
Dually, each right aisle between $\cd_{\geq -1}$ and $\cd_{\geq 0}$ is of
the form $(\Sigma \cf)\star \cd_{\geq 0}$.

Of course, in the situation of section~\ref{ss:der-eq}, the heart
$\mu_k^+(\ca)$ is the right tilt of $\ca$ with respect to
$(\add S_k, S_k^\perp)$ and $\mu_k^-(\ca)$ is its left tilt.
The following lemma will allow us to iterate such mutations.
The point is that given a suitable torsion pair in a right tilt
$\ca'$ of $\ca$, either the left tilt or the right tilt of $\ca'$
with respect to this pair is again a right tilt of $\ca$,
cf. figure~\ref{fig:compositions-right-tilt-with-left-tilt}.

\begin{lemma}[Bridgeland \cite{Nagao10}] \label{lemma:composition-of-tilts}
Let $(\ct,\cf)$ be a torsion pair in $\ca$ and $(\ct',\cf')$ a torsion pair
in $\ca'=\ca(\cf,\Sigma^{-1}\ct)$.
\begin{itemize}
\item[a)] If $\ct'\subset \cf$, then  the right tilt $\ca''=\ca'(\cf', \Sigma^{-1} \ct')$
equals the right tilt
\[
\ca(\ct\star\ct', \cf\cap\cf').
\]
\item[b)] If $\cf' \subset \Sigma^{-1}\ct$, then the left tilt $\ca''=\ca'(\Sigma \ct', \cf')$
equals the right tilt
\[
\ca(\Sigma \cf'\star\ct, \ct'\cap \Sigma^{-1}\ct).
\]
\end{itemize}
\end{lemma}

The lemma is not hard to check. The following easy lemma is
a key point for the main conjecture~\ref{conj:main}:

\begin{figure}
\begin{center}
\begin{tikzpicture}[scale=0.4]
%\draw [very thin, gray] (0,0) grid (25,15);
\draw (1,5)--(21,5);
\draw (1,10)--(21,10);
\draw [fill=gray!40] (4,5)--(7,5)--(7,10)--(4,10)--cycle;
\draw [fill=gray!20] (7,5)--(9,5)--(9,10)--(7,10)--cycle;
\draw [fill=gray!40] (13,5)--(16,5)--(16,10)--(13,10)--cycle;
\draw [fill=gray!20] (16,5)--(18,5)--(18,10)--(16,10)--cycle;
\draw [decorate, decoration={brace,amplitude=5pt}] (4,10.5)--(6.9,10.5)
node [midway, above=6pt] {$\ct$};
\draw [decorate, decoration={brace,amplitude=5pt}] (7.1,10.5)--(12.9,10.5)
node [midway, above=6pt] {$\cf$};
\draw [decorate, decoration={brace,amplitude=5pt}] (13.1,10.5)--(16,10.5)
node [midway, above=6pt] {$\Sigma^{-1}\ct$};
% \draw [decorate, decoration={brace,amplitude=5pt}] (7,12)--(16,12)
% node [midway, above=6pt] {$\cb=\ca(\cf,\Sigma^{-1}\ct)$};
\draw [decorate, decoration={brace,amplitude=5pt}] (8.9,4.5)--(7,4.5)
node [midway, below=6pt] {$\ct'$};
\draw [decorate, decoration={brace,amplitude=5pt}] (15.9,4.5)--(9.1,4.5)
node [midway, below=6pt] {$\cf'$};
\draw [decorate, decoration={brace,amplitude=5pt}] (17.9,4.5)--(16.1,4.5)
node [midway, below=6pt] {$\Sigma^{-1}\ct'$};
% \draw [decorate, decoration={brace,amplitude=5pt}] (18,3)--(9,3)
% node [midway, below=6pt] {$\cb(\cf',\Sigma^{-1}\ct')=\ca(\ct*\ct',\cf\cap\cf')$};
\draw[<->] (4,8.5)--(13,8.5) node [midway,fill=gray!20] {\makebox[0.5ex]{$\ca$}};
\draw[<->] (7,7.5)--(16,7.5) node [midway,fill=white] {\makebox[1ex]{$\ca'$}};
\draw[<->] (9,6.5)--(18,6.5) node [midway,fill=gray!40] {\makebox[1ex]{$\ca''$}};
\end{tikzpicture}
\makebox[6cm]{\mbox{ }}
\begin{tikzpicture}[scale=0.4]
%\draw [very thin, gray] (0,0) grid (25,15);
\draw (1,5)--(21,5);
\draw (1,10)--(21,10);
\draw [fill=gray!40] (4,5)--(7,5)--(7,10)--(4,10)--cycle;
\draw [fill=gray!20] (7,5)--(9,5)--(9,10)--(7,10)--cycle;
\draw [fill=gray!40] (13,5)--(16,5)--(16,10)--(13,10)--cycle;
\draw [fill=gray!20] (16,5)--(18,5)--(18,10)--(16,10)--cycle;
\draw [decorate, decoration={brace,amplitude=5pt}] (4,10.5)--(8.9,10.5)
node [midway, above=6pt] {$\ct$};
\draw [decorate, decoration={brace,amplitude=5pt}] (9.1,10.5)--(12.9,10.5)
node [midway, above=6pt] {$\cf$};
\draw [decorate, decoration={brace,amplitude=5pt}] (13.1,10.5)--(18,10.5)
node [midway, above=6pt] {$\Sigma^{-1}\ct$};
% \draw [decorate, decoration={brace,amplitude=5pt}] (7,12)--(16,12)
% node [midway, above=6pt] {$\cb=\ca(\cf,\Sigma^{-1}\ct)$};
\draw [decorate, decoration={brace,amplitude=5pt}] (8.9,4.5)--(7,4.5)
node [midway, below=6pt] {$\Sigma^{-1}\cf'$};
\draw [decorate, decoration={brace,amplitude=5pt}] (15.9,4.5)--(9.1,4.5)
node [midway, below=6pt] {$\ct'$};
\draw [decorate, decoration={brace,amplitude=5pt}] (17.9,4.5)--(16.1,4.5)
node [midway, below=6pt] {$\cf'$};
% \draw [decorate, decoration={brace,amplitude=5pt}] (18,3)--(9,3)
% node [midway, below=6pt] {$\cb(\cf',\Sigma^{-1}\ct')=\ca(\ct*\ct',\cf\cap\cf')$};
\draw[<->] (4,8.5)--(13,8.5) node [midway,fill=gray!20] {\makebox[0.5ex]{$\ca$}};
\draw[<->] (9,7.5)--(18,7.5) node [midway,fill=gray!40] {\makebox[1ex]{$\ca'$}};
\draw[<->] (7,6.5)--(16,6.5) node [midway,fill=white] {\makebox[1ex]{$\ca''$}};
\end{tikzpicture}
\end{center}
\caption{Composition of a right tilt with a right tilt (top)  resp. a left tilt (bottom)}
\label{fig:compositions-right-tilt-with-left-tilt}
\end{figure}
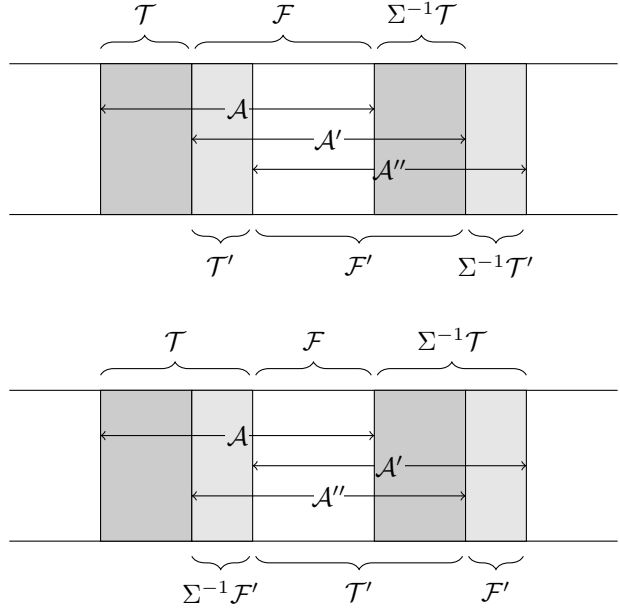

\begin{lemma} \label{lemma:dichotomy-for-simples}
Let $\ca'=\ca(\cf, \Sigma^{-1}\ct)$ be the right
tilt of $\ca$ with respect to a torsion pair $(\ct,\cf)$. Then
each simple object of $\ca'$ either lies in $\ca$ or in $\Sigma^{-1}\ca$.
\end{lemma}

Indeed, let $S$ be the given simple object. Since $(\cf, \Sigma^{-1}\ct)$
forms a torsion pair in $\ca'$, we have the exact sequence
\[
\xymatrix{0 \ar[r] & S_{\cf} \ar[r] & S \ar[r] & S^{\Sigma^{-1}\ct} \ar[r] & 0}
\]
where $S_\cf$ belongs to $\cf\subset\ca$ and $S^{\Sigma^{-1}\ct}$ to $\Sigma^{-1}\cf\subset \Sigma^{-1}\ca$.
Since $S$ is simple, we either have $S_\cf \iso S$ or $S\iso S^{\Sigma^{-1}\ct}$.

Now let $(Q,W)$ and $(Q',W')$ be two quivers with reduced potentials
and let $\Gamma$ and $\Gamma'$ be the associated Ginzburg dg algebras.
Suppose that
\[
\Phi: \cd(\Gamma') \to \cd(\Gamma)
\]
is a triangle equivalence. Let $(\cd_{\leq 0}, \cd_{\geq 0})$ be
the natural $t$-structure on $\cd_{fd}(\Gamma)$ and let
$\ca$ be its heart. Similarly, let $(\cd'_{\leq 0}, \cd'_{\geq 0})$ be the
natural $t$-structure on $\cd_{fd}(\Gamma')$
and let $\ca'$ be its heart. Let us denote
by $H^p_\ca$, $p\in\Z$, the homology functors with respect to the natural
$t$-structure on $\cd_{fd}(\Gamma)$.

\begin{proposition} \label{prop:width-one-equivalence}
The following are equivalent:
\begin{itemize}
\item[(i)] the subcategory $\Phi(\ca')\subset \cd_{fd}(\Gamma)$ is the right
tilt of $\ca$ with respect to a torsion pair $(\ct,\cf)$;
\item[(ii)] the object $\Phi(\Gamma')$ is $(\Sigma^{-1}\Gamma)$-presentable
(cf. the end of section~\ref{ss:der-cat-Ginzburg-alg});
\item[(iii)] we have $\cd_{\leq 0} \subset \Phi(\cd'_{\leq 0}) \subset \cd_{\leq 1}$.
\end{itemize}
If these conditions hold, then the torsion subcategory $\ct$ of (i)
is formed by the finite-dimensional quotients of the object $H^1_{\ca}(\Phi(\Gamma'))$.
\end{proposition}

\subsection{Patterns of tilts and decategorification} \label{ss:patterns}\index{pattern!of tilts}
As in section~\ref{ss:der-eq}, let $Q$ be a finite quiver
without loops nor $2$-cycles with vertex set $\{1, \ldots, n\}$
and let $W$ be a generic potential on $Q$. Let $\Gamma$
denote the Ginzburg algebra associated with the
{\em opposite} quiver with potential $(Q^{op}, W^{op})$.

As we have seen in section~\ref{ss:quivers-with-potential},
we can indefinitely mutate $(Q,W)$. Thus, with each vertex
$t$ of the $n$-regular tree, we can associate a quiver
with potential $(Q(t),W(t))$ such that $(Q,W)$ is
associated with $t_0$ and, whenever $t$ and $t'$ are linked
by an edge labeled $k$, the corresponding quivers with potential
are linked by a mutation. We write $\Gamma(t)$ for the
Ginzburg dg algebra associated with the {\em opposite}
of $(Q(t),W(t))$.

Now we will use induction on the distance of a vertex $t$ of
$n$-regular tree from the root $t_0$ to define a
triangle equivalence
\[
\Phi(t): \xymatrix{\cd(\Gamma(t)) \ar[r] & \cd(\Gamma)}
\]
such that $\Phi(t)$ satisfies the equivalent conditions of
Proposition~\ref{prop:width-one-equivalence}. The construction
follows an idea of Bridgeland \cite{Nagao10}.
By definition, $\Phi(t_0)$ is the identity. Now suppose that
$\Phi(t)$ has been defined and that $t'$ is linked to $t$
by an edge labeled $k$ and is at greater distance from $t_0$
than $t$. Let $\ca(t) \subset \cd_{fd}(\Gamma)$ be the image under
$\Phi(t)$ of the heart of the natural $t$-structure of
$\cd_{fd}(\Gamma(t))$ and let $S_i(t)$, $1\leq i\leq n$, be
the simple objects of $\ca(t)$. By assumption, the subcategory
$\ca(t)$ is the right tilt of $\ca=\ca(t_0)$ with respect to some
torsion theory $(\ct(t), \cf(t))$. Thus, by Lemma~\ref{lemma:dichotomy-for-simples},
the simple object $S_k(t)$ either lies in $\cf(t)\subset \ca$
or in $\Sigma^{-1}\ct(t) \subset \Sigma^{-1}\ca$. In the first
case, we put
\[
\Phi(t') = \Phi(t)\circ \Phi_{k,+} \mbox{ and in the second case }
\Phi(t') = \Phi(t)\circ \Phi_{k,-} \, .
\]
Then in the first case, as we have seen in section~\ref{ss:der-eq},
$\ca(t')$ is the right tilt of $\ca(t)$ with
respect the torsion pair
\[
(\add(S_k(t)), S_k(t)^\perp)
\]
and in the second case, it is the left tilt with respect to
\[
(\mbox{}^\perp S_k(t), \add(S_k(t))).
\]
In both cases, Lemma~\ref{lemma:composition-of-tilts} shows that
$\ca(t')$ is again a right tilt of $\ca$ and so $\Phi(t')$
again satisfies the conditions of Proposition~\ref{prop:width-one-equivalence}.

Notice that at the same time, this construction produces a
sign $\eps(e)$ for each edge $e: t \to t'$ of the $n$-regular tree.
For each vertex $t$ of $\mathbb{T}_n$,
and for $1\leq i\leq n$, let $T_i(t)$ be the image of $e_i \Gamma(t)$
under $\Phi(t)$.

\begin{theorem} \label{thm:decategorification} Let $t$ be a vertex of the $n$-regular tree and
let $1\leq j\leq n$.
\begin{itemize}
\item[a)] The $j$th column of the $c$-matrix $C(t)$ contains the
coordinates of $[S_j(t)]$ in the basis $[S_1]$, \ldots, $[S_n]$ of
$K_0(\cd_{fd}(\Gamma))$.
\item[b)] The object $S_j(t)$ lies in $\ca$ or $\Sigma^{-1}\ca$.
Therefore, each $c$-vector is non zero and has either all its components
non negative or all its components non positive (i.e. the main conjecture~\ref{conj:main}
holds for $Q$).
\item[c)] The $j$th column of the $g$-matrix $G(t)$ contains
the coordinates of $[T_j(t)]$ in the basis $[P_1]$, \ldots, $[P_n]$
of $K_0(\per(\Gamma))$.
\item[d)] The (left) $J(Q,W)$-module $H^1_\ca(T_j(t))$ is finite-dimensional
and the $F$-polyno\-mial $F_j(t)$ equals
\[
\sum_{e} \chi(\Gr_e(H^1_\ca(T_j(t))))\, y^e \ko
\]
where $e$ runs through $\N^n$, $\Gr_e$ denotes the variety of
submodules whose quotient has dimension vector $e$, $\chi$ is
the Euler characteristic (with respect to singular cohomology
with rational coefficients) and
\[
y^e = \prod_{i=1}^n x_{n+i}^{e_i}.
\]
\end{itemize}
\end{theorem}

To make sure that our conventions are coherent, let us consider
the example of the quiver $Q:1 \to 2$ and the vertex $t$ linked
to $t_0$ by the mutation at $1$. We have to consider the
Ginzburg algebra $\Gamma$ associated with $Q^{op}: 2 \to 1$ and perform
a right mutation at the vertex $1$. We get the exchange triangle
\begin{equation} \label{eq:exchange-triangle}
\Sigma^{-1} P_1 \to T_1(t) \to P_2 \to P_1  .
\end{equation}
Thus, the class of $T_1(t)$ in $K_0(\per(\Gamma))$ equals $-[P_1]+[P_2]$,
which does correspond to the $g$-vector
\[
g_1(t) = \left[\begin{array}{c} -1 \\ 1 \end{array} \right].
\]
The new simple modules are $S_1(t) = \Sigma^{-1} S_1$ and
$S_2(t)$ given by the universal extension
\[
S_2 \to S_2(t) \to S_1 \to \Sigma S_2.
\]
So in $K_0(\cd_{fd}(\Gamma))$, we have $[S_1(t)]=-[S_1]$ and
$[S_2(t)]=[S_1]+[S_2]$, which does correspond to the $c$-matrix
\[
C(t) = \left[ \begin{array}{cc} -1 & 1 \\ 0 & 1 \end{array} \right].
\]
Using the exchange triangle~(\ref{eq:exchange-triangle}), we easily
check that $\Hom(\Sigma^{-1}P_1, T_1(t))=\C$ and $\Hom(\Sigma^{-1}P_2, T_1(t))=0$
so that the module $H^1_\ca(T_1(t))$ is the simple at the vertex~$1$.
The associated generating series of Euler characteristics is indeed
equal to
\[
F_1(t)= 1+ y^{e_1}.
\]

Notice that since each $S_j(t)$ belongs to $\ca$ or $\Sigma^{-1}\ca$
(Lemma~\ref{lemma:dichotomy-for-simples}), part a) implies part b).
Thanks to parts a) and c), the duality between the bases formed
by the $[T_i(t)]$ and the $[S_j(t)]$ corresponds to the first part of the
tropical duality theorem~\ref{ss:tropical-duality}.

Parts a), b) and c) are proved in Nagao's \cite{Nagao10}, and part d)
is proved there under an additional technical assumption.
Parts a), b), c) and d) follow from the results of Plamondon \cite{Plamondon11a},
cf. section~\ref{ss:proof-of-decategorification} (and when $H^0(\Gamma)$ is
finite-dimensional from \cite{Palu08a}).
Using his dictionary between objects of the cluster category and decorated
representations, the theorem also can also be deduced from the results of
\cite{DerksenWeymanZelevinsky10}.

For acyclic quivers $Q$, part d) was extended to the
quantum case by Qin \cite{Qin10} and (for prime powers $q$)
by Rupel \cite{Rupel10a} \cite{Rupel11}, who also obtained
an analogous result for acyclic valued quivers.
Under certain technical assumptions, Efimov \cite{Efimov11}
has recently been able to extend part d) to the quantum
case  for arbitrary quivers (without loops nor $2$-cycles).
He mainly builds on the work of
Kontsevich-Soibelman \cite{KontsevichSoibelman08}
\cite{KontsevichSoibelman10} and
Nagao \cite{Nagao10}.

\subsection{Reign of the tropics} \label{ss:reign-of-the-tropics}
The following theorem and its corollary are the basis of the `tropicalization
method' which is used in applications of cluster algebras
to discrete dynamical systems\index{discrete dynamical system} and to dilogarithm identities,
cf. \cite{Nakanishi11b} \cite{InoueIyamaKellerKunibaNakanishi10a}
\cite{InoueIyamaKellerKunibaNakanishi10b} \cite{Nakanishi11} \cite{Nakanishi11a}.

\begin{theorem}[\cite{Plamondon11}] \label{thm:rigidity}
Let $\Gamma'$ and $\Gamma''$ be two Ginzburg dg algebras and let
$\Phi': \cd(\Gamma') \to \cd(\Gamma)$ and $\Phi'': \cd(\Gamma'') \to \cd(\Gamma)$
be triangle equivalences satisfying the conditions of Proposition~\ref{prop:width-one-equivalence}.
For $1\leq j\leq n$, let us write $S'_j$ for the image of the $j$th simple module under $\Phi'$ and
$P'_j$ for the image of the module $e_j\Gamma'$. Similarly, we define $S''_j$ and $P''_j$.
Suppose that for each $1\leq j\leq n$, we have $[S'_j]=[S''_j]$ in $K_0(\cd_{fd}(\Gamma))$.
Then for each $1\leq j\leq n$, we have
$
P'_j \cong P''_j \mbox{ and } S'_j \cong S''_j \ .
$
\end{theorem}

We first notice that for each $1\leq j\leq n$, we have the
equality $[P'_j]=[P''_j]$ in $K_0(\per(\Gamma))$. Indeed, this
follows from the duality between the bases
$[P'_i]$ and $[S'_j]$ in $K_0(\per(\Gamma)]$ and $K_0(\cd_{fd}(\Gamma))$.
Now the first isomorphism follows from the fact, proved in
section~3.1 of \cite{Plamondon11a}, that an object $X$ of
$\pr(\Gamma)$ which is rigid\index{object!rigid}, i.e. $\Hom(X,\Sigma X)=0$,
is determined by its class in $K_0(\per(\Gamma))$. The
objects $S'_i$ are the simples of the abelian subcategory of $\cd(\Gamma)$
formed by the objects $U$ such that $\Hom(P'_j,\Sigma^p U)$
vanishes for all $p\neq 0$ and all $1\leq j\leq n$. Among
these simples, $S'_i$ is the only one which receives
a non zero morphism from $P'_i$. Now it is clear
that the isomorphisms for the $P'_j$ imply those
for the $S'_j$.

\begin{corollary} \label{cor:tropicalization}
Let $t$ and $t'$ be vertices of the
$n$-regular tree such that there is a permutation $\pi$
of $\{1, \ldots, n\}$ with $C(t')= P_\pi C(t)$, where
$P_\pi$ is the permutation matrix associated with $\pi$.
Then we have $G(t) = P_\pi G(t')$, the permutation $\pi$ yields
an isomorphism $Q(t) \iso Q(t')$ and for each $1\leq j\leq n$, we have
\begin{itemize}
\item[a)] $T_j(t') = T_{\pi(j)}(t)$ and $S_j(t')=S_{\pi(j)}(t)$;
\item[b)] $F_j(t') = F_{\pi(j)}(t)$;
\item[c)] $x_j(t') = x_{\pi(j)}(t)$ and $y_j(t')=y_{\pi(j)}(t)$.
\end{itemize}
In particular, the seeds associated with $t$ and $t'$ are
isomorphic via $\pi$.
\end{corollary}

To prove the corollary, we apply the theorem to the
equivalences $\Phi(t)$ and $\Phi(t')$. We immediately obtain
part a). This implies the statement on the $g$-matrices and
the quivers. By Theorem~\ref{thm:decategorification}, it
also implies parts b) and c).

A proof of the corollary based on
the study of stability conditions on $\cd_{fd}(\Gamma)$ is
given in section~4.2 of \cite{Nagao10}. It can also
be deduced from the results of \cite{DerksenWeymanZelevinsky10}.

\subsection{Rigid objects and cluster monomials} \label{ss:rigid-objects-and-cluster-monomials}
Let $\tilde{Q}$ be an ice
quiver (equally valued). Let $\Gamma$ be the Ginzburg dg algebra
associated with the opposite of $(\tilde{Q}, W)$, where $W$ is
a generic potential.

For each object $L$ of $\pr(\Sigma^{-1}\Gamma)$ such
that $H^1(L)$ is finite-dimensional, we define a Laurent polynomial
\[
X(L)= \sum_{e} \chi(\Gr_e(H^1_\ca(L)))\, \hat{y}^{\,e} \ko
\]
where $\hat{y}_l = \prod_{i=1}^m x_i^{\tilde{b}_{il}}$, $1 \leq l \leq n$, and
$\hat{y}^{\,e} = \prod_l \hat{y}_l^{\,e_l}$.
By part c) of Theorem~\ref{thm:decategorification} and by the
separation formula of Theorem~\ref{thm:separation-formula}, when
$L=T_i(t)$ for some $1\leq i\leq n$ and some vertex $t$ of the $n$-regular
tree, then $X(L)$ equals the cluster variable $x_i(t)$. It is not
hard to check that for two objects $L$ and $L'$ of $\pr(\Gamma)$,
we have
\[
X(L\oplus L') = X(L) X(L').
\]
So if we apply the map $L\mapsto X(L)$ to direct sums of objects
$T_i(t)$, $1\leq i\leq n$, for a fixed vertex $t$, we recover the
cluster monomials associated with $t$.

Let us call a rigid object $L$ of $\pr(\Sigma^{-1}\Gamma)$ {\em reachable}\index{object!reachable rigid} if
it there is a vertex $t$ of the $n$-regular tree such that $L$ is
a direct sum of copies of the objects $T_i(t)$, $1\leq i\leq n$.

\begin{theorem}[\cite{Plamondon11a} \cite{CerulliKellerLabardiniPlamondon12}]
\begin{itemize}
\item[a)] If $L_1$, \ldots, $L_N$ are pairwise non-isomorphic
reachable rigid objects, then the Laurent polynomials $X(L_1)$, \ldots, $X(L_N)$
are linearly independent.
\item[b)] The map $L \mapsto X(L)$ induces a bijection from the
set of isomorphism classes of reachable rigid objects onto the
set of cluster monomials.
\end{itemize}
\end{theorem}

The surjectivity in b) is proved by Plamondon \cite{Plamondon11a}.
The linear independence in a), and hence the injectivity in b),
is proved in \cite{CerulliKellerLabardiniPlamondon12}.

\subsection{Proof of decategorification} \label{ss:proof-of-decategorification}
We will sketch a proof of Theorem~\ref{thm:decategorification}. We prove a) and b)
simultaneously by induction on the distance of $t$ from $t_0$. For $t=t_0$, there is
nothing to prove. Now assume we have proved the claim for $t$ and that $t'$
is at greater distance from $t_0$ and linked to $t$ by an edge labeled $k$.
Then the coefficients of the $c$-matrix at $t'$ can be computed as
\[
c_{ij}(t') = \left\{ \begin{array}{ll} -c_{ij}(t) & \mbox{ if $j=k$;} \\
c_{ij}(t) + c_{ik}(t) [\eps b_{kj}(t)]_+ + [-\eps c_{ik}(t)]_+ b_{kj}(t) & \mbox{ else}\ko
\end{array} \right.
\]
where $1\leq i,j\leq n$ and $\eps$ is any sign, cf. Prop.~5.8 of \cite{FominZelevinsky02} and formula~(3.3)
in \cite{Nakanishi11a}. We know that $b_{kj}(t)$ equals the number of arrows
from $k$ to $j$ in $Q(t)$ minus the number of arrows from $j$ to $k$ in $Q(t)$.
Thus, we have
\[
b_{kj}(t) = \dim \Ext^1_{\Gamma}(S_k(t),S_j(t)) - \dim \Ext^1_{\Gamma}(S_j(t),S_k(t)).
\]
By the induction hypothesis, the coordinates of $[S_k(t)]$ in the basis of the
$[S_i]$ are the components $c_{ik}(t)$, $1\leq i\leq n$, of the $c$-vector $C(t) e_k$.
By Lemma~\ref{lemma:dichotomy-for-simples}, they are all of the same sign. Let us
choose $\eps$ equal to this sign. Then the formula for the $c_{ij}(t')$ simplifies
as follows:
\[
c_{ij}(t') = \left\{ \begin{array}{ll} -c_{ij}(t) & \mbox{ if $j=k$;} \\
c_{ij}(t) + c_{ik}(t) [\eps b_{kj}(t)]_+  & \mbox{ else.} \end{array} \right.
\]
Now assume that $\eps=1$. This means that $S_k(t)$ lies in $\ca$ and that
$S_k(t')$ is $\Sigma^{-1}S_k(t)$. Let us put $m=b_{kj}(t)$. If we have $m\leq 0$, then the space
$\Ext^1_\Gamma(S_k(t), S_j(t))$ vanishes and we have $S_j(t')=S_j(t)$. If we have
$m\geq 0$, then we get $m=\Ext^1_\Gamma(S_k(t), S_j(t))$
and the object $S_j(t')$ is constructed as a universal extension:
\[
(\Sigma^{-1} S_k(t))^m \to S_j(t) \to S_j(t') \to S_k(t)^m.
\]
In both cases, the formula for $c_{ij}(t')$ gives the correct
multiplicity of $[S_i]$ in $[S_j(t')]$. Now suppose that $\eps=-1$.
Then $S_k(t)$ belongs to $\Sigma^{-1}\ca(t)$ and $S_k(t')$ is $\Sigma S_k(t)$.
Let us put $m=-b_{kj}(t)=b_{jk}(t)$.
If we have $m\leq 0$, then the space $\Ext^1_{\Gamma}(S_j(t), S_k(t))$ vanishes
and $S_j(t')=S_j(t)$.  If we have $m\geq 0$, then
we get $m=\dim\Ext^1_{\Gamma}(S_j(t), S_k(t))$ and $S_j(t')$ is constructed
as a universal extension
\[
S_k(t)^m \to S_j(t') \to S_j(t) \to \Sigma S_k(t)^m.
\]
Again, in both cases, the formula for $c_{ij}(t')$ gives the correct
multiplicity of $[S_i]$ in $[S_j(t')]$.

We get part c) as a consequence: Indeed, by part b) the main
conjecture~\ref{conj:main} holds for $Q$ and so we have
$G(t)^T C(t)=I$ for all vertices $t$ of the $n$-regular
tree, by the tropical duality theorem~\ref{thm:tropical-duality}.
On the other hand, the basis of the $[P_i(t)]$ is dual
to that of the $[S_j(t)]$. Clearly, this implies c).

We deduce d) from Plamondon's results \cite{Plamondon11}.
Indeed, both $T_j(t)$ and $\Sigma^{-1}\Gamma$ belong to $\pr(\Sigma^{-1}\Gamma)$.
Thus, by Proposition~\ref{prop:equiv-presentable}, we have
\[
H^1_{\ca}(T_j(t)) = \Hom_{\per(\Gamma)}(\Sigma^{-1}\Gamma, T_j(t)) \iso
\Hom_{\cc_\Gamma}(\pi(\Sigma^{-1} \Gamma), \pi(T_j(t))) \ko
\]
where $\cc_\Gamma = \per(\Gamma)/\cd_{fd}(\Gamma)$ is the cluster
category and $\pi$ the projection functor. Let us omit this
functor from the notations. Since $T_j(t)\in \cc_\Gamma$ is
obtained by iterated mutation from $\Gamma$, it belongs
to Plamondon's category $\cd\subset \cc_\Gamma$ and so the
space
\[
\Hom_{\cc_\Gamma}(T_j(t), \Sigma \Gamma)
\]
is finite-dimensional. By Prop.~2.16 of \cite{Plamondon11}, this
space is in duality with
\[
\Hom_{\cc_\Gamma}(\Sigma^{-1} \Gamma, T_j(t))
\]
which therefore also finite-dimensional. So we find that $H^1_{\ca}(T_j(t))$
is finite-dimen\-sional and in duality with $\Hom_{\cc_\Gamma}(T_j(t), \Sigma \Gamma)$.
Now let $P \mapsto P^\vee$ denote the canonical equivalence
\[
\per(\Gamma)^{op} \iso \per(\Gamma^{op}) \ko P \mapsto P^\vee =\RHom_\Gamma(P,\Gamma).
\]
It induces an equivalence $\cc^{op}_\Gamma \iso \cc_{\Gamma^{op}}$ still
denoted by the same symbol. We have
\[
\Hom_{\cc_\Gamma}(T_j(t), \Sigma \Gamma) \iso \Hom_{\cc_{\Gamma^{op}}}(\Sigma \Gamma)^\vee, T_j(t)^\vee)
= \Hom_{\cc_{\Gamma^{op}}}(\Sigma^{-1} \Gamma, T_j(t)^\vee).
\]
Notice that $\Gamma^{op}= \Gamma(Q^{op},W^{op})^{op})=\Gamma(Q,W)$.
So we get that the left $J(Q,W)$-module $H^1_{\ca}(T_j(t))$ is in duality
with the right $J(Q,W)$-module
\[
\Hom_{\cc_{\Gamma^{op}}}((\Sigma \Gamma)^\vee, T_j(t)^\vee) = \Hom_{\cc_{\Gamma^{op}}}(\Sigma^{-1} \Gamma, T'_j(t)) \ko
\]
where $T'_j(t)$ denotes the object obtained from $\Gamma$ in $\cc_{\Gamma^{op}}$
by the sequence of mutations linking $t_0$ to $t$.
Thus, the Grassmannian of $e$-dimensional quotients of $H^1_{\ca}(T_j(t))$
identifies with the Grassmannian of $e$-dimensional submodules of the
above $J(Q,W)$-module. The generating series of their Euler characteristics
is the $F$-polynomial associated with $T'_j(t)$ in Def.~3.14 of
\cite{Plamondon11a} and it equals the $F$-polynomial $F_j(t)$, as shown
in section~4.2 of \cite{Plamondon11a}.

\subsection{Proof of the quantum dilogarithm identities} \label{ss:proof-quantum-dilog}
We will sketch a proof of Theorem~\ref{thm:quantum-dilog-identity}.
We start with part a). We prove the stronger statement given in
Remark~\ref{rem:combinatorial-DT-invariant}. So suppose that,
in the notations of the remark, we have $P C(t) = C(t')$ for
the permutation matrix $P=P_\pi$ associated with a permutation
$\pi$ of $\{1,\ldots, n\}$.  By Corollary~\ref{cor:tropicalization}, we find that the
seeds associated with $t$  and $t'$ are isomorphic via $\pi$ in any cluster algebra associated with a matrix $\tilde{B}$ whose
principal part $B$ corresponds to $Q$.
Now by Theorem~\ref{thm:quantum-exchange-graph}, 
we find that the quantum seeds associated with $t$ and $t'$ 
are isomorphic via $\pi$ in any quantum cluster algebra
associated with a compatible pair $(\tilde{B},\Lambda)$,
where the principal part of $\tilde{B}$ is the given matrix $B$.
Thus, in the notations of Theorem~\ref{thm:quantum-separation},
we have $\Phi(\mathbf{i})=\Phi(\mathbf{i}')$. Now by the duality 
theorem~\ref{thm:tropical-duality}, we also have $P_\pi G(t)=G(t')$ 
for the same permutation $\pi$. By the equality~(\ref{eq:separation-many}) in
Theorem~\ref{thm:quantum-separation}, we obtain
\[
\Ad'(\E(\mathbf{i}))=\Ad'(\E(\mathbf{i'})).
\]
Now let us choose $\tilde{B}=B_{pr}$. We find that the 
power series $\E(\mathbf{i})\E(\mathbf{i}')^{-1}$ in the
variables $x_1$, \ldots\ , $x_n$ commutes with
the variables $x_{n+i}$, $1\leq i\leq n$.
Now by our choice of $\tilde{B}=B_{pr}$, we have
\[
x_{n+i} x_j = q^{\delta_{ij}} x_j x_{n+i}
\]
for all $1\leq i,j \leq n$.
This implies that for any power series $f(x_1, \ldots, x_n)$,
we have
\[
x_{n+i} f(x_1, \ldots, x_n) x_{n+i}^{-1} = f(x_1, \ldots,q x_i, \ldots, x_n).
\]
So a power series in $x_1$, \ldots, $x_n$ which commutes with
all the $x_{n+i}$, $1\leq i\leq n$, must be constant. Since the
constant term of $\E(\mathbf{i})\E(\mathbf{i}')^{-1}$ is $1$, 
we find $\E(\mathbf{i})\E(\mathbf{i}')^{-1}=1$
as claimed.

For part b), we have to invoke Donaldson-Thomas theory\index{Donaldson-Thomas theory} in its form pioneered by Kontsevich-Soibelman \cite{KontsevichSoibelman08}
\cite{KontsevichSoibelman10}. This theory is not yet completely
developed for formal potentials\protect{\footnote{cf. section~4.8 in
\cite{Efimov11} and the discussion in section~7.1
of \cite{KontsevichSoibelman10}}} and therefore, for the moment,
does not apply to arbitrary quivers (cf.  \cite{BehrendBryanSzendroi09}
\cite{Mozgovoy11} \cite{Nagao11a}
for recent progress on special classes). However, in its final
form, the theory should yield the following:
Let $\aff_{Q}$ denote the completed
quantum affine space associated
with $Q$. Let $\ca$ be the category of finite-dimensional (hence nilpotent)
right modules over the completed Jacobian algebra $J(Q,W)$ of
the quiver $Q$ endowed with a generic potential. Let $\ct_1$ and
$\ct_2$ be torsion subcategories of $\ca$. Following the
explanation after Remark~21 on page~90 of \cite{KontsevichSoibelman08}
we define $\ct_1$ to be {\em constructibly less than or equal to} $\ct_2$
if $\ct_1$ is contained in $\ct_2$ and for each dimension vector
$d$, the subset of the variety of (contravariant) representations of $J(Q,W)$
with dimension vector $d$ formed by the points corresponding
to modules in $\ct_1^\perp\cap \ct_2$ is constructible.
In this case, following \cite{KontsevichSoibelman08} we write
\begin{equation} \label{eq:constructibly-less}
\ct_1 \leq_{constr} \ct_2.
\end{equation}
What the fully fledged version of Kontsevich-Soibelman's theory
should yield is a {\em DT-character on $\ca$}, i.e. the datum 
of an element
$A_{\ct_1,\ct_2}$ of the group $\aff_Q^{\times}$ for each pair
of torsion theories $\ct_1$, $\ct_2$ satisfying (\ref{eq:constructibly-less}) such
that the following hold
\begin{itemize}
\item[a)] whenever we have three torsion theories $\ct_1$, $\ct_2$ and $\ct_3$
such that
\[
\ct_1 \leq_{constr} \ct_2\ko \ct_2 \leq_{constr} \ct_3 \mbox{ and } \ct_1 \leq_{constr} \ct_3 \ko
\]
we have
\begin{equation} \label{eq:factorization-property}
A_{\ct_1,\ct_2} A_{\ct_2, \ct_3} = A_{\ct_1, \ct_3}  \; ;
\end{equation}
\item[b)] if we have $\ct_2 = \ct_1 \star \add(L)$, where $L$ is
a module in $\ct_1^\perp$ satisfying $\Hom(L,L)=\C$ and
$\Ext^1(L,L)=0$, we have
\begin{equation} \label{eq:single-factor}
A_{\ct_1, \ct_2} = \E(x^\alpha) \ko
\end{equation}
where $\alpha$ is the dimension vector of $L$.
\end{itemize} 
The {\em non commutative DT invariant} associated with $Q$ and
the given DT-character is then the power series
\begin{equation} \label{eq:non-com-DT-invariant}
DT_{Q}=A_{0,\ca} \in \aff_{Q}.
\end{equation}
Via the duality functor $\Hom_\C(?,\C)$ and the canonical
isomorphism $\aff_{Q}^{op} \iso \aff_{Q^{op}}$
taking $x^\alpha$ to $x^\alpha$, a DT-character for $Q$
yields one for $Q^{op}$ and $DT_Q \in \aff_Q$ is
mapped to $DT_{Q^{op}} \in \aff_{Q^{op}}$.

Now assume that we have a quiver $Q$ whose {\em non commutative
DT-invariant is defined}, i.e. it admits a DT-character. Then this
also holds for $Q^{op}$.
Suppose that we are in the situation of part b) of
Theorem~\ref{thm:quantum-dilog-identity} so that $-C(t_N)$
is a permutation matrix. Then by Theorem~\ref{thm:rigidity},
the simples of $\ca(t_N)$ lie in $\Sigma^{-1}\ca$ and so we must have
$\ca(t_N)=\Sigma^{-1}\ca$ and $\ct(t_N)=\ca$. Now the
torsion subcategories
\[
\{0\}=\ct(t_0) \ko \ct(t_1) \ko \ldots \ko \ct(t_N)=\ca
\]
form a sequence such that for each $1\leq s \leq N$, we either
have
\begin{itemize}
\item[(1)] $S_{i_s}(t_{s-1}) \in \ca$ and then
\[
\ct(t_{s-1}) \leq_{constr} \ct(t_s) \mbox{ and } \ct(t_s)=\ct(t_{s-1})\star \add(S_{i_s}(t_{s-1}))
\]
or
\item[(2)] $S_{i_s}(t_{s-1}) \in \Sigma^{-1}\ca$ and then
\[
\ct(t_{s}) \leq_{constr} \ct(t_{s-1}) \mbox{ and } \ct(t_{s-1})=\ct(t_{s})\star \add(\Sigma S_{i_s}(t_{s-1})).
\]
\end{itemize}
depending on the sign of the $c$-vector $\beta_s=C(t_{s-1})e_{i_s}$, which is
just the (signed) dimension vector of $S_{i_s}(t_{s-1})$, by
Theorem~\ref{thm:decategorification}. By induction on $s$, one now
proves that
\[
\E(\eps_1 \beta_1)^{\eps_1}  \ldots\, \E(\eps_s \beta_s)^{\eps_s} = A_{0,\ct(t_s)}.
\]
For $s=N$, we obtain the equality
\[
\E(\eps_1 \beta_1)^{\eps_1}  \ldots\, \E(\eps_s \beta_N)^{\eps_N} = A_{0,\ca} = DT_{Q^{op}}
\]
in $\aff_{Q^{op}}$. Its image under the canonical isomorphism
$\aff_{Q}^{op} \iso \aff_{Q^{op}}$ is the claimed equality in $\aff_Q$.

\def\cprime{$'$} \def\cprime{$'$}
\providecommand{\bysame}{\leavevmode\hbox to3em{\hrulefill}\thinspace}
\providecommand{\MR}{\relax\ifhmode\unskip\space\fi MR }
\providecommand{\MRhref}[2]{%
  \href{http://www.ams.org/mathscinet-getitem?mr=#1}{#2}
}
\providecommand{\href}[2]{#2}

%\frenchspacing
%\bibliographystyle{amsplain}
%\bibliography{stanKeller}

%\printindex

 \end{document}